\documentclass[12pt]{article}
\usepackage{amsmath}
\usepackage{amssymb}
\usepackage{amsthm}
\usepackage[totalwidth=17cm, totalheight=25cm]{geometry}
\usepackage{graphicx}
\usepackage{mathabx}
\usepackage{tabularx}
	\newcolumntype{C}[1]{>{\centering\arraybackslash}m{#1}} 
	\newcolumntype{R}[1]{>{\raggedleft\arraybackslash}m{#1}} 

\newtheoremstyle{boldplain}
{9pt}
{9pt}
{\itshape}
{}
{\bfseries}
{.}
{.5em}
{\thmname{#1}\thmnumber{ #2}\thmnote{ (#3)}}%

\newtheoremstyle{bolddefinition}
{9pt}
{9pt}
{}
{}
{\bfseries}
{.}
{.5em}
{\thmname{#1}\thmnumber{ #2}\thmnote{ (#3)}}%

\theoremstyle{boldplain}
\newtheorem{add}[equation]{Addendum}

\newtheorem{cor}[equation]{Corollary}

\newtheorem{lem}[equation]{Lemma}

\newtheorem{prop}[equation]{Proposition}

\newtheorem{thm}[equation]{Theorem}
\newtheorem{theorem}[equation]{Theorem}

\newtheorem{fact}[equation]{Fact}

\theoremstyle{bolddefinition}
\newtheorem{dfn}[equation]{Definition}

\newtheorem{defn}[equation]{Definition}

\newtheorem{ques}[equation]{Question}

\newtheorem{rem}[equation]{Remark}

\newfont{\bigbf}{cmbx10 scaled\magstep1}
\normalbaselines
\numberwithin{equation}{section}

\def\no{\noindent}

\def\R{{\mathbb R}}

\def\N{{\mathbb N}}

\def\al{\alpha}

\def\Ga{\Gamma}
\def\de{\delta}
\def\De{\Delta}
\def\eps{\epsilon}
\def\la{\lambda}
\def\La{\Lambda}
\def\si{\sigma}
\def\Si{\Sigma}

\def\om{\omega}

\def\3{\ss}

\def\acts{\curvearrowright}
\def\amod{a_{mod}}

\def\B{\operatorname{B}}

\def\D{\partial}
\def\dH{\mathop{\hbox{dist}_{Haus}}}

\def\diam{\mathop{\hbox{diam}}}
\def\diamo{\diamondsuit}
\def\diamon{\diamondsuit_n}

\def\diamoom{\diamondsuit_{\om}}

\def\diamot{\diamondsuit_{\tau_{mod}}}
\def\diamoTh{\diamondsuit_{\Theta}}
\def\diamotTh{\diamoTh}
\def\Dt{\D_{\taumod}}
\def\DtX{\D_{\taumod}X}
\newcommand {\dist}{\mathrm{dist}}

\def\Flag{\operatorname{Flag}}

\def\Flagt{\Flag_{\tau_{mod}}}
\def\Fmod{F_{mod}}

\def\geo{\partial_{\infty}}

\def\half{\frac{1}{2}}

\def\id{\mathop{\hbox{id}}}

\def\inte{\operatorname{int}}

\def\Isom{\operatorname{Isom}}
\def\Isomth{\operatorname{Isom}_{\theta}}

\def\lra{\longrightarrow}

\def\oa{\overrightarrow}

\def\pihalf{\frac{\pi}{2}}

\def\2pithird{\frac{2\pi}{3}}

\def\rad{\operatorname{rad}}

\def\rank{\mathop{\hbox{rank}}}

\def\pO{\operatorname{pO}}
\def\pOt{\operatorname{pO}^{\taumod}}
\def\bO{\operatorname{bO}}
\def\bOt{\operatorname{bO}^{\taumod}}
\def\tiAt{\operatorname{\tilde A}^{\taumod}}
\def\barpO{\operatorname{p\bar O}}
\def\barpOt{\operatorname{p\bar O}^{\taumod}}
\def\barbO{\operatorname{b\bar O}}
\def\barbOt{\operatorname{b\bar O}^{\taumod}}
\def\barXt{\operatorname{\bar X}^{\taumod}}
\def\pSh{\operatorname{pSh}}
\def\pSht{\operatorname{pSh}^{\taumod}}
\def\bSh{\operatorname{bSh}}
\def\bSht{\operatorname{bSh}^{\taumod}}
\def\simod{\si_{mod}}
\def\st{\operatorname{st}}
\def\stTh{\operatorname{st}_{\Theta}}
\def\stThm{\operatorname{st}_{\Theta_-}}
\def\stThp{\operatorname{st}_{\Theta_+}}
\def\ost{\operatorname{ost}}

\def\tangle{\angle_{Tits}}
\def\taumod{\tau_{mod}}

\def\tits{\partial_{Tits}}

\def\ulim{\mathop{\hbox{$\om$-lim}}}

\def\Vmod{V_{mod}}

\def\8{\infty}
\def\<{\langle}
\def\>{\rangle}

\title{A Morse Lemma for quasigeodesics in symmetric spaces and euclidean buildings}
\author{Michael Kapovich, Bernhard Leeb, Joan Porti}
\date{July 4, 2015}

\begin{document}

\maketitle

\begin{abstract}
\no We prove a Morse Lemma for coarsely regular quasigeodesics 
in nonpositively curved symmetric spaces and euclidean buildings $X$.
The main application is a simpler
coarse geometric characterization of 
{\em Morse subgroups} of the isometry groups $\Isom(X)$
as undistorted subgroups which are coarsely uniformly regular.
We show furthermore that they must be word hyperbolic. 
We introduced this class of discrete subgroups in our earlier paper \cite{morse}, 
in the context of symmetric spaces, 
where various equivalent geometric and dynamical characterizations of 
word hyperbolic Morse subgroups were established,
including the Anosov subgroup property.
This version of the paper is essentially identical to the first version \cite{mlem}
except for Remarks~\ref{rem:hist} and~\ref{rem:hist2} related to the later preprint \cite{GGKW}.
\end{abstract}

\tableofcontents

\section{Introduction}

One of the important features of $\delta$-hyperbolic geodesic metric spaces is the {\em Morse Lemma}, also known as the {\em Stability of Quasigeodesics}: Every (uniform) quasigeodesic is (uniformly) close to a geodesic. This property of hyperbolic spaces is used, among other things, to show that hyperbolicity is a quasiisometry invariant and that quasiisometries between hyperbolic spaces extend to the ideal boundaries. Stability of quasigeodesics is known to fail in CAT(0) metric spaces: Already euclidean plane contains quasigeodesics which are not Hausdorff-close to any geodesic. Some versions of the Morse lemma are known for CAT(0) spaces: In the case of maximal quasiflats \cite{LS} and in the case of {\em Morse quasigeodesics} (also known as {\em hyperbolic} or {\em rank one} quasigeodesics), see \cite{Sultan}. Nevertheless, the real understanding of what should constitute a true analogue of the Morse lemma in the CAT(0) setting, remains elusive. 

The main goal of this paper is to prove an analogue of the Morse Lemma for {\em regular} quasigeodesics in nonpositively curved symmetric spaces and euclidean buildings. In order to unify the terminology, we refer to nonpositively curved symmetric spaces and euclidean buildings as {\em model spaces} throughout the paper. 
 
Instead of concluding that regular quasigeodesics are uniformly close to geodesics 
(which is far from true, see section \ref{sec:examples}), we will prove that they are contained in uniform neighborhoods of certain convex subsets of the model space: {\em Diamonds} in the case of finite quasigeodesics, {\em Weyl cones} in the case of quasigeodesic rays and {\em parallel sets} (more precisely, unions of opposite Weyl cones therein) in the case of complete quasigeodesics.  

Note that the question about regular quasigeodesics reduces to the case of model spaces without flat factors. For, if a model space has a nontrivial flat de Rham factor, 
then all diamonds and Weyl cones split off this flat factor,
and the canonical projection to the complementary (model space) factor 
preserves $\taumod$-regularity of segments and paths.
We therefore restrict our discussion to model spaces without flat factors.

Our main motivation for these results comes from the theory of discrete isometric group actions on model spaces $X$,  more specifically, the desire to give a clean coarse-geometric  characterization of {\em Morse actions}, which have been introduced in \cite{morse}. 

The notion of {\em regularity} used in our paper is defined relative to a certain face $\taumod$ of the model chamber $\simod$ of the Tits boundary of $X$. The definition is the easiest in the case when $\taumod=\simod$ and we first present our results in this setting. 

A quasigeodesic $q$ in $X$ (which might be finite or infinite) 
is {\em coarsely uniformly regular} if any two points $x, y$ in $q$ which are sufficiently far apart ($d(x,y)\ge D$), define a {\em uniformly regular segment} $xy$ in $X$, i.e., a geodesic segment whose direction belongs to a fixed compact subset $\Theta$ of the interior of $\simod$. 
A {\em diamond} $\diamo(x,y)$ in $X$ is a generalization of a geodesic segment.
In the case when the segment $xy$ is regular, $\diamo(x,y)$ is a (convex) subset of a flat $F\subset X$ containing $xy$, 
namely the intersection of two Weyl chambers 
$V(x,\si)\cap V(y, \hat\si)$ with the tips at $x$ and $y$ respectively, over opposite chambers $\si,\hat\si$ in the Tits boundary $\tits X$ of $X$.  

\begin{thm}[Morse Lemma, regular case]
\label{main-sigma}
(i) Every finite coarsely uniformly regular quasigeodesic path $q$ in $X$ with endpoints $x$ and $y$
of distance $d(x,y)\geq D$
is contained in a neighborhood of the diamond $\diamo(x,y)$. 

(ii) Every coarsely uniformly regular quasigeodesic ray $q$ in $X$ with initial point $x$
is contained in a neighborhood of a unique euclidean Weyl chamber $V=V(x,\si)$.

(iii) Every coarsely uniformly regular complete quasigeodesic $q$ in $X$ 
is contained in a neighborhood of a unique maximal flat $F$ and, moreover, is contained in a neighborhood of 
the union 
$V(z,\si)\cup V(z, \hat\si)\subset F$
of two opposite euclidean Weyl chambers with common tip $z\in F$.

Furthermore, the distance from  $q$ to $\diamo(x,y)$, $V$ respectively $V(z,\si)\cup V(z, \hat\si)$
is bounded above in terms of the quasiisometry constants of $q$,
the coarseness scale $D$ and the regularity set $\Theta$.  
In case (iii), the common tip $z$ can be chosen uniformly close to any point on $q$.
\end{thm}
In other words, 
each coarsely uniformly regular quasigeodesic in $X$ 
is a {\em Morse quasigeodesic} in the sense of \cite{morse}. 

We deduce from this result that 
{\em coarsely uniformly regularly quasiisometrically embedded subspaces} 
in model spaces $X$ must be Gromov-hyperbolic. 
A quasiisometric embedding $f$ from a geodesic metric space $Y$ into $X$ 
is {\em coarsely uniformly regular} 
if the images of geodesic segments in $Y$ are 
{\em coarsely uniformly regular} in $X$. 
(Uniformity here refers to the constant $D$ and the subset $\Theta$.)
 
\begin{thm}[Hyperbolicity of domain and boundary map, regular case]
If $f: Y\to X$ is a coarsely uniformly regular quasiisometric embedding, 
then the space $Y$ is Gromov-hyperbolic and,
if it is also locally compact, the map $f$ extends to a topological embedding 
from the Gromov boundary of $Y$ into the F\"urstenberg boundary of $X$.  
\end{thm}
Our work is primarily motivated by the study of discrete subgroups 
of isometry groups of nonpositively curved symmetric spaces; before discussing these, we explain how the theorems above generalize to 
coarsely $\taumod$-regular quasigeodesics 
and coarsely $\taumod$-regular quasiisometric embeddings. 

\medskip 
{\bf $\taumod$-regularity.} The role of the compact $\Theta\subset \inte(\simod)$ which appeared above, 
is now played by a ``Weyl-convex'' compact subset $\Theta\subset \simod$ which intersects the boundary of $\simod$ only in the open faces containing the open simplex $\inte(\taumod)$. 
For instance, if $\taumod$ is a vertex, then $\Theta$ is required to be disjoint from the top-dimensional face of $\simod$ not containing $\taumod$. With this modification, the  definition of coarsely uniformly regular quasigeodesics 
generalizes to the one of coarsely uniformly $\taumod$-regular quasigeodesics. 

Let $\tau$ be a simplex of the Tits boundary of $X$ which has {\em type} $\taumod$. 
We next describe the replacement for the euclidean Weyl chambers $V(x, \si)$ in $X$ with $\tau$ playing the role of $\si$. 
They are replaced by the {\em Weyl cones} $V(x, \st(\tau))$: These convex subsets of $X$ are unions of geodesic rays $x\xi$  in $X$ connecting $x$ to ideal boundary points $\xi\in \geo X$, which belong to a certain  subcomplex $\st(\tau)\subset \geo X$. 
This subcomplex is the union of chambers $\si$ containing $\tau$. The cones $V(x, \st(\tau))$ are no longer contained in maximal flats in $X$ (unless $\taumod=\simod$); instead, each 
$V(x, \st(\tau))$ is a subset of the parallel set in $X$ of a geodesic $l$ through $x$, 
asymptotic to a generic point in $\tau$. Such parallel sets are said to have the {\em type $\taumod$}. 
The diamonds $\diamo_{\taumod}(x,y)$ are again defined as intersections
$$
V(x, \st(\tau)) \cap V(y, \st(\hat\tau))
$$
for opposite simplices $\tau, \hat\tau$ in the Tits boundary of $X$. 

Now, we are ready to state our results. 
The main result is a Morse Lemma for coarsely $\taumod$-regular quasigeodesics in model spaces
(see Theorem~\ref{thm:mlemt} and Corollary~\ref{cor:infqgclcp}):
\begin{thm}[Morse Lemma]
\label{main-tau}
(i) Every finite coarsely uniformly $\taumod$-regular quasigeodesic path $q$ in $X$ with endpoints $x$ and $y$
is contained in a neighborhood of the diamond $\diamot(x,y)$. 

(ii) Every coarsely uniformly $\taumod$-regular quasigeodesic ray $q$ in $X$ with initial point $x$
is contained in a neighborhood of a unique Weyl cone $V=V(x,\st(\tau))$ of type $\taumod$.

(iii) Every coarsely uniformly $\taumod$-regular complete quasigeodesic $q$ in $X$ 
is contained in a neighborhood of a unique parallel set $P$ of type $\taumod$ and, moreover, is contained in a neighborhood of the union 
$V(z,\st(\tau))\cup V(z, \st(\hat\tau))\subset P$ 
of opposite Weyl cones of type $\taumod$ with common tip $z\in P$.

Furthermore, the distance from  $q$ to $\diamot(x,y)$, $V$ respectively $V(z,\st(\tau))\cup V(z, \st(\hat\tau))$
is bounded above in terms of the quasiisometry constants of $q$,
the scale $D$ and the subset $\Theta$.  
In case (iii), the common tip $P$ can be chosen uniformly close to any point on $q$.
\end{thm}

In order to help the reader to appreciate the relation of this theorem to Theorem \ref{main-sigma}, we note that the regularity assumptions in Theorem \ref{main-tau} are weaker (directions of segments $xy$ are allowed to belong to larger subsets of $\simod$), while the conclusions are weaker as well, since we can only conclude that quasigeodesics lie close to certain sets which are larger than the ones in Theorem \ref{main-sigma}. 

\medskip 
Applying these results about coarsely regular quasigeodesics 
to quasi-isometric embeddings we obtain
(see Theorems~\ref{thm:regQIembeddings} and~\ref{thm:bdcoaregmap}):

\begin{thm}[Hyperbolicity of domain and boundary map]
\label{thm:main1}
Suppose that $q: Z\to X$ is a coarsely uniformly $\taumod$-regular
quasiisometric embedding from a quasigeodesic metric space into a model space. Then:

(i) $Z$ is Gromov hyperbolic. 

(ii) If $Z$ is locally compact, the map $q$ extends to a map 
$$\bar q: \bar Z \to \barXt$$
from the (visual) Gromov compactification $\bar Z=Z\cup\geo Z$,
which is continuous at $\geo Z$ 
and sends distinct ideal boundary points to
antipodal elements of the flag space $\DtX=\Flagt(\geo X)$. 
\end{thm}

The main application of this theorem is a simpler
coarse-geometric characterization of 
{\em Morse subgroups} of the isometry groups 
$G=\Isom(X)$ of model spaces $X$. 
This class of discrete subgroups of $G$ was defined in \cite{morse} (in the context of symmetric spaces), 
where various equivalent characterizations of 
word hyperbolic Morse subgroups were established 
(including the {\em Anosov subgroups} of Labourie, Guichard and Wienhard).  
We obtain (see Corollary~\ref{cor:regimpmoract} and Theorem~\ref{thm:frhypgp}):

\begin{thm}\label{thm:morse subgroups}
The following are equivalent for a finitely generated group $\Ga$ and a 
homomorphism $\Ga\to G$:  

1. The group $\Ga$ is hyperbolic and the homomorphism $\Ga\to G$ is $\taumod$-Morse.

2. The orbit maps $\Ga \to X$ are coarsely $\taumod$-regular quasiisometric embeddings. 
\end{thm}

The following corollary is a higher rank analogue of one of the standard characterizations of convex-cocompact subgroups of rank 1 Lie groups as finitely-generated undistorted subgroups. The regularity condition in this corollary is necessary already for subgroups of $SL(2,\R)\times SL(2,\R)$, see Example 6.34 in \cite{morse}. 

\begin{cor}
A finitely generated subgroup $\Ga< G$ is word hyperbolic and $\taumod$-Morse 
if and only if $\Ga$ is undistorted in $G$ and 
asymptotically uniformly $\taumod$-regular.  
\end{cor}

Note that asymptotic uniform $\simod$-regularity 
of a discrete subgroup $\Ga< G$ means that the geometric limit set $\La(\Ga)\subset \geo X$ (the accumulation set of a $\Ga$-orbit in the ideal boundary of $X$) contains no singular points. 

\medskip
{\bf Strategy of the proof.} The main idea behind the proof of Theorems \ref{main-sigma} and \ref{main-tau} is
inspired by trying to follow 
the proof of the Morse Lemma for $\delta$-hyperbolic metric spaces (which goes back to M.~Morse
himself): If a coarsely uniformly $\taumod$-regular quasigeodesic path $q$ connecting points $x, y\in X$ 
strays too far from the
diamond $\diamo=\diamo_{\taumod}(x,y)$, we use the nearest-point projection to $\diamo$ 
to show that it is a uniformly inefficient connection of its endpoints.
This leads to a conflict,
because sufficiently long quasigeodesics have long, arbitrarily efficient
(i.e. almost distance minimizing) subpaths.
In the setting of a
$\delta$-hyperbolic space $X$, 
one is helped by the fact that the nearest-point projection $\pi_\diamo$ to a geodesic
interval $\diamo\subset X$ contracts distances in metric balls $B(z, R)\subset X$ by an exponentially large factor,
in terms of the minimal distance between $B(z,R)$ and $\diamo$. This fails in the setting of higher-rank symmetric
spaces and euclidean buildings $X$. Instead, we define a certain length metric (of Carnot--Finsler type)
$d_{\diamo}$ on diamonds $\diamo\subset X$ by restricting to a certain class of piecewise-geodesic paths in
$\diamo$, which we call {\em non-longitudinal}. The definition of such paths is, again, quite technical, but (if
$\taumod=\simod$) the reader can think of piecewise-geodesic paths where each subsegment is a singular geodesic.

We then prove (Theorem \ref{thm:contrprojdiamo}):
\begin{theorem}
For each euclidean building $X$ (equipped with its standard CAT(0) metric $d$), the projection 
\begin{equation*}
(X,d)\buildrel \pi_{\diamo}\over\lra(\diamo, d_{\diamo})
\end{equation*}
is locally 1-Lipschitz outside $\diamo$. 
\end{theorem}

In order to appreciate the strength of this statement, we note that the
pseudometric $d_\diamo$ is strictly larger than the metric $d$ when
restricted to {\em longitudinal} segments in $\diamo$. 
Therefore the 
above theorem establishes constraints on the behavior of 
rectifiable regular paths in $X$ and, in particular, of regular bilipschitz paths. 
As an application, we prove
(see section~\ref{sec:rectpaths}):

\begin{thm}
Suppose that $c$ is an ``almost length-minimizing'' uniformly $\taumod$-regular path 
in a euclidean building $X$,  connecting points $x$ and $y$. Then $c$ has to meet  the diamond  
$\diamot(x,y)$ in one more point, besides $x$ and $y$. 
\end{thm}

The condition that $c$ is ``almost length-minimizing'' is, actually, not very restrictive, 
since each rectifiable path in $X$ contains such subpaths. As an application of this result we then show
(see Theorem~\ref{thm:rrpathsinflats}):
\begin{thm}\label{thm:bilip-close}
Every rectifiable uniformly $\taumod$-regular path $c$ in $X$ is entirely contained in the diamond 
$\diamot(x,y)$ determined by the endpoints $x, y$ of $c$. 
\end{thm}

This theorem is reminiscent of the fact that each topological path in a real tree $T$ is contained in a geodesic
segment in $T$. This fact is used for proving the Morse Lemma in the case of $\delta$-hyperbolic geodesic spaces
via ultralimits (cf. \cite{DrutuKapovich}). Our argument then proceeds roughly along the same lines as that
proof.
Namely, assuming that Theorem \ref{main-tau} fails, we construct a sequence of coarsely uniformly
$\taumod$-regular $(L,A)$-quasigeodesic paths $q_n$ in a model space $X$; after passing to suitable ultralimits,
we obtain a uniformly $\taumod$-regular bilipschitz path $q_{\om}$ in an asymptotic cone $X_\om$ of $X$, 
which violates Theorem
\ref{thm:bilip-close}.  
(The ultralimit $X_\om$ is a euclidean building.) Along the way, we have to overcome yet
another difficulty. One of the steps in proving the Morse Lemma in the hyperbolic setting via asymptotic cones, is to
show that each asymptotic cone is a uniquely geodesic space (i.e., every geodesic in the cone is the ultralimit
of a sequence of geodesics). Similarly, in our proof, we have to show that the ultralimit of a sequence of
parallel sets in $X$ is a parallel set in the cone $X_\om$ (a priori, it is just a proper subset of such a
parallel set),
and analogous statements for Weyl cones and diamonds.

{\bf Organization of the paper.} In section \ref{sec:prelim} we review basic notions in the theory of model spaces
(nonpositively curved symmetric spaces and buildings), as well as ultralimits and asymptotic cones. Section
\ref{sec:modgeom} is long and technical, it contains the bulk of technical results of the paper. In this section
we define and analyze properties of Weyl cones and diamonds.  
We then define two key notions in the paper:
Regularity and longitudinality of broken segments and paths in model spaces, as well as their coarse analogues. We
establish a preliminary analogue of 
Theorem \ref{main-tau} for a certain class of broken paths in euclidean
buildings, called {\em straight paths}. Furthermore, we prove results about ultralimits of parallel sets, Weyl
cones and diamonds.  
In section \ref{sec:Carnot-Finsler} we define the modified metric $d_\diamo$ on
 diamonds and prove a contraction
theorem for the nearest-point projections to diamonds in euclidean buildings (Theorem \ref{thm:contrprojdiamo}).
In section \ref{sec:Regular and coarse regular} we prove the main result of our paper, Theorem \ref{main-tau} and
its corollaries, including continuous behavior of ends of Morse quasigeodesic rays. Lastly, in section
\ref{sec:appl} we establish structural results for regular and coarsely regular subsets of model spaces, and prove
Theorems \ref{thm:main1} and \ref{thm:morse subgroups}.  

\begin{rem}
\label{rem:hist}
After the first version \cite{mlem} of this paper was posted on the mathematics
arXiv in November 2014, the first version of the preprint \cite{GGKW}  appeared in February  2015.
A large part of the main results in this preprint follows from  \cite{mlem}
and our earlier preprint \cite{morse}:  

1) The first characterizations of Anosov representations of a hyperbolic group
that do not involve the geodesic flow of the group have been given in
\cite[section 6.5]{morse}. In \cite{morse} we have in fact a list of several equivalent
characterizations, analogous to characterizations of convex--cocompactness for Kleinian groups, see Theorem 1.7 in \cite{morse}. A further coarse geometric
characterization of Anosov subgroups, namely as coarsely $\taumod$-regular undistorted finitely generated
subgroups, has been given in \cite{mlem}, see Theorem 1.5 and Corollary 1.6, based on the
characterization of Anosov subgroups as Morse subgroups, which is part of
 \cite[Theorem 1.7]{morse}.  Note also that, unlike \cite{GGKW}, some of our characterizations do not require
the discrete group to be a priori Gromov--hyperbolic; the hyperbolicity turns
out to be a consequence. The existence of the equivariant continuous antipodal
boundary map has also been established in \cite[Theorem 1.4 and Remark 6.16]{mlem}. 

2) Part (2) of Theorem 1.1 in \cite{GGKW}, along with the positive answer to  
Question 5.5 on page 40 in \cite{GGKW},  
is a special case of \cite[Theorem 1.5]{mlem} (where $\Gamma$ is assumed to be
only finitely generated but not a priori hyperbolic), again based on the
characterization of Anosov subgroups as Morse subgroups.

3) Theorem 1.3 in \cite{GGKW} is not new. It follows entirely from \cite{morse,
mlem}.  The equivalence of the conditions (1), (2) and (3) in Theorem 1.3 follows from \cite[Theorem 6.57]{morse}. 
In fact, we deduce the Anosov property already from a weaker expansion
condition, cf. \cite[Definition  6.46]{morse}. The 
 CLI condition (4) in  Theorem 1.3, for hyperbolic groups, is easily seen (by
 elementary euclidean geometry) 
 to be equivalent to our condition of being {\em coarsely $\taumod$-regular and
undistorted}.  The latter condition is equivalent to the $\taumod$-Morse condition by
\cite[Theorem 1.5]{mlem}, and hence to the $\taumod$-Anosov condition, which is
condition (1) in Theorem 1.3. 

4) According to Remark 1.5 (1) in \cite[page 4]{GGKW},  part 4 of Theorem 1.3 in
\cite{GGKW} 
provides a definition of Anosov representations into $p$-adic groups. In fact,
such a definition has been provided before 
in \cite[Definition 5.35]{mlem}. There we have a definition of a Morse action on
a euclidean building (not necessarily discrete). Moreover, we prove Gromov hyperbolicity
of finitely generated groups admitting Morse actions on buildings, as well as 
the existence of an equivariant continuous antipodal boundary map for such actions (Theorem
6.15 and Remark 6.16 of \cite{mlem}). 

5) We would like to emphasize that the main results of our paper are coarse geometric in nature, describing the geometry of coarsely regular quasigeodesics and, more generally, of quasigeodesically connected coarsely regular  subsets of symmetric spaces and euclidean buildings. Results concerning discrete subgroups appear as direct applications. In contrast, the corresponding results of \cite{GGKW} are restricted to an equivariant setting and 
deal only with actions of hyperbolic groups on symmetric spaces.  
\end{rem}
\begin{rem}
\label{rem:hist2}
Concerning the misleading Remark 1.4(b) in the second version of \cite{GGKW} from March 2015,
the ``sharing of ideas'' during the year 2014 consisted of the following: 
In January 2014, the second author of this paper explained in a colloquium talk in Heidelberg in the presence of Wienhard our definition of Anosov representations of hyperbolic groups
not involving geodesic flows, given in \cite[section 6.5]{morse}. In July 2014, 
the first author explained to Gu\'eritaud and Kassel in detail the results of our paper \cite{morse}. 
In August 2014, the third author was present at a conference talk by Kassel where she announced some of the results which became available 
in February 2015 in the first version of the preprint \cite{GGKW}.
\end{rem}

\medskip 
{\bf Acknowledgements.} 
The first author was supported by the NSF grant DMS-12-05312. 
The last author was supported by the grant Mineco MTM2012-34834.

\section{Preliminaries}\label{sec:prelim}

This section contains some background material on metric geometry, geometry of buildings and ultralimits. 
We refer the reader to \cite{Burago-Burago-Ivanov} and \cite{Ballmann} for further reading on metric and CAT(0) geometry, 
and to \cite[ch. 3]{cone}, \cite[ch. 2.4]{qirigid} or \cite[ch. 7]{DrutuKapovich} for a discussion of the notion of ultralimits, asymptotic cones and their basic properties. 

\subsection{Metric spaces}

Let $(Z,d)$ be a metric space. 
We let denote $B(z,r)$ and $\bar B(z,r)$ the open and closed $r$-balls respectively, centered at a point $z\in Z$.
For a subset $A\subset Z$,
we denote by $\rad(A,z)$ its {\em radius} with respect to the center $z$,
i.e. the minimal $r\in[0,+\infty]$ such that $A\subset\bar B(z,r)$.
For a subset $Z'\subset Z$, we let $N_D(Z')$ denote the open $D$-neighborhood of $Z'$ in $Z$. 
We will use the notation $L(c)$ for the length of a (rectifiable) path $c$ in $Z$. 
\begin{dfn}[Almost distance minimizing path]
\label{dfn:almdmin}
We say that a path $c:[a,b]\to Z$ is {\em $\eps$-distance minimizing} if 
\begin{equation*}
L(c)\leq(1+\eps)\cdot d(c(a),c(b))
\end{equation*}
\end{dfn}
\begin{lem}
\label{lem:distminsubp}
Every rectifiable path contains, for arbitrarily small $\eps>0$, 
subpaths which are $\eps$-distance minimizing.
\end{lem}
\proof
Suppose that the path $c:[a,b]\to Z$ is rectifiable.
We choose a subdivision $a=t_0<t_1<\dots<t_k=b$
which almost yields the length of the path,
$$(1+\eps)\cdot\sum_{i=1}^k d(c(t_{i-1}),c(t_i)) \geq L(c)=\sum_{i=1}^k L(c|_{[t_{i-1},t_n]})$$
Then one of the subpaths $c|_{[t_{i-1},t_i]}$ is $\eps$-distance minimizing.
\qed

\medskip
We will use the term {\em pseudo-metric} for a distance function 
where different points may have infinite distance (however not distance zero).

\subsection{Spaces with curvature bounded above}

If $(Z,d)$ is a CAT(1) space, then a subset $C\subset Z$ is called  {\em convex}
if for any two points  
$\zeta_1, \zeta_2\in C$ with $d(\zeta_1, \zeta_2)<\pi$ 
the unique geodesic in $Z$ connecting $\zeta_1, \zeta_2$ is contained in $C$. 

\medskip
Suppose now that $X$ is a CAT(0) space. 

We will use the notation $xy$ for the unique geodesic segment in $X$ connecting $x$ to $y$. 
We will usually regard it as an {\em oriented} segment,
equipped with its natural orientation from the initial point $x$ to the endpoint $y$.
Similarly, given an ideal boundary point $\xi\in\geo X$ and a point $x\in X$, 
we let $x\xi$ denote the unique geodesic ray from $x$ asymptotic to $\xi$,
and $\xi x$ the same ray with the reversed orientation. 
We will denote by $x_0x_1\dots x_k$ the broken geodesic path 
which is the concatenation of the segments $x_{i-1}x_i$ for $i=1,\dots k$.
Similarly, we will denote by $x_0\dots x_k\xi_+$,
$\xi_-x_0\dots x_k$ and $\xi_-x_0\dots x_k\xi_+$ 
semi- and biinfinite paths obtained by attaching one or two rays at the ends of $x_0x_1\dots x_k$.

We will use the notation $\geo X$ for the {\em ideal} or {\em visual boundary} of $X$, 
equipped with the {\em visual topology}. 
It carries in addition another natural topology, called the {\em Tits topology},
which is finer than the visual topology; 
it is induced by a metric $\tangle$ on $\geo X$, called the {\em Tits metric}. For a subset $Y\subset X$ we let 
$\geo Y\subset\geo X$ denote the accumulation set of $Y$ in $\geo X$. For an oriented geodesic line $l$ in $X$, we let $l(\pm\infty)\in\geo X$ denote its ideal endpoints,
$$
l(\pm\infty)=\lim_{t\to\pm\infty} l(t),
$$
where $l: \R\to X$ is a (unit speed) parameterization of $l$ consistent with the orientation. Then 
$\geo l= \{l(-\infty), l(\infty)\}$. 
Similarly,
we denote the ideal endpoint of a ray $r\subset X$ by $r(+\infty)$.

For an ideal point $\xi\in\geo X$,
we denote by $b_{\xi}$ a {\em Busemann function} at $\xi$,
and by $Hb_{\xi,x}$ the {\em horoball} centered at $\xi$ and containing $x$ in its boundary horosphere,
i.e. $Hb_{\xi,x}=\{b_{\xi}\leq b_{\xi}(x)\}$.

We will say that two segments $xy$ and $x'y'$ are {\em oriented $r$-close}
if their initial and endpoints are $r$-close, i.e.\ $d(x,x')\le r$ and $d(y,y')\leq r$. 
In view of the convexity of the distance function of CAT(0)-spaces, any two segments which are oriented $r$-close, 
are also within Hausdorff distance $\le r$ from each other. 

We will denote by $\De(x,y,z)$ the geodesic {\em triangle} with vertices $x,y,z\in X$,
i.e.\ the {\em one}-di\-men\-sio\-nal object $xy\cup yz\cup zx$.
If $\eta,\zeta\in\geo X$,
we denote by $\De(x,y,\zeta)$ the {\em ideal triangle} with vertices $x,y,\zeta$,
that is, the union $x\zeta\cup xy\cup y\zeta$,
and by $\De(x,\eta,\zeta)$ the {\em ideal hinge} $x\eta\cup x\zeta$ with vertices $x,\eta,\zeta$.
We say that a triangle (ideal triangle, hinge) is {\em rigid} 
or can be {\em filled in by a flat triangle (half-strip, sector)}
or {\em spans a flat triangle (half-strip, sector)},
if it is contained in a convex subset which is isometric to a convex subset of euclidean plane. 

We will denote by $\Si_xX$ the {\em space of directions} at a point $x\in X$; this is a replacement of the unit tangent sphere in a Riemannian manifold (see \cite{Burago-Burago-Ivanov} and \cite{qirigid} for the precise definition). 
The space $\Si_xX$ is a CAT(1) space
equipped with the {\em angular metric} denoted $\angle(\xi,\eta)$. Each geodesic segment $xy$ determines a 
direction $\oa{xy}\in \Si_xX$. We will use the notation  $\angle_x(y,z)$ for the {\em angle} at $x$ between the segments 
$xy$ and $xz$ in $X$, i.e.\ the distance in $\Si_xX$ between the directions $\oa{xy}$ and $\oa{xz}$. 
This notation extends to the case of semi-infinite geodesics in $X$: 
For a point $\xi\in\geo X$, we denote by $\angle_x(y,\xi)$
the angle between $xy$ and the geodesic ray $x\xi$. 
Furthermore, for a subset $A$ containing $x$,
we denote by $\angle_x(y,A)$ the angular distance between $\oa{xy}$ and $\Si_xA$ in $\Si_xX$.

The {\em initial velocity} $\dot\rho\in \Si_xX$ 
of a geodesic $\rho: \R_+\to X$ in $X$ is the direction of $\rho$ at the point $x=\rho(0)$.

\medskip 
For a closed convex subset $C\subset X$ we have the {\em nearest point projection} 
$$\pi_C: X\to C .$$ 
This projection is a 1-Lipschitz map.  

Consider the special situation 
when $X$ is a Riemannian CAT(0) space (a Hadamard manifold) 
and $C\subset X$ is a totally-geodesic subspace. 
Then the distance function 
$$
d(x, C)=\min_{y\in C} d(x,y) 
$$
is 1-Lipschitz and smooth outside of $C$; the gradient lines of this function 
are the geodesics $x\bar{x}$, $\bar{x}=\pi_C(x)$. 
Suppose that $r: [0,+\infty)\to X-C$ is a unit speed geodesic ray 
with ideal endpoint $r(+\infty)=\xi\in \geo X$. 
Then the function $f(t)= d(r(t), C)$ is smooth with derivative 
\begin{equation}
\label{eq:derivative}
f'(t)= -\cos(\angle_{r(t)}(\bar r(t), \xi)) ,
\end{equation}
where $\bar r=\pi_C\circ r$ denotes the projection of the ray.

\subsection{Buildings and symmetric spaces} 

In the paper we will be using {\em nonpositively curved symmetric spaces}, {\em spherical} and {\em euclidean} buildings. We regard Riemannian symmetric spaces of noncompact type, respectively, euclidean buildings 
as the smooth ``archimedean", respectively, the singular ``non-archimedean" members 
of the family of CAT(0) ``model spaces" with rigid geometry. Both symmetric spaces and euclidean buildings 
will usually be denoted by $X$, while spherical buildings will be denoted by $B$. 
We will only consider symmetric spaces and euclidean buildings $X$ of {\em noncompact type}, which means that  $X$ is $CAT(0)$ and has {\em no flat factor},
i.e.\ is not isometric to the direct product of metric spaces $X'\times \R^k$ with $ k\ge 1$.

\begin{dfn}
By a {\em model space}, we mean a symmetric space of noncompact type or a euclidean building of noncompact type.
\end{dfn}

 We rule out flat factors for our model spaces, in part, 
because, as far as the results discussed in this paper are concerned, 
the case of spaces with a flat factor immediately reduces to the case without. However, many arguments in the paper use {\em parallel sets} of geodesics or flats in model spaces: These parallel sets 
do have flat factors and, hence, are CAT(0) symmetric spaces and euclidean buildings 
which do not have noncompact type.  

The two types of model spaces are connected via {\em asymptotic cones}; this connection will be explained in section \ref{sec:ulim}. 

For a treatment of buildings from the metric perspective of spaces with curvature bounded above,
we refer to \cite[ch.\ 3-4]{qirigid}.
Some notions needed in this paper or closely related to it, 
have been discussed in the case of symmetric spaces in \cite[ch.\ 2+5.1]{morse},
and the discussion in the building case is very similar, and often simpler. It is important to stress here that 
the euclidean buildings we are considering are allowed to be {\em nondiscrete}
and in particular {\em not locally compact}; such buildings appear as asymptotic cones of symmetric spaces of noncompact type.

\subsection{Spherical buildings}
\label{sec:sphbuil}

Instead of giving the precise definitions of spherical buildings (and euclidean buildings in the following section), we will describe below some of their important features. 
Part of this section is a review of the material in \cite[2.4.1-2]{morse}, to which we refer the reader for more details.

From the metric viewpoint, spherical buildings are CAT(1) spaces;
we will denote their metrics by $\angle$. 

A spherical building $\B$ has an associated {\em spherical Coxeter complex} $(\amod, W)$, 
where the {\em spherical model apartment} 
$\amod$ is a euclidean unit sphere and $W$ is a finite reflection group acting on $\amod$, called the {\em Coxeter} or {\em Weyl group} of $\B$. The quotient $\simod\cong\amod/W$ is called the {\em model chamber}. 
We identify it with a chamber in the model apartment, $\simod\subset\amod$.
We will say that the building $\B$ has {\em type} $\simod$.

As long as $W$ has no fixed point on $\amod$,
the model simplex $\simod$ is a spherical simplex in $\amod$ and has diameter $\leq\pihalf$.
We will use the notation $\taumod$ for faces of $\simod$ and $W_{\taumod}$ for the stabilizer of $\taumod$ in $W$. 
The {\em longest element} of the group $W$ is the unique element $w_o\in W$ which sends $\si_{mod}$ to $-\si_{mod}$ (the latter is also a chamber in $\amod$). The composition $\iota=-w_o$
preserves the model chamber $\simod$. (For some Weyl groups $W$, $w_o=id$, then $\iota=id$.) 

Each spherical building has a natural structure of a polysimplicial cell complex. Facets (top-dimensional faces) of this complex are called {\em chambers} of $\B$. 
Each building $\B$ comes equipped with a system (``atlas") of isometric embeddings $\amod\to \B$, 
whose images are called (spherical) {\em apartments}. Any two points of $\B$ belong to an apartment. 
It is important to stress that the spherical buildings in this paper are {\em not assumed to be thick},
i.e.\ a codimension one face may be adjacent to only two chambers.

A splitting of the model chamber 
as a spherical join $\simod=\simod^1\circ\simod^2$,
equivalently, a splitting $(\amod, W)=(\amod^1, W_1)\circ(\amod^2, W_2)$ of the spherical Coxeter complex, 
induces {\em splittings} of all buildings $\B$ of type $\simod$ 
as spherical joins $\B=\B_1\circ\B_2$ of spherical buildings $\B_i$ of types $\simod^i$.

Two faces $\bar\tau_+, \bar\tau_-\subset\amod$ are called {\em antipodal} or {\em opposite} 
if $- \bar\tau_- =\bar\tau_+$. Similarly, two points 
$\bar\xi, \bar\xi'\in \amod$ are antipodal if $\bar\xi'= -\bar\xi$. These definitions extend to the entire building $\B$ since any two faces (and any two points) are contained in an apartment in $\B$.  

\medskip
In a general simplicial complex $\Si$, 
we define the {\em interior} $\inte(\tau)$ of a simplex $\tau$
as the corresponding open face. 
We define the {\em star} $\st(\tau)\subset\Si$ of $\tau$
as the union of all (closed) faces containing $\tau$.
We note that the star is also known as the {\em residue}; this notion of the star should not be confused with the smallest subcomplex of $\Si$ consisting of faces which have nonempty intersection with $\tau$.
We define the {\em open star} of $\tau$, 
$$\ost(\tau)\subset\st(\tau),$$ 
as the union of all open faces whose closure contains $\tau$.
Furthermore, we define the {\em boundary} of the star,
$$\D\st(\tau):=\st(\tau)-\ost(\tau);$$ 
it is the union of all (closed) faces of the star, which do not contain $\tau$. 
If the simplex $\tau$ is maximal,
i.e.\ not contained in a simplex of larger dimension,
then $\st(\tau)=\tau$, $\ost(\tau)=\inte(\tau)$ is the open face, 
and $\D\st(\tau)=\D\tau$ is the topological frontier of $\tau$. 
We will apply these notions to spherical buildings and their model chambers,
which both carry natural structures as simplicial complexes.

\medskip
There exists a canonical projection 
$$\theta:\B\to\simod$$  called the {\em type map}. The type map restricts to an isometry on each chamber of $\B$ and, hence, is  1-Lipschitz.   A {\em type} is a point in $\simod$, 
and a {\em face type} is a face of $\simod$. The type of a point $\xi\in \B$ is $\theta(\xi)$,
and the type of a face $\tau\subset \B$ is $\theta(\tau)$. 
If the simplices $\tau_{\pm}$ in $\B$ are opposite to each other, then $\theta(\tau_-)=\iota\theta(\tau_+)$. 
We call a type $\bar\xi\in\simod$ a {\em root type} if the ball $\bar B(\bar\xi,\pihalf)\subset\amod$ is a subcomplex,
equivalently, if the great sphere $S(\bar\xi,\pihalf)\subset\amod$ is a wall.

\medskip 
Throughout the paper, we will denote by $\taumod\subset \simod$ a face type.

We denote by $\Flagt(\B)$ the {\em flag space} of type $\taumod$ simplices in $\B$.
It is a discrete space.
If $\B$ carries an additional structure as a {\em topological building},
as do Tits boundaries of model spaces, compare below,
then the flag spaces inherit a topology.

A point $\xi\in \B$ is called {\em $\taumod$-regular} if $\theta(\xi)\in\ost(\taumod)$
and {\em $\taumod$-singular} if $\theta(\xi)\in\D\st(\taumod)$.
We call the $\simod$-regular points simply {\em regular};
these are the points with type in $\inte(\simod)$.
The $\taumod$-regular points in $\B$ form an open subset, 
whose connected components are the open stars $\ost(\tau)$ of the type $\taumod$ faces $\tau$.
For a $\taumod^{\pm}$-regular point $\xi\in \B$ we define $\tau_{\pm}(\xi)$ 
as the type $\taumod^{\pm}$ face such that $\xi\in\ost(\tau_{\pm}(\xi))$;
we set $\tau(\xi)=\tau_+(\xi)$. 

A subset $A\subset \simod$ is called {\em $\taumod$-convex} (or {\em Weyl-convex}) 
if its symmetrization $W_{\tau_{mod}}A\subseteq\st(\si_{mod})$ 
is a convex subset of $a_{mod}$, cf.\ \cite[Def.\ 2.15]{morse}.
By $\Theta, \Theta', \Theta''$ we will always denote 
compact $\taumod$-convex subsets of 
$\ost(\taumod)\subset\simod$. 
Note that $\taumod$ is determined by such a subset $\Theta$,
namely, as the smallest face whose interior intersects $\Theta$. 
When we use several such subsets $\Theta,\Theta',\Theta''$,
we will always assume that $\Theta\subset\inte(\Theta')$ and $\Theta'\subset\inte(\Theta'')$. 

Since $\diam(\simod)\leq\pihalf$,
for every type $\bar\xi\in\taumod$ there exists a radius $\rho=\rho(\Theta,\bar\xi)<\pihalf$
such that:
\begin{equation}
\label{eq:thetinba}
\Theta\subset\bar B(\bar\xi,\rho)
\end{equation}
The following constant will frequently occur:
\begin{equation}
\label{eq:eps0}
\eps_0(\Theta):=\angle(\Theta,\D\st(\taumod))= \min \{\angle(\eta, \zeta): \eta\in \Theta, \zeta\in  \D\st(\taumod)\}>0
\end{equation}
Sometimes we will also use:
\begin{equation}
\label{eq:eps0th}
\eps_0(\Theta,\Theta'):=\angle(\Theta,\st(\taumod)-\Theta')
>0
\end{equation}
A point $\xi\in \B$ is called {\em $\Theta$-regular}, if $\theta(\xi)\in\Theta$.
We define the {\em $\Theta$-star} of a type $\taumod$ simplex $\tau\subset\B$
as the set of $\Theta$-regular points in its star,
$\stTh(\tau)=\st(\tau)\cap\theta^{-1}(\Theta)$.
We will often use the fact that the $\Theta$-stars are uniformly separated from each other:
\begin{lem}
For any two distinct type $\taumod$ simplices $\tau_1,\tau_2\subset\B$,
the (nearest point) distance between $\st_{\Theta}(\tau_1)$ and $\st(\tau_2)$ 
is $\geq\eps_0(\Theta)$.
\end{lem}
\proof
Since the open stars are disjoint,
any path connecting a point in $\st_{\Theta}(\tau_1)$ to a point in $\st(\tau_2)$
must exit $\st(\tau_1)$ at its boundary. 
It therefore has a subpath which projects via the type map $\theta$ to a path in $\simod$ 
connecting a point in $\Theta$ to a point in $\D\ost(\taumod)$.
The assertion follows because $\theta$ is 1-Lipschitz.
\qed

\medskip
We will always use the conventions
$$\taumod^+:=\taumod, \quad \taumod^-:=\iota\taumod$$ 
and
$$\Theta_+:=\Theta ,\quad \Theta_-:=\iota\Theta.$$
A {\em singular sphere} in a spherical building $\B$ is an  isometrically embedded 
(eulicdean unit) sphere $s\subset\B$ which is, at the same time, a subcomplex of $\B$. Each singular sphere equals the intersection of some (possibly one) 
aparatments in $\B$. 

For an ordered pair of opposite simplices $\tau_{\pm}\subset \B$,
we denote by $s(\tau_-,\tau_+)\subset\B$ the singular sphere spanned by $\tau_{\pm}$,
i.e.\ containing them as top-dimensional simplices. Equivalently, $s(\tau_-,\tau_+)$ is the smallest (with respect to inclusion) isometrically embedded sphere in $\B$ containing $\tau_+\cup \tau_-$. Each singular sphere $s\subset \B$ has the form $s=s(\tau_-,\tau_+)$ for a pair of antipodal simplices $\tau_{\pm}$.

Given a singular sphere $s\subset \B$, we let ${\mathcal B}(s)\subset \B$  denote the subbuilding 
which is the union of all apartments containing $s$. 
There is a natural decomposition 
\begin{equation}
\label{eq:sphjdec}
{\mathcal B}(s)\cong s\circ \Si_s\B 
\end{equation}
as the spherical join of the sphere $s$ and its {\em link} $\Si_s\B$ in $\B$.
In the case when $s=s(\tau_-,\tau_+)$, 
we will use the notation ${\mathcal B}(\tau_-,\tau_+)$ for ${\cal B}(s)$. 
When we want to specify the ambient building ${\mathrm B}$,
we put it as a subscript and write ${\mathcal B}_{\mathrm B}(s)$.

The following properties will be often used:

(i) Each apartment $a\subset {\mathcal B}(\tau_-, \tau_+)$ contains $s=s(\tau_-,\tau_+)$. 

(ii) $\st(\tau_{\pm})\subset {\mathcal B}(\tau_-,\tau_+)$. 

(iii) $\ost(\tau_{\pm})$ is open in $\B$; in particular, $\ost(\tau_{\pm})$ is open in ${\mathcal B}(\tau_-,\tau_+)$. 

In view of the spherical join decomposition,
it is clear that every point in $s$ has inside ${\mathcal B}(s)$ a unique antipode,
and this antipode lies in $s$.
\begin{lem}
\label{lem:antipbdpar}
All antipodes $\xi_-\in{\mathcal B}(s)$ of a point $\xi_+\in\st(\tau_+)$
are contained in $\st(\tau_-)$.
Moreover, if $\xi_+\in\ost(\tau_+)$, then $\xi_-\in\ost(\tau_-)$.
\end{lem}
\proof
Let $\xi_+\in\st(\tau_+)$, and let $\xi_-\in {\mathcal B}(s)$ be an antipode of $\xi_+$.
Since ${\mathcal B}(\tau_-,\tau_+)$ is a subbuilding, 
the pair of antipodes $\xi_{\pm}$ is contained in an apartment $a\subset {\mathcal B}(s)$.
As for all apartments in ${\mathcal B}(s)$,
we have that $\tau_{\pm}\subset a$.
There exists  a chamber $\si_+\subset a$ containing $\xi_+$ with face $\tau_+$.
The opposite chamber $\si_-$ in $a$ contains $\xi_-$ and has $\tau_-$ as a face.
Thus $\xi_-\in\st(\tau_-)$.
The assertion for open stars follows.
\qed

\medskip
The last observation extends to almost antipodes in a quantitative manner.
\begin{lem}
\label{lem:almantipbdpar}
Let $\xi_+\in\stTh(\tau_+)$ and $\eta_-\in {\mathcal B}(s)$ 
be points such that $\angle(\xi_+,\eta_-)>\pi-\eps_0(\Theta)$.
Then $\eta_-\in\ost(\tau_-)$.
\end{lem}
\proof
We only need to treat the case when $\xi_+$ and $\eta_-$ are not opposite.
The geodesic arc $\eta_-\xi_+$ extends
to an arc $\eta_-\xi_+\eta_+$ of length $\pi$.
It connects $\eta_-$ to an antipode $\eta_+$.
Since $\angle(\xi_+,\eta_+)<\eps_0(\Theta)$,
the arc $\xi_+\eta_+$ is too short to leave $\ost(\tau_+)$ and 
therefore $\eta_+\in\ost(\tau_+)$.
The previous lemma then implies that $\eta_-\in\ost(\tau_-)$.
\qed

\begin{cor}
\label{cor:unifangsusp}
Let $\xi_+\in\stTh(\tau_+)$ and let $\zeta_-\in \B$ 
be an antipode of $\xi_+$ outside ${\mathcal B}(s)$.
Then $\angle(\zeta_-,{\mathcal B}(s))\geq\eps_0(\Theta)$.
\end{cor}
\proof
Suppose that $\angle(\zeta_-,{\mathcal B}(s))<\eps_0(\Theta)$
and let $\bar\zeta_-\in {\mathcal B}(s)$ be the nearest point projection of $\zeta_-$ to ${\mathcal B}(s)$.
(Note that, as a subbuilding, ${\mathcal B}(s)$ is a closed convex subset 
of $\B$,
and the nearest point projection to ${\mathcal B}(s)$ is well-defined on the open $\pihalf$-neighborhood.)
Since $\ost(\tau_-)\subset {\mathcal B}(s)$ is open in $\B$, 
it cannot contain the projection 
of a point outside ${\mathcal B}(s)$,
and hence $\bar\zeta_-\not\in\ost(\tau_-)$.
On the other hand, 
we have 
$\angle(\xi_+,\bar\zeta_-)>\pi-\eps_0(\Theta)$,
which leads to a contradition with the previous lemma.
\qed

\medskip
It has been proven in \cite[2.5.2]{morse} that stars and $\Theta$-stars of simplices are convex.
This follows from the fact that they can be represented as intersections of balls with radius $\pihalf$.
More precisely, one has:
\begin{prop}[Convexity of stars, cf. {\cite[Lemma 2.12]{morse}}]
\label{prop:stconv}
Let $\tau\subset\B$ be a simplex. 

(i) For every simplex $\hat\tau\subset\B$ opposite to $\tau$,
the star $\st(\tau)$ is the intersection of ${\mathcal B}(\hat\tau,\tau)$ 
and the simplicial $\pihalf$-balls whose interior contains $\inte(\tau)$
and whose center lies in ${\mathcal B}(\hat\tau,\tau)$.

(ii) $\stTh(\tau)$ equals the intersection of all $\pihalf$-balls containing it. 
\end{prop}

\subsection{CAT(0) model spaces}
\label{sec:modsp}

Similarly to spherical buildings, each model space $X$ has an associated {\em euclidean Coxeter complex} $(\Fmod, W_{aff})$, where the {\em model flat} (respectively, {\em apartment}) $\Fmod$ is a euclidean space and $W_{aff}$ is a, possibly nondiscrete, group of  isometries of  $\Fmod$ generated  by reflections. 
The linear part of this group is a finite reflection group, called the {\em Weyl group} $W$ of $X$;
we pick a base point $0\in \Fmod$ and think of $W$ as acting on $\Fmod$ fixing $0$.
The dimension of $\Fmod$ is called the {\em rank} of $X$. The quotient $\Fmod/W$ will be denoted $\Delta$ or $\De_{euc}$ or $\Vmod$; it is called the {\em euclidean model Weyl chamber} of $X$. 
We identify it with a euclidean Weyl chamber with tip $0$ in the model flat, $\De\subset\Fmod$.

Each model space $X$ comes equipped with a system (``atlas") of isometric embeddings 
$$
\kappa^{-1}: \Fmod \to X.
$$
The images of the maps $\kappa^{-1}$ are the {\em maximal flats} in $X$.
(In the case when $X$ is a euclidean building, they are also called {\em apartments}.)
The inverse maps $\kappa$ are called {\em charts} for the maximal flats (respectively, apartments).
The charts are {\em compatible} in the sense
that for any two charts $\kappa_1,\kappa_2$ the transition function
$\kappa_1\circ\kappa_2^{-1}$ 
is the restriction of an element in $W_{aff}$.

Any two points in $X$ are contained in a maximal flat.

In addition to its usual distance function $d$, each model space comes equipped with 
a {\em $\Delta$-valued distance function} 
or {\em $\Delta$-distance}, denoted  $d_{\De}$. The function $d_\Delta$ is defined on $\Fmod$ by 
\begin{equation*}
d_\Delta(x,y)= proj(y-x)\in\De
\end{equation*}
where $proj: \Fmod/W\cong\Delta$ is the quotient map. The function $d_\De$ extends to the entire model space $X$ 
due to the compatibility of apartment charts
and the fact that any two points are contained in a maximal flat.

The $\De$-distance satisfies the following symmetry property:
$$
d_\De(y,x)= \iota d_\De(x,y)
$$
If $X$ is a symmetric space, 
then $d_\De(x,y)$ completely determines the $\Isom_o(X)$-congruence class of the pair $(x,y)$,
i.e.\ $d_\De(x,y)= d_\De(x',y')$ if and only if there exists $g\in \Isom_o(X)$ such that  
$g(x)=x'$ and $g(y)=y'$.

The projection
\begin{equation}
\label{eq:dele}
X\times X\to\De, \quad (x,y)\mapsto d_\De(x,y)
\end{equation}
is 1-Lipschitz in each of the two variables, which implies the {\em triangle inequality} for $\De$-lengths: 
\begin{equation}
\label{eq:tridele}
\|d_\De(x, y) - d_\De(x, y')\|\le \| d_\De(y, y')  \|= d(y, y')   
\end{equation}
and 
\begin{equation*}
\|d_{\De}(x,y)-d_{\De}(x',y)\|\leq \| d_\De(x,x')  \|= d(x,x'),
\end{equation*}
where the differences of $\De$-lengths are taken in $\Fmod$,
viewed as a vector space with origin $0$,
see \cite{KLM}.

\medskip
Spherical buildings appear naturally 
when one looks at the geometry at infinity of a model space and, in the euclidean building case,
at the infinitesimal geometry:

(i) The visual boundary $\geo X$ of a model space $X$, equipped with the Tits metric $\tangle$, 
has a natural structure of a spherical building; we will refer to  this spherical building as the {\em Tits boundary} $\tits X$ of $X$. The Weyl group of $X$ is canonically isomorphic to the Weyl group of $\tits X$; the dimension of $\tits X$ equals $\rank(X)-1$. The euclidean Weyl chamber $\De$ of $X$ is canonically isometric to the complete euclidean cone over $\si_{mod}$.  If $X$ is a symmetric space then the building $\tits X$ is always thick, 
while if $X$ is a euclidean building then $\tits X$ is thick provided that $X$ is thick. 
We will say that the model space $X$ is of {\em type} $\simod$. The chamber $\simod$ determines the Coxeter complex $(\Fmod, W_{aff})$ of $X$ if $W_{aff}$ acts transitively on $\Fmod$ (which is the case of symmetric spaces and their asymptotic cones); in general, $\simod$ determines $\Fmod$ and the Weyl group $W$. 

(ii) In the same vein, for each euclidean building  $X$ and each point $x\in X$, the space of directions $\Si_xX$, 
equipped with the angle metric $\angle_x$, 
has a natural structure of a spherical building of the {\em same type} $\simod$,
equivalently, with the {\em same associated Coxeter complex} $(\amod, W)$ as 
$\tits X$. 
Note that in general the spherical building $\Si_xX$ is {\em not thick}.
(For instance, 
if $X$ is a discrete euclidean building and $x$ is not a vertex.)

We denote by $\theta:\tits X\to\simod$ and $\theta_x:\Si_xX\to\simod$ 
the natural type maps, and by
$$\log_x:\tits X\to\Si_xX$$ 
the natural 1-Lipschitz {\em logarithm} map, sending an ideal point 
$\xi$ to the direction $\oa{x\xi}$. This map sends faces isometrically onto faces and satisfies
$$
\theta=\theta_x\circ\log_x.
$$

(ii') If $X$ is a symmetric space,
then the spaces of directions $\Si_xX$ are unit spheres 
and the logarithm maps $\log_x$ are bijective and homeomorphisms 
with respect to the visual topology on $\geo X$.
One can pull back the Tits metric and the spherical building structure to $\Si_xX$ 
and then also speak of simplices, chambers, apartments etc.\ in $\Si_xX$.

\medskip
Along with these spherical buildings associated to $X$,
we have the {\em flag spaces at infinity} 
$\Dt X=\Flagt(\geo X)$
and the spaces $\Flagt(\Si_xX)$ of {\em infinitesimal flags},
cf. section~\ref{sec:sphbuil}.
The visual topology on $\geo X$ induces visual topologies on the flag spaces at infinity.
(This is emphasized by the notation $\Flagt(\geo X)$ instead of $\Flagt(\tits X)$.)

For a type $\bar\xi\in\inte(\taumod)$,
the natural identification $$\Flagt(\geo X)\cong\theta^{-1}(\bar\xi)\subset\geo X,$$
which assigns to a type $\taumod$ simplex the point of type $\bar\xi$ in its interior,
is a topological embedding.
The infinitesimal flag spaces are discrete in the euclidean building case,
while in the symmetric space case, they inherit natural (manifold) topologies 
from the unit tangent spheres. 
In the symmetric space case, 
these flag spaces at infinity are flag {\em manifolds};
$\Flagt(\geo X)$ is naturally homeomorphic to the (generalized partial) flag manifold $G/P$, 
where $G=Isom_o(X)$ and $P$ is a parabolic subgroup stabilizing a simplex of type $\taumod$ in $\geo X$. 
The infinitesimal flag manifolds are homeomorphic to the flag manifolds at infinity of the corresponding types. 

\medskip
A spherical join splitting $\simod=\simod^1\circ\simod^2$ of the model chamber 
induces {\em splittings} 
of all model spaces $X$ of type $\simod$ 
as products $$X=X_1\times X_2$$
of model spaces $X_i$ of types $\simod^i$,
compare section~\ref{sec:sphbuil}.

\medskip
If $xy\subset X$ is a nondegenerate segment, 
then we call $\theta(\oa{xy})$ its {\em type}.
Similarly, an oriented geodesic $l\subset X$ is said to have {\em type} $\theta(l(+\infty))$.

A nondegenerate segment $xy$, respectively, 
a pair $(x,y)$ of distinct points is called
{\em $\taumod$-regular}, respectively {\em $\Theta$-regular},
if its direction $\oa{xy}$ is. 
In this case, we define its {\em $\taumod$-direction} $\tau(xy)$ at $x$
as the type $\taumod$ face $\tau(\oa{xy})\subset\Si_xX$; 
then $\oa{xy}\in\ost(\tau(xy))$.
Analogously, we denote by $\tau_{\pm}(zw)$ 
the $\taumod^{\pm}$-direction of a $\taumod^{\pm}$-regular segment $zw$.

We denote by $$\Isomth(X)<\Isom(X)$$
the subgroup of {\em type preserving} isometries,
i.e. isometries which preserve the types of segments and ideal boundary points.
Note that $\Isomth(X)$ has finite index in $\Isom(X)$, 
because $X$ has no flat factor.

Since there is a unique geodesic segment connecting any two points in $X$,
we can identify the space of oriented segments in $X$ 
with the space $X\times X$, equipped with the product topology. 
We observe that $\taumod$-regularity is an {\em open} condition for oriented segments, 
because the type of a segment varies continuously with it.

\medskip
The phenomenon of {\em angle rigidity} is specific to euclidean buildings, see \cite[\S 4.1]{qirigid}.
In the case of symmetric spaces, one only encounters it at infinity, in the Tits boundary.
It is useful to keep in mind the following basic consequences of angle rigidity.

Two nondegenerate segments $xy_1,xy_2\subset X$ with the same initial point initially span a flat triangle,
i.e.\ there exist points $x\neq y'_i\in xy_i$
such that the geodesic triangle $\De(x,y'_1,y'_2)$ can be filled in by a flat triangle.
In particular, 
if the initial directions of the segments agree,
$\angle_x(y_1,y_2)=0$,
then the segments initially agree,
i.e.\ $xy_1\cap xy_2$ is a nondegenerate segment. 

More is true:
For any ray $x\eta_1$ and any nondegenerate segment $xy_2$ with the same initial point 
there exists a point $x\neq y'_2\in xy_2$ such that the ideal triangle $\De(x,\eta_1,y'_2)$ 
can be filled in by a flat half-strip.
Furthermore, $xy'_2$ can be extended to a ray $x\eta'_2$ such that the ideal hinge $\De(x,\eta_1,\eta'_2)$ 
can be filled in by a flat sector.

\medskip
We return to the discussion of model spaces in general.

The logarithm maps send stars {\em onto} stars:
\begin{lem}
\label{lem:starsontostars}
For each point $x\in X$ and simplex $\tau\subset\geo X$, it holds: 
$$
\log_x\st(\tau)=\st(\log_x\tau),  \quad
\log_x\ost(\tau)=\ost(\log_x\tau)
$$
\end{lem}
\proof
In the symmetric space case, the assertion is tautological, since the logarithm maps are homeomorphisms.

In the euclidean building case, 
the assertion is a consequence of angle rigidity.
Only the surjectivity requires an argument.

Let $v=\oa{xy_2}\in\st(\log_x\tau)$,
and let $\eta_1\in\inte(\tau)$.
According to our discussion of angle rigidity, there exists $\eta_2\in\geo X$
such that $\oa{x\eta_2}=v$
and the ideal hinge $\De(x,\eta_1,\eta_2)$ 
can be filled in by a flat sector.
This means that 
$\angle_x(\eta_1,\eta_2)=\tangle(\eta_1,\eta_2)$
and the restriction of $\log_x$ to the arc $\eta_1\eta_2$ is an isometric embedding.
Since logarithm maps restrict to isometries on simplices,
and since $\oa{x\eta_1}$ and $v=\oa{x\eta_2}$ are contained in one chamber,
it follows that also $\eta_1\eta_2$ must be contained in one chamber,
i.e.\ $\eta_2\in\st(\tau)$. 
This shows the assertion for closed stars.

The assertion for open stars follows, because logarithm maps are type preserving. 
\qed

\medskip
Each apartment $a\subset\tits X$ is the ideal boundary of a unique maximal flat $F\subset X$.
More generally, 
each (isometrically embedded) unit sphere $s\subset \tits X$ is the ideal boundary of a flat $f\subset X$.
If $s$ is not an apartment, then the flat $f$ is not maximal and not unique. 
If the sphere $s$ is singular, then also the flat $f$ is singular,
i.e. is the intersection of some maximal flats in $X$. 

\medskip
{\bf Parallel sets in model spaces and spherical joins at infinity.} 
One defines the {\em parallel set} $P(s)\subset X$ of a unit singular sphere $s\subset\tits X$ 
as the union of the (parallel) flats 
with ideal boundary $s$.
Parallel sets are totally geodesic subspaces, respectively, euclidean subbuildings,
depending on whichever $X$ is,
and as such they carry themselves natural structures as symmetric spaces, respectively, euclidean buildings 
with the same associated Coxeter complex and 
of the same type $\simod$ as $X$.

As a consequence, 
geodesic segments in parallel sets are extendible to complete geodesics,
and tangent directions to parallel sets are represented by segments in the parallel set.

However, parallel sets are not model spaces in our sense, because they have flat factors.
The parallel set $P(s)$ splits isometrically as 
\begin{equation}
\label{eq:splitpar}
P(s)\cong f\times CS(s)
\end{equation}
where the slices $f\times pt$ are the flats with ideal boundary sphere $s$,
and the {\em cross section} $CS(s)$ is a symmetric space or euclidean building
with corank $\dim(f)=\dim(s)+1$,
$$\rank(X)=\dim(f)+\rank(CS(s)) .$$
The visual boundary of $P(s)$ is (underlying) 
the subbuilding ${\cal B}_{\tits X}(s)$ 
of $\tits X$ associated to the sphere $s$,
$$\tits P(s)={\cal B}_{\tits X}(s).$$
Accordingly, there is the natural spherical join decomposition
$$\tits P(s)\cong s\circ\tits CS(s),$$ 
where $\tits CS(s)$ is canonically identified with the link $\Si_s(\tits X)$ of $s$ in $\tits X$,
compare (\ref{eq:sphjdec}).

\medskip
Let $\tau_{\pm}\subset s$ be a pair of opposite simplices spanning $s$,
i.e. $s=s(\tau_-,\tau_+)$.
The subset $\ost(\tau_+)\subset P(s)$ is open in $\geo X$ with respect to the Tits topology,
but in general not with respect to the visual topology.
However:
\begin{lem}
$\ost(\tau_+)$ is open in $\geo P(s)$ also with respect to the visual topology.
\end{lem}
\proof
Let $\xi_+\in\ost(\tau_+)$, and let $\xi_-\in\ost(\tau_-)$ be an antipode.
Any ideal point $\eta_+\in\geo P(s)$ sufficiently close to $\xi_+$ is almost opposite to $\xi_-$
because of the lower semicontinuity of the Tits metric with respect to the visual topology.
Lemma~\ref{lem:almantipbdpar} then implies that $\eta_+\in\ost(\tau_+)$.
\qed

\medskip
As for the visual boundary,
we have an analogous description and splitting of the spaces of directions of parallel sets
as subbuildings of the spaces of directions of $X$.
(In the symmetric space case, this refers to the spherical building structures 
on the spaces of directions pulled back from the visual boundary by the logarithm maps, and is tautological.)
\begin{lem}
For $x\in P$, it holds that 
$$
\Si_xP={\mathcal B}_{\Si_xX}(\log_xs).$$
\end{lem}
\proof
We only need to consider the case when $X$ is a euclidean building. 

Every direction in $\Si_xP$ is tangent to a maximal flat $F\subset P$ through $x$.
Since the apartment $\geo F\subset\geo P$ contains the sphere $s$,
the apartment $\Si_xF\subset\Si_xP$ 
contains the sphere $\log_xs$.
Therefore 
$$v\in\Si_xF\subset {\mathcal B}_{\Si_xX}(\log_xs)$$
and, hence
$$
\Si_xP\subset{\mathcal B}_{\Si_xX}(\log_xs).$$

Vice versa, 
let $\tau_{\pm}\subset s$ be a pair of opposite simplices spanning $s$,
i.e. $s=s(\tau_-,\tau_+)$.
Since 
$$
\tits (P(s))={\mathcal B}_{\tits X}(s)\supset\st(\tau_{\pm}),
$$
it follows with Lemma~\ref{lem:starsontostars} that 
$$
\st(\log_x\tau_{\pm})=\log_x\st(\tau_{\pm})\subset\Si_xP.
$$
Since $\Si_xP$ is a subbuilding of $\Si_xX$, 
it must therefore contain all apartments containing $\log_xs$.
This shows the reverse inclusion.
\qed

\medskip
Note that for spheres $s\subset s'\subset\tits X$,
we have that $P(s)\supset P(s')$. 
If $s$ is not singular and $s'$ is the unique smallest singular sphere containing $s$,
then there is equality.

For a flat $f\subset X$, we define its parallel set as $P(f):=P(\geo f)$;
it is the union of all flats parallel to $f$.
For flats $f\subset f'$, it holds that $P(f)\supset P(f')$. 
Again, if $f$ is not singular and $f'$ is the unique smallest singular flat containing $f$,
then equality holds. 

When $s=s(\tau_-,\tau_+)$, we will use the notation $P=P(\tau_-,\tau_+)$ for $P(s)$. In this notation we emphasize 
that we regard $P$ as a parallel set together with a choice 
of an ordered pair $(\tau_-,\tau_+)$ of antipodal simplices in $\geo P$.  
One can think of this choice as a higher rank analogue of 
an orientation of a geodesic. 
We will say that the parallel set $P(\tau_-,\tau_+)$ has {\em type} $\theta(\tau_+)$.

Each parallel set of a flat (or a sphere at infinity) can also be represented as the parallel set of a geodesic line. 
Namely,
$P(\tau_-,\tau_+)=P(l)$ for every line $l$ 
with $l(\pm\infty)\in\inte(\tau_{\pm})$.

\medskip
Two ideal points $\xi_{\pm}\in \geo X$ are opposite, i.e.\ $\tangle(\xi_-, \xi_+)=\pi$,
if and only if there exists a geodesic line $l\subset X$ asymptotic to $\xi_{\pm}$,
i.e.\ $l(\pm\infty)=\xi_{\pm}$.
(Note that this is not true for general CAT(0) spaces.)
Two simplices $\tau_{\pm}\subset\geo X$ are opposite
if and only if there exists a line $l\subset X$ such that $l(\pm\infty)\in\inte(\tau_{\pm})$.
\begin{dfn}
[$x$-opposite]
\label{def:xopp}
We say that two (opposite) simplices $\tau_{\pm}\subset\geo X$ are {\em $x$-opposite}
if the simplices $\log_x\tau_{\pm}\subset\Si_xX$ are opposite.  
\end{dfn}
If $X$ is a symmetric space, this condition means that the differential $ds_x$ 
of the point reflection at $x$ (Cartan involution) $s_x: X\to X$ swaps $\tau_+$ and $\tau_-$. 
In this case, for every simplex there exists a unique $x$-opposite simplex.

\begin{lem}
Two opposite simplices $\tau_{\pm}\subset\geo X$ are  $x$-opposite if and only if 
$x\in P(\tau_-,\tau_+)$. 
\end{lem}
\proof 
If the simplices $\log_x\tau_{\pm}\subset\Si_xX$ are opposite,
then they contain a pair of opposite directions $\log_x\xi_{\pm}\in\inte(\log_x\tau_{\pm})$.
Hence, there exists a pair of antipodes $\xi_{\pm}\in\inte(\tau_{\pm})$ 
such that $\xi_-x\xi_+$ is a geodesic line.
It follows that $x\in P(\tau_-,\tau_+)$. 
The converse is clear.
\qed

\medskip
A spherical join splitting $\simod=\simod^1\circ\simod^2$ 
induces {\em splittings} of all model spaces $X$ of type $\simod$ 
as metric products $X=X_1\times X_2$ of model spaces $X_i$ of types $\simod^i$.

\medskip
{\bf Cones.} For a  subset $A\subset \geo X$ and a point $x\in X$ 
we let $V(x, A)\subset X$ be the complete {\em cone} over $A$ with tip $x$, i.e.\ the union of 
the geodesic rays $x\xi$ for all $\xi\in {A}$. If $A$ is closed with respect to the visual topology on $\geo X$, 
then the subset $V(x,A)$ is closed in $X$. The cones $V(x,A)$, in general, are not isometric to (euclidean) 
metric cones. However, if $A$ is contained in an apartment in $\geo X$, then 
$V(x,A)$ is canonically isometric to the complete euclidean cone over the set $A$,
equipped with the Tits metric.

In the special case when 
$\tau\subset\geo X$ is a simplex, the cone $V(x,\tau)$ is called a {\em euclidean Weyl sector} in $X$,
and if $\si\subset\geo X$ is a chamber, then $V(x,\si)$ is called a {\em euclidean Weyl chamber}. 
The {\em open} sector 
$\inte(V(x,\tau)):=V(x,\inte(\tau))-\{x\}$
is the subset of points 
where $V(x,\tau)$ is locally isometric to euclidean space (of dimension $\dim\tau+1$).
It is the {\em interior} of the sector $V(x,\tau)$
inside any minimal singular flat containing it.

For a simplex $\tau\subset \geo X$, the cone $V(x, \st(\tau))$ is called a {\em Weyl cone} in $X$. 
It is the union of the euclidean Weyl chambers $V(x,\si)$ over all chambers $\si\subset \geo X$ 
containing $\tau$ as a face.
If $\hat\tau$ is a simplex $x$-opposite to $\tau$, then 
$$V(x,\st(\tau))\subset P(\hat\tau,\tau).$$
We call such a parallel set an {\em ambient parallel set} for the Weyl cone. 
We will refer to the subset $V(x,\ost(\tau))-\{x\}\subset V(x,\st(\tau))$ as the {\em open} Weyl cone. 
It is the subset of points $y\in V(x,\st(\tau))$ 
whose spaces of directions $\Si_yV(x,\st(\tau))$ are spherical buildings. 

Another class of cones which we will use 
are the {\em $\Theta$-cones} $V(x, \st_\Theta(\tau))$.

\subsection{Trees}
\label{sec:tree}

We recall the geometric notion of tree:
\begin{dfn}[Metric tree]
A {\em metric tree} is a 0-hyperbolic geodesic metric space.
\end{dfn}
Note that euclidean buildings of rank one are metric trees.

We will use the following fact:
\begin{lem}
\label{lem:bilipttr}
Every path metric space bilipschitz homeomorphic to a metric tree is itself a metric tree.
\end{lem}
\proof
Suppose that $(T,d)$ is a metric tree, and that $d'$ is another path metric on $T$ 
which is bilipschitz equivalent to $d$.
Any two points in $T$ are connected by an embedded path,
and this path is unique up to reparametrization.
Moreover, it is $d$-rectifiable and therefore $d'$-rectifiable. 
Any non-embedded path with the same endpoints is at least as $d'$-long,
because its image contains the image of the embedded connecting path. 
It follows that $d'$-geodesics coincide, up to reparametrization, with $d$-geodesics.
Thus, any two points in $T$ can be connected by a unique $d'$-geodesic
and $d'$-geodesic triangles are tripods.
\qed

\subsection{Ultralimits}
\label{sec:ulim}

We let $\om$ denote a {\em nonprincipal ultrafilter} on the set $\N$ of natural numbers. For a map 
$h: \N\to K$ from $\N$ to a compact Hausdorff space, one defines the {\em ultralimit}
$$
\ulim h(n) =k\in K, 
$$
as the unique point $k\in K$ such that for every neighborhood $U$ of $k$ in $K$, the subset $h^{-1}(U)$ belongs to $\om$. 

Consider a sequence of  pointed metric spaces $(X_n,\star_n)$ parameterized by $\N$; 
we use the notation $\dist_{X_n}$ for the metric on $X_n$. 
The {\em ultralimit} 
$$
 (X_{\om},\star_{\om}) = \ulim_n (X_n,\star_n)
$$
of the sequence of pointed metric spaces $(X_n, \star_n)$ is a pointed metric space defined as follows. 
Define a pseudo-distance $\dist_\om$ on the product space $\prod_{n\in \N} X_n$ by the formula
$$
\dist_{\om}\left((x_n),(y_n)\right):=\ulim \left(n\mapsto \dist_{X_n}(x_n,y_n)\right)$$
where we take the {\em ultralimit} of the function
$n\mapsto \dist_{X_n}(x_n,y_n)$
with values in the compact space $[0,\infty]$. The function
$\dist_{\om}$ takes values in $[0,\infty]$. In order to convert this function to a metric, we first consider the subset 
$$
X^o_\om \subset \prod_{n\in \N} X_n
$$
consisting of sequences $(x_n)_{n\in \N}$ such that 
$$
\dist_\om((x_n), (\star_n))<\infty. 
$$
Then $\dist_\om$ restricted to $X^o_\om \times X^o_\om$ takes only finite values. Lastly, take the quotient of  
$X^o_\om$, where we identify points with zero $\dist_{\om}$--distance. The result is the ultralimit $X_\om$; we retain the notation 
$\dist_\om$ for the projection of the pseudo-distance from $X^o_\om$ to $X_\om$. 
Points $x_{\om}\in X_{\om}$ are thus represented by sequences $(x_n)$ of points $x_n\in X_n$;
abusing notation,
we will sometimes write $x_{\om}=(x_n)$. 
The natural base point of $X_{\om}$ is $\star_{\om}=(\star_n)$.

\medskip 
The ultralimits that we will be using in the paper are of very special kind. They are defined by starting with a fixed metric space $(X, \dist_X)$, taking sequences of base points $\star_n\in X$ and 
of scale factors $\la_n>0$ converging to $0$, 
and then setting 
$$
X_n=X, \quad \dist_{X_n}= \la_n \dist_X. 
$$ 
Such ultralimits are called {\em asymptotic cones} of $(X, \dist_X)$. By abusing the notation, we will abbreviate 
$(X, \la_n \dist_X)$ to $\la_n X$.

\medskip 
We will need a basic construction, which relates quasi-isometries and asymptotic  cones. 
Suppose that $(Y_n, \star'_n),(X_n, \star_n)$ are sequences of pointed metric spaces and that 
$$
f_n: Y_n\to X_n
$$
are $(L,A)$-quasiisometric embeddings such that
$$
\ulim \la_n\dist_{X_n}(f_n(\star'_n), \star_n)< +\infty. 
$$ 
Suppose that $(\la_n)$ is a sequence of positive numbers satisfying $\ulim \la_n=0$ 
and consider the ultralimits
$$
(Y_\om,\star'_{\om}) = \ulim (Y_n, \la_n \dist_{Y_n}, \star'_n), \quad 
(X_\om,\star_\om) = \ulim (X_n, \la_n \dist_{X_n}, \star_n).
$$
Then the induced map 
$$
f_\om: Y_\om \to X_\om, \quad f_\om((y_n))= (f_n(y_n))
$$
is well-defined.
The map $f_\om$ is called the {\em ultralimit} of the sequence of maps $(f_n)_{n\in \N}$. 

Since, with respect to the rescaled metrics, 
the maps $f_n$ are $(L,\la_nA)$-quasiisometric embeddings,
their ultralimit is a $(L,0)$-quasiisometric embedding:
\begin{lem}
\label{lem:ulimqiebil}
The map $f_\om$ is an $L$-bilipschitz embedding:
$$
L^{-1} \dist_{Y_\om}(y_\om, y'_\om)\le \dist_{X_\om}(f_\om(y_\om), f_\om(y'_\om)) \le L\, \dist_{Y_\om}(y_\om, y'_\om). 
$$ 
\end{lem}

We will use this lemma primarily to conclude that the ultralimit of a sequence of uniform quasigeodesics in a symmetric space (or a building) is a bilipschitz path in the asymptotic cone, while ultralimits of seqeucnes of flats  are flats.

The following construction  is a special case of the lemma. Suppose that 
$(X_n, \star_n)$ is a sequence of pointed  metric spaces 
with ultralimit $(X_\om, \star_\om)=\ulim(X_n, \star_n)$
and that $Y_n\subset X_n$ are subsets such that 
$$
\ulim \dist_{X_n}(\star_n, Y_n)<+\infty. 
$$ 
Define the {\em ultralimit of the sequence of subsets} $Y_n$, 
$$Y_\om = \ulim Y_n \subset X_\om ,$$ 
as the subset consisting of all points $y_\om\in X_\om$ represented by sequences $(y_n)_{n\in\N}, y_n\in Y_n$. 
Alternatively, one can describe $Y_\om$ as follows. 
For any sequence of base points $\star'_n\in Y_n$ with 
$\ulim dist_{X_n}(\star_n,\star'_n)<+\infty$,
there is a natural isometric embedding of ultralimits 
$$
\ulim (Y_n, \star'_n) \to\ulim (X_n, \star'_n) =(X_\om,\star'_\om) 
$$
where $\dist_{Y_n}$ is the restriction of the distance function from $X_n$ to $Y_n$,
and the image of the embedding coincides with $Y_\om$.

Since the ultralimit of any sequence of metric spaces is a {\em complete} metric space 
(cf. Lemma I.5.53 in \cite{BH} or Proposition 7.44 in \cite{DrutuKapovich}), 
it follows that the ultralimit of any sequence of subspaces is {\em closed}.

\section{Geometry of CAT(0) model spaces}
\label{sec:modgeom}

Throughout this chapter, $X$ denotes a model space.
When parts of the discussion apply only to euclidean buildings or symmetric spaces,
this will be indicated explicitly.

\subsection{Regularity and coarse regularity}
\label{sec:regcreg}

The regularity of pairs of points, equivalently, of segments has been defined in section~\ref{sec:modsp}. 

We call a {sequence} $(x_n)$ {in $X$} {\em $\Theta$-regular} if all pairs $(x_m,x_n)$ for $m<n$ are $\Theta$-regular; a path $c:I\to X$ is {\em $\Theta$-regular} if all pairs of points $(c(t_1),c(t_2))$ for $t_1<t_2$ are $\Theta$-regular. When we do not want to specify $\Theta$, we say that a sequence $(x_n)$ or a path $c$ is {\em uniformly $\taumod$-regular} if it is $\Theta$-regular for some $\Theta$. 

A weaker version of uniform regularity is {\em regularity}: A sequence $(x_n)$, resp. a path $c$ is {\em $\taumod$-regular} if  all pairs $(x_m,x_n)$ for $m<n$ are $\taumod$-regular, resp. 
all pairs of points $(c(t_1),c(t_2))$ for $t_1<t_2$ are $\taumod$-regular. Note that $\taumod$-regularity does {\em not} imply local uniform $\taumod$-regularity.

\medskip 
If $\taumod$ and $\Theta$ are $\iota$-invariant,
then the order of the points does not matter:
A segment is $\taumod$- or $\Theta$-regular if and only if the reversely oriented segment is.
Freed of the orientation issues,
we then say that a {\em subset} $R\subset X$ is {\em $\Theta$-regular} 
if any pair of distinct points in $R$ is $\Theta$-regular,
and more generally, that a {\em map} $Z\to X$ into $X$ is {\em $\Theta$-regular} 
if it sends any pair of distinct points in $Z$ to a $\Theta$-regular pair of points in $X$.
In the same way, we define the $\taumod$-regularity of subsets of and maps into $X$.
Note that regular maps are necessarily injective,
and their images are regular subsets.
Vice versa, injective maps into regular subsets are regular.

\medskip
A natural way to {\em coarsify} the notion of regularity is as follows.

Let $B\geq 0$.
We say that a {\em pair} $(x,y)$ of (not necessarily distinct) points is {\em $(\Theta,B)$-regular} 
if it is oriented $B$-close to some $\Theta$-regular pair of points $(x',y')$,
i.e. $d(x,x')\leq B$ and $d(y,y')\leq B$.
Since we are working in a CAT(0) setting,
this is equivalent to the property that the segment $xy$ is oriented 
$B$-Hausdorff close to the $\Theta$-regular segment $x'y'$,
and we say also that the {\em segment} $xy$ is $(\Theta,B)$-regular.

We say that a {\em sequence} $(x_n)$ in $X$ is $(\Theta,B)$-regular
if all pairs $(x_m,x_n)$ for $m<n$ are $(\Theta,B)$-regular.
Similarly, 
we say that a (not necessarily continuous) {\em path} $p:I\to X$ is $(\Theta,B)$-regular,
if for every subinterval $[a',b']\subset I$, the segment $p(a')p(b')$ is $(\Theta,B)$-regular.
We will primarily use this definition in the case of quasigeodesics (finite or infinite). 

If $\taumod$ and $\Theta$ are $\iota$-invariant,
then we say that a {\em subset} of $X$ is $(\Theta,B)$-regular
if every pair of points in the subset has this property,
and more generally, that a {\em map} into $X$ is $(\Theta,B)$-regular 
if it sends any pair of points to a $(\Theta,B)$-regular pair of points in $X$.
Note that the images of $(\Theta,B)$-regular maps are $(\Theta,B)$-regular subsets.
We say that the subset or map is {\em coarsely $\Theta$-regular}
if it is $(\Theta,B)$-regular for some constant $B$.
We say that an {\em isometric group action} on $X$ is coarsely $\Theta$-regular
if some (every) orbit map is.

A path, map, subset or action is said to be {\em coarsely uniformly $\taumod$-regular} 
if it is coarsely $\Theta$-regular for some $\Theta$.

\medskip
Let $\taumod$ and $\Theta$ again be unrestricted.

Here is a useful weakening of the notion of coarse uniform regularity:
\begin{dfn}[Asymptotically regular sequence, cf. {\cite[Def. 5.1]{morse}}]
We say that a sequence $x_n\to\infty$ in $X$ is {\em asymptotically $\Theta$-regular},
if for some (any) basepoint $x\in X$
the set of accumulation points of the sequence of direction types $\theta(\oa{xx_n})\in\simod$ 
is contained in $\Theta$, equivalently,
if the set of accumulation points of the sequence of $\De$-lengths $d_{\De}(x,x_n)\in\Vmod$
is contained in $\Theta\subset\simod\cong\geo\Vmod$.

A sequence in $X$ is called {\em asymptotically uniformly $\taumod$-regular}
if it is asymptotically $\Theta$-regular for some $\Theta$.
\end{dfn}

\begin{lem}\label{lem:ASREGULAR}
(i) The set of accumulation points in $\geo X$ 
of an asymptotically $\Theta$-regular sequence $x_n\to\infty$
is contained in the $\Theta$-regular part $\theta^{-1}(\Theta)\subset\geo X$
of the ideal boundary. 
If $X$ is locally compact, then the converse holds as well.  

(ii) If $x_n\to\infty$ is an asymptotically $\Theta$-regular sequence, then 
for every point $x\in X$ the segments $xx_n$ are $\Theta'$-regular for all sufficiently large $n$.

(iii) $(\Theta,B)$-regular sequences in $X$ are asymptotically $\Theta$-regular.
\end{lem}
\proof 
The first assertion of part (i) is clear.
For the second, suppose that $X$ is locally compact
and consider a sequence $x_n\to\infty$ 
which accumulates at a subset of $\theta^{-1}(\Theta)\subset\geo X$
and such that for a point $x\in X$,
after passing to a subsequence, 
the direction types $\theta(\oa{xx_n})$ converge,
$\theta(\oa{xx_n})\to\bar\xi\in\simod$.
After passing to a subsequence again, 
we may assume that also the sequence $(x_n)$ converges at infinity, $x_n\to\xi\in\geo X$.
It follows that $\bar\xi=\theta(\xi)\in\Theta$.
Thus $(x_n)$ is asymptotically $\Theta$-regular.

Parts (ii) and (iii) follow from the triangle inequality for $\De$-lengths (\ref{eq:tridele}).
\qed

\begin{dfn}[Asymptotically regular subset]
\label{dfn:rsubs}
We call a subset $R\subset X$ {\em asymptotically $\Theta$-regular} 
if all diverging sequences in $R$ have this property. 
\end{dfn}
We suppose again that $\taumod$ and $\Theta$ are $\iota$-invariant
and consider the concepts introduced so far in the context of discrete subgroups.
\begin{dfn}[Asymptotically regular subgroup and action]
\label{dfn:rsubgp}
We say that a discrete {\em subgroup} $\Ga<\Isom(X)$ is {\em asymptotically $\Theta$-regular}
if its orbits in $X$ have this property.
More generally, 
we call a properly discontinuous isometric {\em action} $\Ga\acts X$ of a discrete group $\Ga$ on $X$ 
{\em asymptotically $\Theta$-regular}
if its orbits in $X$ have this property.
\end{dfn}
\begin{rem}
\label{rem:asyregact}
(i) If $X$ is locally compact, then the asymptotic uniform $\taumod$-regularity of $\Ga$ is 
equivalent to the property that the limit set of $\Ga$ is contained in the $\taumod$-regular part of the 
visual boundary, $\theta(\La(\Ga))\subset\ost(\taumod)$, see Lemma \ref{lem:ASREGULAR},  
cf. \cite[Def. 5.1]{morse}. 
We recall that the limit set $\La(\Ga)\subset\geo X$ of $\Ga$ 
is the accumulation set of a $\Ga$-orbit $\Ga x\subset X$.

(ii) Coarsely $\Theta$-regular actions are also asymptotically $\Theta$-regular.
Asymptotically $\Theta$-regular actions are coarsely $\Theta'$-regular.
\end{rem}

\medskip
The next observations relate coarse regularity to regularity.
\begin{lem}[Long coarsely regular implies regular]
\label{lem:lcoaregureg}
There is a constant $c=c(\Theta,\Theta')>0$ such that
every $(\Theta,B)$-regular segment of length 
$\geq cB$ is $\Theta'$-regular.
\end{lem}
\proof Suppose that the segment $xy$ is oriented $B$-close to the $\Theta$-regular segment $x'y'$; 
define 
$$
D:= \max(d(x,y), d(x',y')).$$
The triangle inequality for $\De$-lengths
(\ref{eq:tridele}) yields that 
$|d_{\De}(x,y)-d_{\De}(x',y')|\leq 2B$.
It follows that the angular distance 
$$
\al= \angle(\theta(xy),\theta(x'y'))
$$ 
between the types of the segments  $xy$ and $x'y'$ satisfies
$$
\sin(\al/2)\leq \frac{B}{D}\leq \frac{B}{d(x,y) - 2B}.$$
The lemma follows. \qed

\medskip
We note that long chords of coarsely regular quasigeodesics are uniformly regular:
\begin{lem}
\label{lem:cregqgeolongsegmreg}
With the constant $c=c(\Theta,\Theta')>0$ from Lemma~\ref{lem:lcoaregureg} 
the following holds:

Suppose that $q:I\to X$  is a $(\Theta,B)$-regular $(L,A)$-quasigeodesic.
Then for every subinterval $[a',b']\subset I$ with length
$\geq L(A+cB)$,
the segment $q(a')q(b')$ is $\Theta'$-regular.
\end{lem}
\proof
The segment $q(a')q(b')$ has length $\geq cB$
and is therefore $\Theta'$-regular by Lemma~\ref{lem:lcoaregureg}.
\qed

\medskip
Similarly, one obtains the same conclusion 
for the projections of coarsely regular quasigeodesics to nearby parallel sets.
Let $P=P(\tau_-,\tau_+)$ be a type $\taumod$ parallel set,
and let $\bar q=\pi_P\circ q$ denote the nearest point projection of the path $q$ to $P$.

\begin{lem}
\label{lem:cregqgeolongsegmregproj}
With the constant $c=c(\Theta,\Theta')>0$ from Lemma~\ref{lem:lcoaregureg} 
the following holds:

Suppose that $q:I\to X$  is a $(\Theta,B)$-regular $(L,A)$-quasigeodesic 
such that $q(I)\subset\bar N_D(P)$.
Then for every subinterval $[a',b']\subset I$ with length 
$\geq L(A+c(B+D))$
the segment $\bar q(a')\bar q(b')\subset P$ is $\Theta'$-regular.
\end{lem}
\proof
The projected quasigeodesic $\bar q:I\to P$ is $(\Theta,B+D)$-regular.
(Its quasiisometry constants are irrelevant.)
As in the proof of the previous lemma
we note that 
the segment $\bar q(a')\bar q(b')$ has length $\geq c(B+D)$
and is therefore $\Theta'$-regular by Lemma~\ref{lem:lcoaregureg}.
\qed

\subsection{Longitudinality and coarse longitudinality}
\label{sec:longi}

Longitudinality is a property of segments and directions in a parallel set,
which is ``oriented'' by the choice of a pair of opposite simplices 
spanning the singular sphere factor of its visual boundary. 
It means that the segments or directions point towards the open stars of these simplices.
To prepare the precise definition,
we first need an observation 
which relates the property of pointing to these stars 
for directions, segments and rays.

Let $P=P(\tau_-,\tau_+)\subset X$ be a type $\taumod$ parallel set.
\begin{lem}
\label{lem:extintcon}
Let $xy\subset P$ be a nondegenerate segment
and let $x\xi_+\subset P$ be a ray.

(i) If $\oa{xy}\in\st(\log_x\tau_+)$, then $y\in V(x,\st(\tau_+))$.
If $\oa{xy}\in\ost(\log_x\tau_+)$, then $y\in V(x,\ost(\tau_+))$.
Moreover, $\oa{xy}\in\st(\log_x\tau_+)$ if and only if $\oa{xy}\in\st(\log_y\tau_-)$.

(ii) If $\oa{x\xi_+}\in\st(\log_x\tau_+)$, then $\xi_+\in\st(\tau_+)$.
If $\oa{x\xi_+}\in\ost(\log_x\tau_+)$, then $\xi_+\in\ost(\tau_+)$.
\end{lem}
\proof
(ii) The direction $\oa{x\xi_+}$ has an antipode $v_-\in\st(\log_x\tau_-)$.
By Lemma~\ref{lem:starsontostars},
$v_-$ is the initial direction $v_-=\oa{x\xi_-}$ of a ray $x\xi_-\subset P$ with $\xi_-\in\st(\tau_-)$.
Since $\xi_{\pm}$ are antipodes, 
Lemma~\ref{lem:antipbdpar}
implies that $\xi_+\in\st(\tau_+)$. 
If $\oa{x\xi_+}\in\ost(\log_x\tau_+)$, then $\xi_+\in\ost(\tau_+)$ 
because $\theta(\oa{x\xi_+})=\theta(\xi_+)$.

(i)
The corresponding assertions for $xy$ follow, 
because segments in $P$ extend to rays in $P$.
Moreover, as in the proof of (ii),
if $\oa{xy}\in\st(\log_x\tau_+)$, 
then $xy$ is contained in a ray $\xi_-y$ with $\xi_-\in\st(\tau_-)$
and hence $\oa{yx}\in\st(\log_y\tau_-)$.
\qed

\begin{rem}
Part (ii) of the last lemma yields a partial converse to Lemma~\ref{lem:antipbdpar};
it implies:
$$\log_x^{-1}(\st(\log_x\tau_+))\cap\geo P=\st(\tau_+)$$
\end{rem}

The lemma
motivates the following notion:
\begin{dfn}[Longitudinal directions and segments in parallel sets]
\label{dfn:longipar}
At a point $x\in P$,
the directions in $\ost(\log_x\tau_+)$ are called {\em longitudinal}
and the directions in $\ost(\log_x\tau_-)$ {\em anti-longitudinal}.
Moreover, 
$\Theta$-regular (anti-)longitudinal directions are called {\em $\Theta$-(anti-)longitudinal}.
A nondegenerate segment $xy\subset P$ is called {\em ($\Theta$-)(anti-)longitudinal}
if $\oa{xy}$ has this property.
\end{dfn}
\begin{rem}
(i) Longitudinal directions and segments are in particular $\taumod$-regular.

(ii) A direction is anti-longitudinal if and only if some, equivalently, all opposite directions tangent to $P$ 
are longitudinal.

(iii) A nondegenerate segment is ($\Theta$-)(anti-)longitudinal,  if and only if all nondegenerate subsegments are. 
\end{rem}
We make analogous definitions for paths:
\begin{dfn}[Longitudinal paths in parallel sets]
\label{dfn:longipp}
We say that a path $c:I\to P$ is {\em ($\Theta$-)(anti-)longitudinal} 
if all segments $c(t_1)c(t_2)$ for $t_1<t_2$ have this property. 
\end{dfn}

Note that if $c:I\to P$ is longitudinal, then
$c(I\cap(t,+\infty))\subset V(c(t),\ost(\tau_+))$
and $c(I\cap(-\infty,t))\subset V(c(t),\ost(\tau_-))$ 
for $t\in I$.

Longitudinal paths are, up to reparametrization, bilipschitz;
they become bilipschitz when parametrized by arc length:
\begin{lem}[Bounded detours]
\label{lem:bdddet}
There exists a constant $L=L(\Theta)\geq1$ such that for every $\Theta$-longitudinal
path $c:[a,b]\to P$ it holds that 
$L(c)\leq L(\Theta)\cdot d(c(a),c(b))$.
\end{lem}
\proof
We choose $\xi_-\in\tau_-$.
By the radius bound (\ref{eq:thetinba}) for $\Theta$, 
there exists $\rho<\pihalf$ such that for every $\Theta$-longitudinal segment $xy\subset P$
it holds that 
$b_{\xi_-}(y)-b_{\xi_-}(x)\geq d(x,y)\cdot\cos\rho$.
It follows that 
$d(c(b),c(a))\geq b_{\xi_-}(c(b))-b_{\xi_-}(c(a))\geq L(c)\cdot\cos\rho$.
\qed

\medskip
In order to be able to speak of openness and closedness of the longitudinality condition, 
we identify, as before for $X$, the space of segments in $P$
with the space $P\times P$ of pairs of points 
which is equipped with a natural topology.

\begin{lem}[Open and closed]
The subset of longitudinal segments in $P$ 
is open in the space of all segments in $P$,
and also closed in the subspace of $\taumod$-regular segments.
\end{lem}
\proof
Let $xy\subset P$ be $\Theta$-longitudinal.
Then $y$ lies in the interior of the cone $V(x,\st_{\Theta'}(\tau_+))$
and its distance from the boundary of this cone is $\geq\sin\eps_0(\Theta,\Theta')\cdot d(x,y)$,
with the constant $\eps_0$ from (\ref{eq:eps0th}).
Therefore $xy'$ is $\Theta'$-longitudinal 
if $d(y,y')<\sin\eps_0(\Theta,\Theta')\cdot d(x,y)$.
Similarly,
$x'y'$ is $\Theta''$-longitudinal 
if $d(x,x')<\sin\eps_0(\Theta',\Theta'')\cdot d(x,y')$.
Hence, 
longitudinality is an open condition for segments in $P$.

The uniform estimates show moreover, 
that if a $\taumod$-regular segment can be arbitrarily well approximated by longitudinal segments,
then it is longitudinal itself.
So, longitudinality is also a closed condition for $\taumod$-regular segments in $P$.
\qed

\begin{cor}[Longitudinality preserved under regular deformation]
\label{cor:remlongi}
A continuous family of $\taumod$-regular segments $x_sy_s\subset P$, $0\leq s\leq1$, 
which contains one longitudinal segment, consists only of longitudinal segments.
\end{cor}

As we did with regularity,
one can also coarsify the notion of longitudinality
and call a segment {\em coarsely longitudinal}
if it is oriented Hausdorff close to a longitudinal segment in the parallel set.
The notion of coarse longitudinality
then applies to segments and paths 
which are close to the parallel set 
but not necessarily contained in it. 

The observation that longitudinality is preserved under regular deformation
implies that a coarsely regular quasigeodesic close to the parallel set
must be coarsely longitudinal 
as soon as some sufficiently long chord of the projected quasigeodesic is longitudinal.
In the following lemma, we again use the notation $\bar q=\pi_P\circ q$.
\begin{lem}[Coarsely longitudinal quasigeodesic]
\label{lem:clqg}
With the constant $c=c(\Theta,\Theta')>0$ from Lemma~\ref{lem:lcoaregureg} 
the following holds:

Suppose that $q:I\to X$ is a $(\Theta,B)$-regular $(L,A)$-quasigeodesic 
such that $q(I)\subset\bar N_D(P)$.
If for some subinterval $[a',b']\subset I$ of length 
$\geq L(A+c(B+D))$
the $\Theta'$-regular segment $\bar q(a')\bar q(b')\subset P$ is longitudinal,
then the same holds also for all other such subintervals.
\end{lem}
\proof
We may assume that the quasigeodesic is continuous.
The subintervals $[a',b']\subset I$ of length $\geq L(A+c(B+D))$ form a connected (possibly empty) family.
That the corresponding segments $\bar q(a')\bar q(b')\subset P$ are $\Theta'$-regular,
is due to Lemma~\ref{lem:cregqgeolongsegmregproj}.
The assertion therefore follows from Corollary~\ref{cor:remlongi}.
\qed

\subsection{Cones}

In this section,
we consider a type $\taumod$ Weyl cone 
along with the corresponding $\Theta$-cones 
and an ambient type $\taumod$ parallel set:
$$ V_{\Theta}=V(x,\stTh(\tau_+)) \subset V=V(x,\st(\tau_+)) \subset P=P(\tau_-,\tau_+)$$
If $X$ is a symmetric space,
then $P$ is determined by $V$;
if $X$ is a euclidean building, it is not.

\begin{lem}[Open Weyl cone]
\label{lem:intcon}
For a point $x\in P$, 
the open Weyl cone $V(x,\ost(\tau_+))-\{x\}$ is 
the interior of $V(x,\st(\tau_+))$ in $P$,
and $V(x,\D\st(\tau_+))$ is its topological boundary.
\end{lem}
\proof
Let $y\in V(x,\ost(\tau_+))-\{x\}$.
Then $\oa{xy}\in\ost(\log_x\tau_+)$.
Since $\ost(\log_x\tau_+)$ is open in $\Si_xX$,
it follows that also $\oa{xy'}\in\ost(\log_x\tau_+)$ for every point $y'\in P$ sufficiently close to $y$,
and Lemma~\ref{lem:extintcon} implies that $y'\in V(x,\ost(\tau_+))$.

Vice versa,
suppose that $y$ lies in the interior of $V(x,\st(\tau_+))$ with respect to $P$,
and let $F\subset P$ be a maximal flat through $x$ and $y$.
(Such a flat exists because also $P$ is a euclidean building.)
Then $y$ lies in the interior, with respect to $F$, 
of the finite union of euclidean Weyl chambers
$F\cap V(x,\st(\tau_+))=V(x,\st(\tau_+)\cap\geo F)$,
and it follows that $y\in V(x,\ost(\tau_+))-\{x\}$.
\qed

\medskip
A basic property of Weyl and $\Theta$-cones is their convexity.
It is deduced from the convexity of stars at infinity (Proposition~\ref{prop:stconv}):
\begin{prop}[Convexity of cones]
\label{prop:thconeconv}
The Weyl cone $V(x,\st(\tau_+))$,
the open Weyl cone $V(x,\ost(\tau_+))-\{x\}$ 
and the $\Theta$-cones $V(x,\stTh(\tau_+))$ 
are convex subsets of $X$. 

More precisely, in the Weyl cone case,
$V(x,\st(\tau_+))$ is the intersection of the parallel set $P(\tau_-,\tau_+)$
and the root type horoballs which are centered at $\geo P(\tau_-,\tau_+)$,
contain $x$ in their boundary and $\st(\tau_+)$ in their visual boundary,
and $V(x,\ost(\tau_+))-\{x\}$ is the intersection of the parallel set and the open horoballs.
\end{prop}
\proof
The assertions for the closed cones 
have been proven in \cite[Props. 2.14 and  2.18]{morse} 
in the case of symmetric spaces; 
the proofs for euclidean buildings are identical and we will omit them. 
The assertion for the open Weyl cone follows from Lemma~\ref{lem:intcon}, 
because it is the interior of $V$ inside $P$.
The open Weyl cone 
is therefore contained in the interior of every horoball $Hb_{\zeta,x}$
for which, at infinity, 
$\st(\tau_+)\subset\bar B(\zeta,\pihalf)$. 
\qed

\medskip
As a consequence, one obtains, compare \cite[Cor 2.19]{morse}:
\begin{cor}[Nested cones]
\label{cor:nestcone}
(i) If $x'\in V(x,\st(\tau_+))$, then $V(x',\st(\tau_+))\subset V(x,\st(\tau_+))$
and $V(x',\ost(\tau_+))-\{x'\}\subset V(x,\st(\tau_+))-\{x\}$.

(ii) If $x'\in V(x,\ost(\tau_+))-\{x\}$, then $V(x',\st(\tau_+))\subset V(x,\ost(\tau_+))$.

(iii) If $x'\in V(x,\st_{\Theta}(\tau_+))$, then 
$V(x',\st_{\Theta}(\tau_+))\subset V(x,\st_{\Theta}(\tau_+))$.   
\end{cor}

Longitudinality in the Weyl cone can be defined independently of the ambient parallel set:
\begin{dfn}[Longitudinal directions in Weyl cones]
\label{dfn:longicon}
At a point $y\in V$,
the directions in $\ost(\log_y\tau_+)$ are called {\em longitudinal}
and the directions opposite to them {\em anti-longitudinal}.
\end{dfn}
As before in the case of parallel sets, see Definition~\ref{dfn:longipar},
we call $\Theta$-regular (anti-)longitudinal directions {\em $\Theta$-(anti-)longitudinal},
and we call a nondegenerate segment {\em ($\Theta$-)(anti-)longitudinal}
if its (initial) direction has this property.
Moreover, we define {\em longitudinal paths} in Weyl cones 
as in the parallel set case, cf. Definition~\ref{dfn:longipp}.

Note that tangent directions to the Weyl cone $V$ and segments in it are longitudinal in the Weyl cone
if and only if they are longitudinal in the ambient parallel set $P$.

\medskip
We next describe the anti-longitudinal directions.
\begin{lem}
\label{lem:spdircon}
(i) If $y\in V(x,\D\st(\tau_+))$, then $\ost(\log_y\tau_-)\cap\Si_yV=\emptyset$.

(ii) If $y\in V(x,\ost(\tau_+))-\{x\}$, then 
$\st(\log_y\tau_-)=\st(\tau_-(yx))\subset\Si_yV$.
\end{lem}
\proof
(i) Suppose that $\ost(\log_y\tau_-)\cap\Si_yV\neq\emptyset$.
Since $\ost(\log_y\tau_-)$ is open,
it must contain a direction 
which is represented by a segment in $V$,
i.e. there exists 
$y\neq z\in V\cap V(y,\ost(\tau_-))$. 
Hence $y\in V(z,\ost(\tau_+))$, 
and Corollary~\ref{cor:nestcone} yields that 
$y\in V(x,\ost(\tau_+))-\{x\}=V-V(x,\D\st(\tau_+))$,
which shows the first assertion. 

(ii) If $y\in V(x,\ost(\tau_+))-\{x\}$,
then $xy$ is longitudinal and $\tau_-(yx)=\log_y\tau_-$.
According to Lemma~\ref{lem:intcon}, 
$y$ lies in the interior of $V$ with respect to $P$,
so $\st(\log_y\tau_-)\subset\Si_yP=\Si_yV$.
\qed

\begin{cor}[Anti-longitudinal directions in Weyl cones]
\label{cor:spdirdiamo1}
Anti-longitudinal directions in $y\in V$ exist if and only if $y\in V(x,\ost(\tau_+))-\{x\}$.
In this case,
the set of anti-longitudinal directions in $y$ equals
$\ost(\tau_-(yx))=\ost(\log_y\tau_-)$. 
\end{cor}
\proof
The space of anti-longitudinal directions in $\Si_yV$ equals $\ost(\log_y\tau_-)\cap\Si_yV$
and is, by its definition, independent of $P$.
The assertion therefore follows from the lemma.
\qed

\subsection{Longitudinal convexity of parallel sets}
\label{sec:strbrokpbuil}

This section prepares the discussion of diamonds.
We consider broken geodesic paths in $X$. 
\begin{dfn}[Straight path]
We say that a broken geodesic path $x_0x_1\dots x_k$ in $X$ is {\em $\taumod$-straight} 
if it is piecewise $\taumod$-regular and if at any vertex $x_i$ for $0<i<k$ 
the $\taumod^{\pm}$-directions of the adjacent segments are opposite, 
i.e. if the simplices $\tau_{\pm}(x_ix_{i\pm1})\subset\Si_{x_i}X$ are opposite.
We call the path {\em $\Theta$-straight}
if in addition it is piecewise $\Theta$-regular. 
\end{dfn}
Note that 
if the directions of the adjacent segments themselves are opposite, 
$\angle_{x_i}(x_{i-1},x_{i+1})=\pi$, 
then the broken geodesic path is geodesic.

The definitions carry over to semi- and biinfinite broken geodesic paths 
$x_0\dots x_k\xi_+$, $\xi_-x_0\dots x_k$ and $\xi_-x_0\dots x_k\xi_+$
for $\xi_{\pm}\in\geo X$.
A finite $\taumod$-straight path $x_0x_1\dots x_k$ 
can always be extended to a biinfinite $\taumod$-straight path $\xi_-x_0\dots x_k\xi_+$
with $\taumod^{\pm}$-regular ideal endpoints $\xi_{\pm}$.
\begin{dfn}[Longitudinal path]
We call a broken geodesic path $x_0x_1\dots x_k$ in a parallel set
{\em ($\Theta$-)(anti-)longitudinal} if all subsegments $x_{i-1}x_i$ have this property.
\end{dfn}
Longitudinal paths in parallel sets are clearly straight.
The next result shows that, conversely, straight paths are longitudinal paths in parallel sets.
This is clear when $X$ is a symmetric space and the parallel set is uniquely determined,
but requires an argument when $X$ is a euclidean building. 
\begin{prop}
\label{prop:longiconvpar}
Each semi-infinite $\Theta$-straight path $x_0x_1\dots x_k\xi_+$ 
is contained in the $\Theta$-cone $V(x_0,\stTh(\tau(\xi_+)))$.
For each biinfinite $\Theta$-straight path $\xi_-x_0x_1\dots x_k\xi_+$,
the simplices $\tau_{\pm}(\xi_{\pm})\subset\geo X$ are opposite, 
the path is contained in the parallel set $P(\tau_-(\xi_-),\tau_+(\xi_+))$
and all segments $x_ix_j$ for $i<j$ are $\Theta$-longitudinal. 
\end{prop}
\proof
Consider a $\taumod$-straight path $x_0x_1\xi_+$.
By straightness, 
the direction $\oa{x_1x_0}$ has an antipode $v_+$ such that $\tau(v_+)=\tau(x_1\xi_+)$.
Using Lemma~\ref{lem:starsontostars},
we can extend $x_0x_1$ to a ray $x_0x_1\eta_+$ such that $\tau(\eta_+)=\tau(\xi_+)$. 
It follows that $x_1\in V(x_0,\stTh(\tau(\eta_+)))=V(x_0,\stTh(\tau(\xi_+)))$.
The assertion for semi-infinite paths follows by induction
using the nestedness of cones (Corollary~\ref{cor:nestcone}).

Consider now a biinfinite path $\xi_-x_0x_1\dots x_k\xi_+$.
From the semi-infinite case we know that 
$x_0x_1\dots x_k\xi_+\subset V(x_0,\stTh(\tau(\xi_+)))$.
In particular,
$\tau(x_0\xi_+)=\tau(x_0x_1)$,
and hence the simplices $\tau_{\pm}(\xi_{\pm})$ are $x_0$-opposite. 
It follows that $x_0\in P=P(\tau_-(\xi_-),\tau_+(\xi_+))$ and furthermore that $x_0x_1\dots x_k\subset P$.
The longitudinality follows from the semi-infinite case.
\qed

\medskip
Since longitudinal paths are not only piecewise regular, but globally regular,
the proposition can be understood as a 
{\em local-to-global} principle for the regularity of broken geodesic paths:
\begin{cor}
\label{cor:loglobgp}
Suppose that the path $x_0x_1\dots x_k$ is $\Theta$-straight.
Then all segments $x_ix_j$ for $i<j$ are $\Theta$-regular,
and for $i<j<k$ it holds that 
$\tau_+(x_ix_j)=\tau_+(x_ix_k)\subset\Si_{x_i}X$
and $\tau_-(x_kx_j)=\tau_-(x_kx_i)\subset\Si_{x_k}X$.
\end{cor}
\proof
We extend the path to a biinfinite $\Theta$-straight path and then apply the proposition.
\qed

\medskip
We next observe an extension of the convexity property for parallel sets.
That parallel sets are convex means, by definition, 
that a geodesic segment is contained in the parallel set if its endpoints are.
This remains true for straight broken geodesic paths 
whose pair of endpoints in the parallel set is longitudinal:
\begin{cor}[Longitudinal convexity of parallel sets]
\label{cor:longconvpar}
Let $x_0x_1\dots x_k$ be a $\Theta$-straight path
with endpoints in the parallel set $P=P(\tau_-,\tau_+)$ 
and suppose that the segment $x_0x_k\subset P$ is longitudinal.
Then $x_0x_1\dots x_k\subset V(x_0,\stThp(\tau_+))\cap V(x_k,\stThm(\tau_-))\subset P$.
\end{cor}
\proof
Let $\xi_{\pm}\in\inte(\tau_{\pm})$.
By assumption, the broken path $\xi_-x_0x_k\xi_+$ is then longitudinal in $P$.
Since $\tau_+(x_0x_1)=\tau_+(x_0x_k)$ and $\tau_-(x_kx_{k-1})=\tau_-(x_kx_0)$
by Corollary~\ref{cor:loglobgp},
the biinfinite path $\xi_-x_0x_1\dots x_k\xi_+$ is also $\Theta$-straight.
Proposition~\ref{prop:longiconvpar} yields that the path 
$x_0x_1\dots x_k$ is contained in $P$
and, more precisely, in the cones $V(x_0,\stThp(\tau_+))$ and $V(x_k,\stThm(\tau_-))$.
\qed

\medskip
We push the longitudinal convexity property slightly further for once broken paths $x_-yx_+$,
replacing the open assumption of straightness by a closed condition.
\begin{lem}
Let $x_-x_+\xi_+$ be a $\taumod$-straight broken path.
Suppose that $y\in X-\{x_-,x_+\}$ and there exists a pair of opposite type $\taumod^{\pm}$
simplices $\tau_y^{\pm}\subset\Si_yX$
such that $\oa{yx_{\pm}}\in\st(\tau_y^{\pm})$.
Then 
$$y\in V(x_-,\st(\tau(\xi_+))).$$
\end{lem}
\proof
Again, the assertion (and the following argument) is trivial if $X$ is a symmetric space.

We first look for simplices $\tau'_{\pm}\subset\geo X$
such that $\tau_y^{\pm}=\log_y\tau'_{\pm}$
and $yx_{\pm}\subset V(y,\st(\tau'_{\pm}))$.
To find them,
we extend the segments $yx_{\pm}$ to rays $y\xi'_{\pm}$
and let $\si_y^{\pm}\supset\tau_y^{\pm}$ be chambers in $\Si_yX$
containing the directions $\oa{yx_{\pm}}$.
According to Lemma~\ref{lem:starsontostars},
there exist chambers $\si'_{\pm}\subset\geo X$
such that 
$\si_y^{\pm}=\log_y\si'_{\pm}$
and $\xi'_{\pm}\in\si'_{\pm}$. 
Then their type $\taumod^{\pm}$ faces $\tau'_{\pm}\subset\si'_{\pm}$ have the desired properties.
Moreover, 
the simplices $\tau'_{\pm}$ are $y$-opposite, because the simplices $\tau_y^{\pm}$ are opposite.
It follows that $x_-yx_+\subset P'=P(\tau'_-,\tau'_+)$.

To see that $x_-x_+$ is longitudinal in $P'$,
we note that $y\in V(x_{\mp},\st(\tau'_{\pm}))$.
By Corollary~\ref{cor:nestcone},
there are the triples of nested cones
$V(x_{\mp},\st(\tau'_{\pm}))\supset V(y,\st(\tau'_{\pm}))\supset V(x_{\pm},\st(\tau'_{\pm}))$.
So, $x_{\pm}\in V(x_{\mp},\st(\tau'_{\pm}))$. 
Since the segment $x_-x_+$ is $\taumod$-regular,
it follows that even $x_{\pm}\in V(x_{\mp},\ost(\tau'_{\pm}))$
and $x_-x_+$ is longitudinal in $P'$. 

Now $\xi_+$ comes in and we show that $\tau'_+$ can be replaced by $\tau(\xi_+)$.
The straightness of $x_-x_+\xi_+$ implies that 
the pair of simplices $(\tau'_-,\tau(\xi_+))$ is $x_+$-opposite. 
Hence $x_-yx_+\subset V(x_+,\st(\tau'_-))\subset P(\tau'_-,\tau(\xi_+))$.
Since $x_-\in V(y,\st(\tau'_-))$, it follows that $y\in V(x_-,\st(\tau(\xi_+)))$.
\qed

\begin{cor}
\label{cor:wloconvbi}
Let $\xi_-x_-x_+\xi_+$ be a $\taumod$-straight broken path.
Suppose that $y\in X-\{x_-,x_+\}$ and there exists a pair of opposite type $\taumod^{\pm}$
simplices $\tau_y^{\pm}\subset\Si_yX$
such that $\oa{yx_{\pm}}\in\st(\tau_y^{\pm})$.
Then
$$y\in V(x_-,\st(\tau_+(\xi_+)))\cap V(x_+,\st(\tau_-(\xi_-)))\subset P(\tau_-(\xi_-),\tau_+(\xi_+)).$$
\end{cor}

Sometimes the following terminology extending Definition~\ref{def:xopp} will be convenient:
\begin{dfn}
[$(x_-,x_+)$-opposite]
\label{def:xyopp}
For a $\taumod$-regular segment $x_-x_+\subset X$,
we say that a pair $(\tau_-,\tau_+)$ of opposite simplices $\tau_{\pm}\subset\geo X$ 
is {\em $(x_-,x_+)$-opposite}
if the pairs of simplices $(\log_{x_{\pm}}\tau_{\pm},\tau_{\mp}(x_{\pm}x_{\mp}))$
are opposite
(for both choices of signs).
\end{dfn}

\begin{lem}
Two opposite simplices $\tau_{\pm}\subset\geo X$ are  $(x_-,x_+)$-opposite if and only if 
$x_-x_+$ is a longitudinal segment in the parallel set $P(\tau_-,\tau_+)$. 
\end{lem}
\proof 
This follows from the fact that straight broken paths are contained in parallel sets 
as longitudinal paths, cf. Proposition~\ref{prop:longiconvpar}.
\qed

\subsection{Diamonds and Weyl hulls}
\label{sec:diamo}

We define diamonds independently of ambient parallel sets:
\begin{dfn}[Diamond]
For a $\taumod$-regular segment $x_-x_+\subset X$,
the {\em $\taumod$-diamond} 
$$\diamot(x_-,x_+)\subset X$$
is the subset consisting of $x_{\pm}$ and all points $y\in X-\{x_-,x_+\}$
such that  $\oa{yx_{\pm}}\in\st(\tau_y^{\pm})$
for some pair of opposite type $\taumod^{\pm}$ simplices $\tau_y^{\pm}\subset\Si_yX$.
\end{dfn}
Longitudinal convexity implies that diamonds are contained in parallel sets 
and yields the following description:
\begin{lem}
\label{lem:diamoinpar}
For any pair $(\tau_-,\tau_+)$ of $(x_-,x_+)$-opposite type $\taumod^{\pm}$ simplices 
$\tau_{\pm}\subset\geo X$,
it holds that 
\begin{equation*}
\diamot(x_-,x_+)=V(x_-,\st(\tau_+))\cap V(x_+,\st(\tau_-))\subset P(\tau_-,\tau_+).
\end{equation*}
\end{lem}
\proof
That the diamond is contained in the intersection of Weyl cones,
follows immediately from Corollary~\ref{cor:wloconvbi}.
The reverse inclusion is clear.
\qed

\medskip
We will refer to $V(x_{\mp},\st(\tau_{\pm}))$ as {\em ambient Weyl cones} 
and to $P(\tau_-,\tau_+)$ as an {\em ambient parallel set} for the diamond.
Again, these are unique if $X$ is a symmetric space, but not if it is a euclidean building.  

It follows in particular that diamonds are {\em convex}.

Around their tips,
diamonds coincide up to a uniform radius with Weyl cones. 
With the constant $\eps_0(\Theta)$ from (\ref{eq:eps0}), we have:
\begin{lem}[Conical around tips]
\label{lem:diamoco}
(i) $\Si_{x_{\pm}}\diamot(x_-,x_+)=\st(\tau_{\mp}(x_{\pm}x_{\mp}))$.

(ii) If $x_-x_+$ is $\Theta$-regular,
then every segment in $\diamot(x_-,x_+)$ with initial point $x_{\pm}$
uniquely extends to a segment in $\diamot(x_-,x_+)$ 
with length $\sin\eps_0(\Theta)d(x_-,x_+)$.
\end{lem}
\proof
Let $(\tau_-,\tau_+)$ be an $(x_-,x_+)$-opposite pair of simplices. 
Since $\diamot(x_-,x_+)$ coincides with $V(x_{\pm},\st(\tau_{\mp}))$ near $x_{\pm}$,
we have that 
$\Si_{x_{\pm}}\diamot(x_-,x_+)=\Si_{x_{\pm}}V(x_{\pm},\st(\tau_{\mp}))=\st(\tau_{\mp}(x_{\pm}x_{\mp}))$.

By triangle comparison, 
$x_{\pm}$ has distance $\geq\sin\eps_0(\Theta)d(x_-,x_+)$ from $\D V(x_{\mp},\st(\tau_{\pm}))$.
It follows that 
$B(x_{\pm},\sin\eps_0(\Theta)d(x_-,x_+))\cap P(\tau_-,\tau_+)\subset V(x_{\mp},\st(\tau_{\pm}))$.
Intersecting with $V(x_{\pm},\st(\tau_{\mp}))$ yields the assertion.
\qed

\medskip
As a consequence of Lemma~\ref{lem:intcon}, 
the {\em interior} of the diamond $\diamo=\diamot(x_-,x_+)$ 
with respect to an ambient parallel set $P(\tau_-,\tau_+)$ is given by
\begin{equation}
\label{eq:intdiamo}
\inte(\diamot(x_-,x_+))=\bigl(V(x_-,\ost(\tau_+))\cap V(x_+,\ost(\tau_-))\bigr) -\{x_-,x_+\}.
\end{equation}
Note that the interior is always nonempty.
For instance, the interior points of the $\taumod$-regular segment $x_-x_+$ belong to it.

For a $\Theta$-regular segment $x_-x_+$,
we define the {\em $\Theta$-diamond} 
$$\diamotTh(x_-,x_+)\subset\diamot(x_-,x_+)$$
in a similar way 
as the subset consisting of $x_{\pm}$ and all points $y\in X-\{x_-,x_+\}$,
for which the segments $yx_{\pm}$ are $\Theta_{\pm}$-re\-gu\-lar
with opposite $\taumod^{\pm}$-directions $\tau_{\pm}(yx_{\pm})$ at $y$.
It follows from Lemma~\ref{lem:diamoinpar} that 
$$\diamotTh(x_-,x_+)=V(x_-,\st_{\Theta}(\tau_+))\cap V(x_+,\st_{\Theta}(\tau_-)).$$
We will need the following semicontinuity property of diamonds:
\begin{lem}[Semicontinuity]
\label{lem:semicontdiamo}
Suppose that the diamond 
$\diamot(x_-,x_+)$ intersects the open subset $O\subset X$.
Then for all pairs of points $(x'_-,x'_+)$ sufficiently close to 
$(x_-,x_+)$,
the diamond $\diamot(x'_-,x'_+)$ still intersects $O$.
\end{lem}
\proof
Suppose first that $X$ is a euclidean building. 
By assumption,
there exists a point 
$y\in O\cap\inte(\diamot(x_-,x_+))$.
The segments $yx_{\pm}$ are $\taumod^{\pm}$-regular.
Therefore, 
if the points $x'_{\pm}$ are sufficiently close to $x_{\pm}$,
then also the segments $yx'_{\pm}$ are $\taumod^{\pm}$-regular
and $\tau_{\pm}(yx'_{\pm})=\tau_{\pm}(yx_{\pm})$.
In particular, the simplices $\tau_{\pm}(yx'_{\pm})$ are opposite 
and it follows that $y\in\diamot(x'_-,x'_+)$.

Suppose now that $X$ is a symmetric space.
Let $P$ denote the unique ambient parallel set of the diamond $\diamot(x_-,x_+)$.
The argument used in the euclidean building case
now proves the assertion for pairs of points in the parallel set,
i.e. there exists $\eps>0$ such that the assertion holds for the pair $(x'_-,x'_+)$ 
if $x'_{\pm}\in B(x_{\pm},\eps)\cap P$.
Using the action of the isometry group $G=\Isom(X)$,
it follows furthermore 
that there exists a neighborhood $U$ of the neutral element in $G$,
such that the assertion holds for all pairs $g\cdot(x'_-,x'_+)$ 
with $g\in U$ and $x'_{\pm}\in B(x_{\pm},\eps)\cap P$.
Since $G$ acts transitively on type $\taumod$ parallel sets,
this finishes the proof.
\qed

\medskip
We will prove later the stronger property
that diamonds depend {\em continuously} on their pair of tips,
see Proposition~\ref{prop:contdiamo} below.

\medskip
In order to define longitudinal directions in diamonds, 
we observe that, whether a direction is longitudinal with respect to an ambient parallel set,
does not depend on the ambient parallel set:
\begin{lem}
\label{lem:spdirdiamo}
Let $y\in\diamo$. 
If the segment $yx_{\pm}$ is $\taumod^{\pm}$-regular, 
then $\st(\log_y\tau_{\pm})=\st(\tau_{\pm}(yx_{\pm}))\subset\Si_y\diamo$.
Otherwise, 
$\ost(\log_y\tau_{\pm})\cap\Si_y\diamo=\emptyset$.
\end{lem}
\proof
The segment $yx_+$ is $\taumod^+$-regular
if and only if $y\in V(x_+,\ost(\tau_-))-\{x_+\}$.

Thus, if $yx_+$ is not $\taumod^+$-regular,
then $\ost(\log_y\tau_+)$ is disjoint from 
$\Si_yV(x_+,\st(\tau_-))\supset\Si_y\diamo$ 
by Lemma~\ref{lem:spdircon}.
On the other hand, 
if $yx_+$ is $\taumod^+$-regular,
then $\diamo=V(x_-,\st(\tau_+))$ near $y$,
and Lemma~\ref{lem:spdircon} yields that $\Si_y\diamo$ contains $\st(\log_y\tau_+)$.
Moreover, $\log_y\tau_+=\tau_+(yx_+)$.
\qed

\begin{cor}
\label{cor:spdirdiamo2}
The intersection
$\Si_y\diamo\cap\ost(\log_y\tau_{\pm})$
does not depend on the ambient parallel set $P(\tau_-,\tau_+)$.
It is nonempty if and only if the segment $yx_{\pm}$ is $\taumod^{\pm}$-regular,
and then equal to $\ost(\log_y\tau_{\pm})=\ost(\tau_{\pm}(yx_{\pm}))$.
\end{cor}
This justifies:
\begin{dfn}[Longitudinal directions in diamonds]
\label{dfn:longidiamo}
In a point $y\in\diamo$,
we call the directions in $\Si_y\diamo\cap\ost(\log_y\tau_+)$ {\em longitudinal}
and the directions in $\Si_y\diamo\cap\ost(\log_y\tau_-)$ {\em anti-lon\-gi\-tu\-di\-nal}.
\end{dfn}
As before in the case of parallel sets and Weyl cones,
we call $\Theta$-regular (anti-)longitudinal directions {\em $\Theta$-(anti-)longitudinal},
and we call a nondegenerate segment {\em ($\Theta$-)(anti-)longitudinal}
if its (initial) direction has this property.
Note that the segment $x_-x_+$ is longitudinal.

Our discussion shows that 
directions and segments in the diamond are longitudinal 
if and only if they are longitudinal in an ambient parallel set.

\medskip
Based on the notion of longitudinality, we can now state:
\begin{lem}[Nested diamonds]
\label{lem:nestdiamo}
If $x'_-x'_+\subset\diamo$ is longitudinal, 
then $\diamot(x'_-,x'_+)\subset\diamo$.

If the segment $x'_-x'_+\subset\diamoTh(x_-,x_+)$ is $\Theta$-longitudinal, 
then $\diamoTh(x'_-,x'_+)\subset\diamoTh(x_-,x_+)$.
\end{lem}
\proof
This is a direct consequence of the nestedness of cones,
cf. Corollary~\ref{cor:nestcone}.
\qed

\medskip
We can also reformulate 
the longitudinal convexity property of parallel sets, 
cf. Corollary~\ref{cor:longconvpar},
for diamonds:
\begin{cor}[Longitudinal convexity of diamonds]
\label{cor:longiconvdiamo}
Each $\Theta$-straight broken geodesic path $x_0x_1\dots x_k$ 
is contained in the $\Theta$-diamond $\diamoTh(x_0,x_k)$,
and all segments $x_ix_j$ for $i<j$ are $\Theta$-longitudinal. 
\end{cor}

We turn to the discussion of Weyl hulls of segments.

Weyl hulls are analogs of diamonds inside singular flats.
We also define them intrinsically without reference to ambient flats:
\begin{dfn}[Weyl hull]
The {\em Weyl hull}
of a nondegenerate segment $x_-x_+\subset X$ 
with type $\theta(x_-x_+)\in\inte(\taumod)$
is the subset 
$$Q(x_-,x_+)\subset X$$
consisting of $x_{\pm}$ and all points $y\in X-\{x_-,x_+\}$
such that  $\oa{yx_{\pm}}\in\tau_y^{\pm}$
for some pair of opposite type $\taumod^{\pm}$ simplices $\tau_y^{\pm}\subset\Si_yX$.
\end{dfn}
Clearly, 
$Q(x_-,x_+)\subset\diamot(x_-,x_+)$.

Applying the description of diamonds, 
it follows that Weyl hulls are cross sections of diamonds by singular flats.
Indeed,
let $(\tau_-,\tau_+)$ be a pair of $(x_-,x_+)$-opposite type $\taumod^{\pm}$ simplices 
$\tau_{\pm}\subset\geo X$.
Then the segment $x_-x_+$ is contained in a singular flat $f$
with ideal boundary sphere $\geo f=s(\tau_-,\tau_+)$,
and we obtain:
\begin{lem}
\label{lem:whinfl}
$Q(x_-,x_+)=\diamot(x_-,x_+)\cap f=V(x_-,\tau_+)\cap V(x_+,\tau_-)\subset f$.
\end{lem}
\proof
Let $y\in Q(x_-,x_+)-\{x_-,x_+\}$.
Then $\theta(x_{\pm}y)\in\taumod^{\mp}$.
In view of $Q(x_-,x_+)\subset\diamot(x_-,x_+)$,
Lemma~\ref{lem:diamoinpar} implies that 
$y\in V(x_{\pm},\tau_{\mp})=V(x_{\pm},\st(\tau_{\mp}))\cap f$.
Conversely,
the intersection of the sectors $V(x_{\pm},\tau_{\mp})$ belongs to $Q(x_-,x_+)$.
\qed

\medskip
We will refer to $V(x_{\mp},\tau_{\pm})$ as {\em ambient Weyl sectors} 
and to $f$ as an {\em minimal ambient flat} for the Weyl hull.
These are unique if $X$ is a symmetric space, but not if it is a euclidean building.  

It follows that Weyl hulls are {\em flat parallelepipeds}.
In particular that Weyl hulls are {\em convex}.
\begin{rem}
Weyl hulls can in fact be characterized by these properties:
One can show that 
$Q(x_-,x_+)$ is the smallest closed convex subset of $X$
which contains the segment $x_-x_+$ and has the property that all spaces of directions are subcomplexes.
\end{rem}

We have the following estimate 
for the size of conical neighborhoods around the tips:
\begin{lem}[Conical around tips]
\label{lem:whcontsect}
(i) $\Si_{x_{\pm}}Q(x_-,x_+)=\tau_{\mp}(x_{\pm}x_{\mp})$.

(ii)
If the segment $x_-x_+$ is $\Theta$-regular,
$\theta(x_-x_+)\in\Theta\cap\taumod\subset\inte(\taumod)$, 
then the intersection 
$Q(x_-,x_+)\cap\bar B(x_{\pm},d(x_-,x_+)\cdot\sin\eps_0(\Theta))$
is a flat cone of height $d(x_-,x_+)\cdot\sin\eps_0(\Theta)$ with tip $x_{\pm}$.
\end{lem}
\proof
This is a consequence of Lemma~\ref{lem:diamoco}.

Alternatively, one can prove this lemma analogously to Lemma~\ref{lem:diamoco} also directly
using Lemma~\ref{lem:whinfl}:
The point $x_{\pm}$ has distance $\geq d(x_-,x_+)\cdot\sin\eps_0(\Theta)$ 
from $\D V(x_{\mp},\st(\tau_{\pm}))$.
It follows that 
$\bar B(x_{\pm},d(x_-,x_+)\cdot\sin\eps_0(\Theta))\cap f\subset V(x_{\mp},\tau_{\pm})$.
Intersecting with $V(x_{\pm},\tau_{\mp})$ yields the assertion.
\qed

\medskip
We define the Weyl hull of a degenerate segment $xx$ as the one point subset $Q(x,x)=\{x\}$.

\subsection{Rays longitudinally approaching parallel sets}

We consider now geodesic rays which are {\em longitudinally asymptotic} to parallel sets
and show that they must approach and, in euclidean buildings, enter the parallel set at a uniform rate. 

Let $P=P(\tau_-,\tau_+)\subset X$ be a type $\taumod$ parallel set.

Suppose first that the model space $X$ is a symmetric space.
In this case, longitudinally asymptotic rays do not enter parallel sets, 
but approach them at a uniform exponential rate.
We will only prove the weaker statement sufficient for our purposes,
that they do so at some uniform rate.
\begin{lem}
For $d>0$ there exists a constant $\de=\de(d,\Theta)>0$ such that the following holds:

If $x\in X$ with $d(x,P)\geq d$ and $\xi\in\stTh(\tau_+)$,
then $\angle_x(\pi_P(x),\xi)\leq\pihalf-\de$.
\end{lem}
\proof
We denote $\bar x=\pi_P(x)$.
If $x'\in x\bar x$ is the point at distance exactly $d$ from $P$, 
then $\angle_{x'}(\bar x,\xi)\geq\angle_x(\bar x,\xi)$.
We may therefore assume that $d(x,P)=d$.

Note that 
$\angle_x(\bar x,\xi)\leq\pihalf$,
because the angle sum of the ideal triangle $\De(x,\bar x,\xi)$ is $\leq\pi$.
Suppose that 
$\angle_x(\bar x,\xi)=\pihalf$.
Then the ray $x\xi$ is parallel to $P$ and extends to a geodesic line $l\not\subset P$ 
parallel to $P$.
This line is forward asymptotic to $l(+\infty)=\xi\in\ost(\tau_+)$,
and its backward ideal endpoint $l(-\infty)$ is therefore contained in $\ost(\tau_-)$, 
cf.\ Lemma~\ref{lem:antipbdpar}.
Since the singular sphere in $\tits X$ spanned by the pair of antipodes $l(\pm\infty)$
contains the simplices $\tau_{\pm}$,
and hence also the singular sphere $s(\tau_-,\tau_+)$,
it follows that $P(l)\subset P$, which is a contradiction to $l\not\subset P$.

Thus the continuous function 
$$(x,\xi)\mapsto\angle_x(\pi_P(x),\xi)$$
on $\D N_d(P)\times\stTh(\tau_+)$ takes values in the open interval $(0,\pihalf)$.
It is invariant under the stabilizer in $\Isom(X)$ of the pair of simplices $(\tau_-,\tau_+)$,
because the stabilizer preserves $P$ and $\stTh(\tau_+)$.
It acts transitively on $P$ and hence cocompactly on $\D N_d(P)\times\stTh(\tau_+)$.
It follows that the range of the function in $(0,\pihalf)$ is compact.
Furthermore, the range does not depend on the parallel set,
because all type $\taumod$ parallel sets are equivalent modulo the action of the isometry group.
\qed

\medskip
We obtain:
\begin{prop}[Rays approaching parallel sets in symmetric spaces]
\label{prop:diveparsymm}
For $d>0$ 
there exists a constant $C=C(\Theta,d)>0$ such that the following holds:

If $x\in X$ and $\xi\in\stTh(\tau_+)$,
then the points on the ray $x\xi$ with distance 
$\geq C\cdot d(x,P)$ from $x$ are contained in $\bar N_d(P)$.
\end{prop}
\proof
Let $r:[0,+\infty)\to X$ be a unit speed parametrization of the the ray $x\xi$.
Then the function $f(t):=d(r(t),P)$
is smooth with derivative 
$f'(t) = -\cos\angle_{r(t)}(\pi_P\circ r(t),\xi)$,
cf.\ (\ref{eq:derivative}).
By the previous lemma,
$f'(t) \leq-\sin\de$
as long as $f(t)\geq d$.
This yields a uniform upper bound for the entry time, linear in $d(x,P)$.
\qed

\medskip
Suppose for the rest of this section that the model space $X$ is a euclidean building.

We obtain the following version of Proposition~\ref{prop:diveparsymm},
where $\eps_0(\Theta)$ is the constant from (\ref{eq:eps0}).
The special case of maximal flats had been proven in \cite[Lemma 4.6.3]{qirigid}.
\begin{prop}[Rays diving into parallel sets in euclidean buildings]
\label{prop:diveparbuil}
Suppose that $x\in X$ and $\xi\in\stTh(\tau_+)$.
Then the ray $x\xi$ enters $P$, 
and its entry point $z$ satisfies
$$\angle_z(x,P)\geq\eps_0(\Theta)>0$$
and 
$$d(x,z)\leq (\sin\eps_0(\Theta))^{-1}\cdot d(x,P).$$
\end{prop}
\proof
We assume that $x\not\in P$ and denote $\bar x=\pi_P(x)$.

Let $y\in\bar x\xi$
and suppose that $d(x,y)\cdot\sin\eps_0(\Theta) > d(x,P)$.
Applying comparison to the triangle $\De(x,\bar x,y)$,
we can bound the angle $\angle_y(\bar x,x)$ by
\begin{equation}
\label{ineq:angbd}
d(x,y)\cdot\sin\angle_y(\bar x,x) \leq d(x,P).
\end{equation}
It follows that 
$\angle_y(\bar x,x)<\eps_0$.
This implies that $\tau_-(yx)=\tau_-(y\bar x)$ and hence $\oa{yx}\in\ost(\log_y\tau_-)$
is tangent to $P$.

Since all tangent directions to $P$ are represented by segments in $P$,
and since segments with angle zero in a euclidean building initially coincide,
cf.\ the discussion of angle rigidity in section~\ref{sec:modsp},
it follows that the segment $yx$ is initially contained in $P$.
Let $z$ denote the interior point on the segment $yx$ where it exits $P$, in other words, 
the point where the segment $xy$ enters $P$.
Then $\oa{zx}\not\in\Si_zP$,
because $zx\cap P=\{z\}$.

As a consequence of (\ref{ineq:angbd}),
given $\Theta'$,
the direction $\oa{yz}=\oa{yx}\in\Si_yP$ becomes $\Theta'_-$-longitudinal as $y\to\infty$,
and, accordingly,
$\oa{zy}\in\st_{\Theta'}(\log_z\tau_+)\subset\Si_zP$.
With Corollary~\ref{cor:unifangsusp},
it follows that the antipodal direction $\oa{zx}\not\in\Si_zP$
cannot have too small angle with $P$,
i.e. $\angle_z(x,P)\geq\eps_0(\Theta')>0$.
Applying comparison to the triangle 
$\De(x,\bar x,z)$ as above then yields a uniform estimate for the entry time of $xy$ into $P$:
$$ d(x,z)\cdot\sin\eps_0(\Theta') \leq d(x,P)$$
The segment $xy$ converges to the ray $x\xi$ as $y\to\infty$,
and the entry point subconverges to a point in $x\xi\cap P$.
This shows that the ray $x\xi$ enters $P$, and that the entry point $\hat z$ satisfies the same estimate:
$$ d(x,\hat z)\cdot\sin\eps_0(\Theta') \leq d(x,P)$$
Since this estimate holds for all $\Theta'$ (containing $\Theta$ in their interior),
we also obtain it for $\Theta$.
\qed

\begin{rem}
The longitudinality assumption (that $\xi\in\stTh(\tau_+)$) is necessary, 
in both Proposition~\ref{prop:diveparsymm} and Proposition~\ref{prop:diveparbuil},
if $\taumod\subsetneq\simod$.
Note that $\geo P$ does in general not contain the stars around the type $\taumod$ simplices 
in $\geo P$ other than $\tau_{\pm}$.
Accordingly,
there may exist $\taumod$-regular rays 
which are asymptotic to $P$, but not strongly asymptotic.
Note that such rays cannot be $\simod$-regular.
\end{rem}

We will later use different versions and consequences of the proposition.
For instance, 
we can also uniformly estimate the entry time into Weyl cones in $P$ asymptotic to $\tau_+$:
\begin{cor}[Rays diving into Weyl cones]
\label{cor:diveconbuil}
Suppose that $x\in X$, $\xi\in\st_\Theta(\tau_+)$
and $\hat x\in P$.
Then the ray $x\xi$ enters the Weyl cone $V(\hat x,\st(\tau_+))\subset P$,
and its entry point $w$ satisfies
$$d(x,w)\leq (\sin\eps_0(\Theta))^{-1}\cdot d(x,\hat x).$$
\end{cor}
\proof Let $z\in P$ be the entry point of $x\xi$ into $P$,
as given by the previous proposition.
Then
$$d(x,z)\leq(\sin\eps_0(\Theta))^{-1}\cdot d(x,\hat x).$$
Let $w'\in z\xi=x\xi\cap P$ be a point at distance 
$d(x,w')>(\sin\eps_0(\Theta))^{-1}\cdot d(x,\hat x)$.
Applying CAT(0) comparison to the triangle $\De(x,\hat x,w')$,
we get that $\angle_{w'}(z,\hat x)=\angle_{w'}(x,\hat x)<\eps_0$.
It follows that 
$\oa{w'\hat x}\in\ost(\log_{w'}\tau_-)$,
i.e. the segment $w'\hat x$ is anti-longitudinal.
Hence $w'\in V(\hat x,\st(\tau_+))$.
\qed

\medskip
The next versions of the proposition
estimate the rate at which rays move away from Weyl cones and sectors:
\begin{cor}[Rays leaving Weyl cones]
\label{cor:leavecone}
Let $\rho:[0,+\infty)\to X$ be a $\Theta$-regular unit speed ray, 
and let $t_0\geq0$ denote the time when $\rho$ exits 
the Weyl cone 
$V=V(\rho(0),\st(\tau_+))$.
Then 
$$\angle_{\rho(t_0)}(\rho(+\infty),V)\geq\eps_0(\Theta)>0$$
and
$$d(\rho(t),V)\geq(t-t_0)\sin\eps_0(\Theta)$$
for $t\geq t_0$.
\end{cor}
\proof
If $t_0=0$,
then $\dot\rho(0)\not\in\Si_{\rho(0)}V=\st(\log_{\rho(0)}\tau_+)$,
and the angle estimate holds due to $\Theta$-re\-gu\-la\-ri\-ty.
If $t_0>0$, 
we may assume that $\rho(0)\in P$.
Then 
$\dot\rho(0)\in\Si_{\rho(0)}V$
and $\rho$ can be extended by a ray in $P$ 
to a line $l:\R\to X$ backward asymptotic to $l(-\infty)\in\stTh(\tau_-)$.
Applying Proposition~\ref{prop:diveparbuil} 
to subrays of $l$ yields the angle estimate also in this case.
As before, triangle comparison based on the angle estimate yields the distance estimate 
$d(\rho(t),V)\geq d(\rho(t),P)\geq(t-t_0)\sin\eps_0(\Theta)$. 
\qed

\begin{cor}[Rays leaving Weyl sectors]
\label{cor:leavesect}
Let $\rho:[0,+\infty)\to X$ be a 
unit speed ray 
of type $\theta(\dot\rho)\equiv\bar\zeta\in\Theta\cap\taumod\subset\inte(\taumod)$, 
and let $t_0\geq0$ denote the time when $\rho$ exits 
the Weyl sector
$\check V=V(\rho(0),\tau_+)$.
Then 
$$\angle_{\rho(t_0)}(\dot\rho(t_0),\check V)\geq\eps_0(\Theta)>0$$
and
$$d(\rho(t),\check V)\geq(t-t_0)\sin\eps_0(\Theta)$$
for $t\geq t_0$.
\end{cor}
\proof
The exit direction $\dot\rho(t_0)$ is not tangent to $\check V$,
$\dot\rho(t_0)\not\in\Si_{\rho(t_0)}\check V$.
Since it has type $\theta(\dot\rho(t_0))\in\inte(\taumod)$,
it spans the simplex $\tau(\dot\rho(t_0))$, 
i.e. $\dot\rho(t_0)\in\inte(\tau(\dot\rho(t_0)))$.
It follows that 
the simplex $\tau(\dot\rho(t_0))$ is not contained in the (finite) subcomplex $\Si_{\rho(t_0)}\check V$,
equivalently, $\ost(\tau(\dot\rho(t_0)))\cap\Si_{\rho(t_0)}\check V=\emptyset$.
This yields the angle estimate, 
and triangle comparison the distance estimate.
\qed

\medskip
We apply the above estimates to show that 
Weyl cones with the same type and tip must coincide up to a certain radius,
if they are close up to a certain larger radius in some uniformly regular direction:
\begin{lem}[Initial coincidence of nearby truncated Weyl cones]
\label{lem:closediamos}
Let $r,D,R\geq0$ be constants with 
$R\sin\eps_0(\Theta)\geq r+D$.
Suppose that for simplices $\tau_+,\tau'_+\in\Flagt(\geo X)$ and a point $x\in X$
it holds that 
$$(V(x,\stTh(\tau_+))-B(x,R))\cap \bar N_D(V(x,\st(\tau'_+)))\neq\emptyset.$$
Then
$$V(x,\st(\tau_+))\cap\bar B(x,r)=V(x,\st(\tau'_+))\cap\bar B(x,r).$$
\end{lem}
\proof
Let $y\in (V(x,\stTh(\tau_+))-B(x,R))\cap \bar N_D(V(x,\st(\tau'_+)))$,
and let $\tau'_-$ be a simplex $x$-opposite to $\tau'_+$.
Then $V'=V(x,\st(\tau'_+))\subset P'=P(\tau'_-,\tau'_+)$.
Furthermore,
the segment $xy$ is $\Theta$-regular and has length $\geq R$.
Let $z\in xy$ denote the point where the segment $xy$ exits the cone $V'$,
i.e. $xy\cap V'=xz$.
Then 
$$\diamot(x,z)\subset V(x,\st(\tau_+))\cap V(x,\st(\tau'_+)).$$
Corollary~\ref{cor:leavecone} yields the estimate 
$d(z,y)\cdot\sin\eps_0\leq D$.
Hence, the $\Theta$-regular segment $xz$ has length 
$\geq R-(\sin\eps_0)^{-1}D\geq(\sin\eps_0)^{-1}r$.
By Lemma~\ref{lem:diamoco}, 
the diamond $\diamot(x,z)$ agrees 
up to radius $d(x,z)\cdot\sin\eps_0\geq r$ around its vertex $x$
with both cones 
$V(x,\st(\tau_+))$ and $V(x,\st(\tau'_+))$.
\qed

\medskip
We give a version of the last result for sectors,
namely that 
Weyl sectors with the same type and tip must coincide up to a certain radius,
if they are close up to a certain larger radius:
\begin{lem}[Initial coincidence of nearby truncated Weyl sectors]
\label{lem:sectclothcoi}
Let $r,D,R\geq0$ be constants with 
$R\sin\eps_0(\Theta)\geq r+D$.
Suppose that for simplices $\tau_+,\tau'_+\in\Flagt(\geo X)$ and a point $x\in X$
the truncated Weyl sectors
$V(x,\tau_+)\cap\bar B(x,R)$ and $V(x,\tau'_+)\cap\bar B(x,R)$ are $D$-Hausdorff close.
Then $$V(x,\tau_+)\cap\bar B(x,r)=V(x,\tau'_+)\cap\bar B(x,r).$$
\end{lem}
\proof
Consider a $\Theta$-regular unit speed ray with initial point $x$ in the Weyl sector $V(x,\tau'_+)$.
Since the ray remains in the $D$-neighborhood of the other sector $V(x,\tau_+)$ up to time $R$, 
it does not exit $V(x,\tau_+)$ before time $R-(\sin\eps_0)^{-1}D$, 
cf. Corollary~\ref{cor:leavesect}.
Therefore the intersection of the sectors $V(x,\tau_+)$ and $V(x,\tau'_+)$
contains a $\Theta$-regular segment $xz$ of length $R-(\sin\eps_0)^{-1}D$.
In view of Lemma~\ref{lem:whinfl},
it follows that the Weyl hull of $xz$ is also contained in this intersection,
$$Q(x,z)\subset V(x,\tau_+)\cap V(x,\tau'_+).$$
By Lemma~\ref{lem:whcontsect}, 
$Q(x,z)$ contains a conical neighborhood of radius 
$d(x,z)\cdot\sin\eps_0\geq r$ around its tip $x$.
So, the sectors $V(x,\tau_+)$ and $V(x,\tau'_+)$ coincide at least up to radius $r$.
\qed

\subsection{Continuity of diamonds}
\label{sec:mdiamo}

Let $X$ be again a model space.
The main result of this technical section is 
that diamonds depend continuously on their tips, 
see Proposition~\ref{prop:contdiamo}.

Let $x_-x_+\subset X$ be a $\Theta$-regular segment
and consider the $\taumod$-diamond 
$$\diamo=\diamot(x_-,x_+)$$ 
spanned by it.
Our first goal is to estimate the {\em inradius} of the diamond.

We represent the diamond as an intersection of Weyl cones in an ambient parallel set,
$$\diamo=V(x_-,\st(\tau_+))\cap V(x_+,\st(\tau_-))\subset P=P(\tau_-,\tau_+) .$$
We recall, see Proposition~\ref{prop:thconeconv},
that inside $P$ each of the Weyl cones $V(x_{\mp},\st(\tau_{\pm}))$ 
is the intersection of a certain family of horoballs,
$$
V(x_{\mp},\st(\tau_{\pm}))= P\cap\bigcap_{\zeta\in Z_{\pm}} \{b_{\zeta}\le 0\},  
$$
with centers $\zeta\in Z_{\pm}\subset\geo P$
and normalized by $b_{\zeta}(x_{\mp})=0$ for $\zeta\in Z_{\pm}$. 
Accordingly, for the convex function
$$
b=\sup_{\zeta\in Z}  b_\zeta,
$$
where $Z=Z_-\cup Z_+$,
we have 
$$
\diamo=\{b|_P\leq 0\}. 
$$
We estimate the decay of these Busemann functions $b_{\zeta}$ along the segment $x_-x_+$.
\begin{lem}
If $\zeta\in Z_{\pm}$, then 
$b_{\zeta}(x_{\mp})-b_{\zeta}(x_{\pm})
\geq d(x_-,x_+)\cdot\sin\eps_0(\Theta)$.
\end{lem}
\proof
We extend $x_-x_+$ to a $\Theta$-longitudinal line $l\subset P$.
Then $\xi_{\pm}:=l(\pm\infty)\in\stTh(\tau_{\pm})$.
If $\zeta\in Z_{\pm}$, 
then the Weyl cone $V(x_{\mp},\st(\tau_{\pm}))$ is contained in a horoball centered at $\zeta$,
and therefore $\st(\tau_{\pm})\subset\bar B(\zeta,\pihalf)$.
Hence 
$\bar B(\xi_{\pm},\eps_0(\Theta))
\subset\st(\tau_{\pm})\subset\bar B(\zeta,\pihalf)$.
It follows that $\tangle(\xi_{\pm},\zeta)\leq\pihalf-\eps_0(\Theta)$.
Thus, the Busemann function $b_{\zeta}$ has slope 
$\leq-\sin\eps_0(\Theta)$ 
along the ray $x_{\mp}\xi_{\pm}\supset x_-x_+$.
\qed

\medskip
We denote by $m$ the midpoint of the segment $x_-x_+$.
\begin{cor}[Thickness of diamonds]
\label{cor:diamothick}
$b(m)\leq-\half d(x_-,x_+)\cdot\sin\eps_0(\Theta)$
\end{cor}
\proof
Since $b_{\zeta}(x_{\mp})=0$ for $\zeta\in Z_{\pm}$,
the convexity of Busemann functions implies that 
$b_{\zeta}(m)\leq-\half d(x_-,x_+)\cdot\sin\eps_0(\Theta)$.
Taking the supremum over $Z$ yields the assertion. 
\qed

\medskip
Next, we discuss {\em product splittings} of diamonds
induced by splittings of the model space.

Suppose that the model chamber splits as a spherical join 
\begin{equation}
\label{eq:splmcha}
\simod=\simod^1\circ\simod^2,
\end{equation}
and let 
$X=X_1\times X_2$ be the corresponding product splitting of the model space,
cf. section~\ref{sec:modsp}.
If $\taumod\subset\simod^1$,
then also the $\taumod$-parallel sets, $\taumod$-Weyl cones and $\taumod$-diamonds split off 
the $X_2$-factor. 
Thus, the diamond $\diamo=\diamot(x_-,x_+)$ splits as 
$$\diamo=\diamo_1\times X_2,$$
and the cross section in the complementary factor is again a diamond, 
$\diamo_1=\diamot^{X_1}(x^-_1,x^+_1)\subset X_1$.
The segment $x^-_1x^+_1\subset X_1$ is $(\Theta\cap\simod^1)$-regular.
It is shorter than the segment $x_-x_+$, 
but of comparable length.
Indeed, the angle between the $\Theta$-regular segment $x_-x_+$ and the $X_1$-factor is
bounded above by $\diam(\Theta)\leq\pihalf-\eps_0(\Theta)$,
and hence:
\begin{equation}
\label{ineq:dicomp}
d(x_-,x_+)\cdot\sin\eps_0(\Theta) \leq d(x^-_1,x^+_1)\leq d(x_-,x_+)
\end{equation}
In the following discussion,
we fix the unique splitting (\ref{eq:splmcha})
such that $\simod^1$ is minimal with the property that it contains $\taumod$.
This includes the possibility of the trivial splitting with $\simod^1=\simod$ and $\simod^2=\emptyset$,
accordingly, $X=X_1$ and $X_2=pt$.
We note that $\simod^2\subset\D\st(\taumod)$.

In general, there is no better diameter bound than $\diam(\simod^1)\leq\pihalf$,
but we do have a uniform radius bound 
\begin{equation}
\label{ineq:betradb}
\rad(\simod^1,\cdot)\leq\rho_0=\rho_0(\Theta)<\pihalf
\end{equation}
on $\Theta_1:=\Theta\cap\simod^1$,
because otherwise $\simod^1$ would not be minimal.

\medskip
We prove now that there is a {\em uniform diameter bound} for the cross section of diamonds.

For a type $\bar\xi\in\Theta\cap\taumod$, 
we define a {\em $\bar\xi$-height} function
$$ h_{\bar\xi}:\diamo\to\R $$
as follows:
For every longitudinal (oriented) line $l_{\bar\xi}\subset P$ of type $\bar\xi$
it holds that $P=P(l_{\bar\xi})\cong l_{\bar\xi}\times CS(\geo l_{\bar\xi})$.
We define $h_{\bar\xi}$ as the restriction of a Busemann function,
$$h_{\bar\xi}:=b_{l_{\bar\xi}(-\infty)}|_{\diamo}.$$
The function $h_{\bar\xi}$ 
has the following properties,
and is determined by them up to an additive constant:
It is 1-Lipschitz, affine (i.e. affine linear along every geodesic segment),
constant on the intersections of $\diamo$ with the cross sections $pt\times CS(\geo l_{\bar\xi})$ of $P$,
and linear with slope $\equiv1$ on the intersections of $\diamo$ with the lines $l_{\bar\xi}\times pt$.
The function $h_{\bar\xi}$ is therefore independent of the ambient parallel set $P$
and well-defined up to an additive constant.

Since $\diam(\simod)\leq\pihalf$, 
we have that $\rad(\simod,\cdot)\leq\pihalf$ in particular on $\taumod$,
and hence, for every simplex $\tau\in\Flagt(\geo X)$, that $\rad(\st(\tau),\cdot)\leq\pihalf$ on $\tau$.
In particular, 
\begin{equation}
\label{eq:anglest}
\angle_{x_{\mp}}(l_{\bar\xi}(\pm\infty),\cdot)\leq\pihalf
\end{equation}
on $\diamo-\{x_{\mp}\}$,
and it follows that 
\begin{equation}
\label{eq:heiest}
h_{\bar\xi}(x_-)\leq h_{\bar\xi}\leq h_{\bar\xi}(x_+) .
\end{equation}
The estimates (\ref{eq:anglest}) and (\ref{eq:heiest}) 
improve, when restricting to the cross section of the diamond.
The angle bound (\ref{ineq:betradb}) yields the estimate
\begin{equation*}
\angle_{x^{\mp}_1}(l_{\bar\xi}(\pm\infty),\cdot)\leq\rho_0
\end{equation*}
on $\diamo_1-\{x^{\mp}_1\}$. 
(Note that the line $l_{\bar\xi}$ is parallel to $X_1$ and $l_{\bar\xi}(\pm\infty)\in\geo X_1$.)
It implies that 
\begin{equation}
\label{eq:heiest1}
h_{\bar\xi}(x^-_1) +d(x^-_1,\cdot)\cos\rho_0 \leq h_{\bar\xi}
\leq h_{\bar\xi}(x^+_1) -d(x^+_1,\cdot)\cos\rho_0
\end{equation}
on $\diamo_1$.
\begin{lem}[Diameter bound] 
\label{lem:diambcr}
The diameter of the cross section of the $\taumod$-diamond $\diamo$ is uniformly bounded by 
$$\diam(\diamo_1)\leq 2 (\cos\rho_0(\Theta))^{-1}d(x^-_1,x^+_1).$$
\end{lem}
\proof
From (\ref{eq:heiest1})
we get the radius bound 
$$\rad(\diamo_1,x^{\pm}_1)\cos\rho_0\leq h_{\bar\xi}(x^+_1)-h_{\bar\xi}(x^-_1)
\leq d(x^-_1,x^+_1)$$
and hence the diameter bound as claimed.
\qed

\medskip
We apply our discussion to prove that diamonds 
{\em depend continuously on their tips} with respect to the Hausdorff topology.
We first consider diamonds inside a fixed parallel set.

Let $\diamo\subset P$ be as above,
and let $\diamo'= \diamot(x'_-, x'_+)$ 
$$
\diamo' = \diamot(x'_-, x'_+)=V(x'_-, \st(\tau_+))\cap V(x'_+, \st(\tau_-)) \subset P
$$
be a second diamond in the same parallel set $P$. 
Then the Weyl cones $V(x'_{\mp},\st(\tau_{\pm}))$ are intersections 
$$
V(x'_{\mp},\st(\tau_{\pm}))= P\cap\bigcap_{\zeta\in Z_{\pm}} \{b'_{\zeta}\le 0\}
$$
of horoballs with the same centers, but the new Busemann functions $b'_{\zeta}$ 
are normalized at the new tips $x'_{\mp}$,
i.e. $b'_{\zeta}(x'_{\mp})=0$ for $\zeta\in Z_{\pm}$.
The second diamond is then defined as the sublevel set
$$\diamo'=\{b'|_P\le 0\}$$
of the convex function 
$$ 
b'=\sup_{\xi\in Z}  b'_{\zeta},
$$
Since the points $x_\pm, x'_\pm$ are the normalization points of the corresponding Busemann functions, 
it follows that for all $\zeta\in Z_{\pm}$ we have 
\begin{equation*}
||b_{\zeta} - b'_{\zeta}||\le d(x_{\pm},x'_{\pm}),
\end{equation*}
and therefore 
\begin{equation}
\label{eq:bb'}
||b - b'||\le  \max(d(x_-,x'_-),d(x_+,x'_+)).
\end{equation}
Here and below, $||\cdot ||$ denotes the supremum-norm of functions $X\to \R$. 

\begin{lem}
There exist constants $c(\Theta),\de(\Theta)>0$
such that the following holds:

If $\max(d(x_1^-,{x'_1}^-),d(x_1^+,{x'_1}^+))\leq d\leq\de(\Theta) d(x_1^-,x_1^+)$,
then $\diamo\subset N_{c(\Theta)d}(\diamo'\cap\diamo)$.
\end{lem}
\proof
We may assume without loss of generality that the splitting (\ref{eq:splmcha}) is trivial,
i.e. $X=X_1$.

Take a point $y\in\diamo$.
We connect $y$ to the midpoint $m$ of $x_-x_+$ by the geodesic segment $ym$
and consider the behavior of the convex function $b$ along $ym$.
This will provide an estimate for the time when the segment $ym$ enters the other diamond $\diamo'$.
In view of (\ref{eq:bb'}),
we have 
$$ \diamo\cap\diamo'\subset \{b|_P\leq -d\} .$$
Since 
$$b(m)\leq -\half\sin\eps_0(\Theta)d(x_-,x_+)$$
by Corollary~\ref{cor:diamothick}, and 
$$d(y,m)\leq\diam(\diamo)\leq  2 (\cos\rho_0(\Theta))^{-1}d(x_-,x_+)$$
by Lemma~\ref{lem:diambcr}, 
the point $z\in ym$ at distance 
$d(y,z)=\frac{4\sin\eps_0(\Theta)}{\cos\rho_0(\Theta)}d$
satisfies $b(z)\leq-d$.
\qed

\begin{cor}
\label{cor:contdiamopar}
There exist constants $c(\Theta),\de(\Theta)>0$
such that the following holds:
If also the segment $x'_-x'_+\subset P$ is $\Theta$-longitudinal,
and if 
$$\max(d(x_-,x'_-),d(x_+,x'_+))\leq d\leq\de(\Theta) d(x_1^-,x_1^+),$$
then $$\dH(\diamo,\diamo')\leq c(\Theta)d.$$
\end{cor}
\proof
Note that $d(x_1^{\pm},{x'_1}^{\pm})\leq d(x_{\pm},x'_{\pm})$.
By the triangle inequality, 
$$
(1-2\de(\Theta))\cdot d(x_1^-,x_1^+) \leq d({x'_1}^-,{x'_1}^+)
\leq (1+2\de(\Theta))\cdot d(x_1^-,x_1^+) .
$$
After replacing $\de$ by $\de(1-2\de)$ and 
switching the roles of $\diamo$ and $\diamo'$,
the previous lemma yields that also 
$\diamo'\subset N_{c(\Theta)d}(\diamo'\cap\diamo)$.
The assertion follows.
\qed

\medskip
Now we extend our results and 
estimate the Hausdorff distance between arbitrary $\taumod$-diamonds
$\diamo=\diamot(x_-,x_+)$ and $\diamo'=\diamot(x'_-,x'_+)$ 
which do not have to lie in the same parallel set.

We first consider the euclidean building case.
There, nearby diamonds have large overlap:
\begin{lem}[Nearby diamonds in euclidean buildings]
Let $X$ be a euclidean building. 
There exist constants $c(\Theta,\Theta'),\de(\Theta,\Theta')>0$
such that the following holds:
If the segment $x_-x_+$ is $\Theta$-regular and if 
$$\max(d(x_-,x'_-),d(x_+,x'_+))\leq d\leq\de(\Theta,\Theta') d(x_1^-,x_1^+),$$
then 
$$\dH(\diamo,\diamo')\leq c(\Theta,\Theta')d.$$
\end{lem}
\proof
According to Lemma~\ref{lem:lcoaregureg},
if $\de(\Theta,\Theta')$ is chosen sufficiently small,
then the segments $x'_-x'_+$ and $x_{\mp}x'_{\pm}$
are $\Theta'$-regular.

Let $P=P(\tau_-,\tau_+)\supset\diamo$ be an ambient parallel set as considered above.
In order to find a point in the intersection $P\cap\diamo'$ close to $x_+$, 
we apply Corollary~\ref{cor:leavecone}
to the ambient Weyl cone $V(x_-,\st(\tau_+))\supset\diamo$
and the $\Theta'$-regular segment $x_-x'_+$ (respectively, a ray extending it).
We obtain a point 
$y_+\in x_-x'_+\cap P$
at distance $\leq (\sin\eps_0(\Theta'))^{-1}d$ from $x'_+$ and, 
consequently, distance  $\leq (1+(\sin\eps_0(\Theta'))^{-1})d$ from $x_+$.

Since $\diamot(x_-,y_+)\subset\diamot(x_-,x'_+)\cap P$,
Corollary~\ref{cor:contdiamopar} yields estimates
for the Hausdorff distances of $\diamot(x_-,y_+)$ 
from the diamonds $\diamot(x_-,x'_+)$ and $\diamo$,
and hence for the Hausdorff distance between the latter two diamonds.
The estimates are linear in $d$ with constants only depending on $\Theta'$.
Note hereby that all diamonds split off the same $X_2$-factor,
and that the quantity $d(x_1^-,x_1^+)$,
which appears as a bound in the hypothesis of Corollary~\ref{cor:contdiamopar}, 
varies continuously with the pair $(x_-,x_+)$.

Similarly, working with an ambient Weyl cone $V(x'_+,\st(\tau'_-))\supset\diamo'$
and the $\Theta'$-regular segment $x_-x'_+$, 
one obtains a point 
$y_-\in x_-x'_+\cap V(x'_+,\st(\tau'_-))$ uniformly close to $x'_-$.
Using the intermediate diamond
$\diamot(y_-,x'_+)\subset\diamot(x_-,x'_+)\cap V(x'_+,\st(\tau'_-))$,
one estimates the Hausdorff distance from 
$\diamot(x_-,x'_+)$ to $\diamo'$.
\qed

\medskip
We return to the general model space case and are now ready to show:
\begin{prop}[Continuity of diamonds]
\label{prop:contdiamo}
The $\taumod$-diamonds in $X$ 
depend continuously,
with respect to the Hausdorff topology,
on their $\taumod$-regular pair of tips.
\end{prop}
\proof
If $X$ is a euclidean building,
this is a direct consequence of the previous lemma.
We assume therefore that $X$ is a symmetric space. 

Consider a diamond $\diamo=\diamot(x_-,x_+)$ and an ambient parallel set $P$.
As a consequence of the Hausdorff distance estimates for diamonds in the same parallel set,
cf. Corollary~\ref{cor:contdiamopar},
there exists $\de>0$ such that $\diamo$ has Hausdorff distance $<\half\eps$ 
from all diamonds $\diamot(x'_-,x'_+)$ with $x'_{\pm}\in B(x_{\pm},\de)\cap P$.

Let $U\subset\Isom(X)$ be a neighborhood of the identity
such that $d(ux,x)<\half\eps$ for all $x\in N_{\eps}(\diamo)$ and all $u\in U$.
Then the diamonds $\diamot(ux'_-,ux'_+)$ 
with $u\in U$ and $x'_{\pm}\in B(x_{\pm},\de)\cap P$
are $\eps$-Hausdorff close to $\diamo$. 
The pairs of tips of these diamonds form a neighborhood of $(x_-,x_+)$
in $X\times X$, because the manifold 
${\cal P}_{\taumod}\subset\D_{\taumod^-}X\times\D_{\taumod^+}X$
of type $\taumod$ parallel sets $P(\tau_-,\tau_+)$, respectively,
of pairs $(\tau_-,\tau_+)$ of opposite simplices $\tau_{\pm}$ of types
$\taumod^{\pm}$,
is a homogeneous space for the Lie group $G=\Isom_o(X)$. 
\qed

\subsection{Topology at infinity and partial bordification}
\label{sec:flagtop}

We will describe the topologies 
on the visual compactification $\bar X=X\cup\geo X$ 
and on the $\taumod$-bordification $\barXt=X\cup\DtX$
in terms of shadows and related ``basic subsets''. 

We need the following notions of shadows at infinity 
in $\geo X$ and $\DtX$.
\begin{dfn}[Shadows at infinity]
\label{dfn:shadgen}
(i) For points $x,y\in X$ we define the {\em shadow}
of the point $y$ as seen from $x$ by 
$$ \pSh_{x,y}:=\{\xi:y\in x\xi\} \subset\geo X,$$
and for $r>0$ the shadow of the open $r$-ball around $y$ by 
$$ \bSh_{x,y,r}:=\{\xi:x\xi\cap B(y,r)\neq\emptyset\} \subset\geo X.$$
(ii) For points $x,y\in X$ we define the {\em $\taumod$-shadow}
of the point $y$ as seen from $x$ by 
$$ \pSht_{x,y}:=\{\tau:y\in V(x,\st(\tau))\} \subset\DtX,$$
and for $r>0$ the $\taumod$-shadow of the open $r$-ball around $y$ by 
$$ \bSht_{x,y,r}:=\{\tau:V(x,\st(\tau))\cap B(y,r)\neq\emptyset\} \subset\DtX,$$
\end{dfn}
By coning off the shadows at infinity at points in $X$ and removing large balls around their tips, 
we obtain the subsets of $\bar X$ and $\barXt$ 
which we will use to describe, respectively, construct the natural topologies.
\begin{dfn}[Basic subsets]
(i) For points $x,y\in X$ and radii $r>0$, we define the subsets
$$\pO_{x,y}:=\{z:z\neq y\in xz\}\subset X,$$
and
$$\bO_{x,y,r}:=\{z:xz\cap B(y,r)\neq\emptyset\}\subset X,$$
and the {\em basic subsets}
$$\barpO_{x,y}:=\pO_{x,y}\cup\pSh_{x,y}\subset\bar X\quad\hbox{ and }\quad
\barbO_{x,y,r}:=\bO_{x,y,r}\cup\bSh_{x,y,r}\subset\bar X.$$
(ii) For points $x,y\in X$ and radii $r>0$, we define the subsets
$$\pOt_{x,y}:=\{z:\hbox{$xz$ $\taumod$-regular and $z\neq y\in\diamot(x,z)$}\}\subset X$$
and
$$\bOt_{x,y,r}:=\{z:\hbox{$xz$ $\taumod$-regular and $\diamot(x,z)\cap B(y,r)\neq\emptyset$}\}\subset X,$$
and the {\em $\taumod$-basic subsets}
$$\barpOt_{x,y}:=\pOt_{x,y}\cup\pSht_{x,y}\subset\barXt\quad\hbox{ and }\quad
\barbOt_{x,y,r}:=\bO_{x,y,r}\cup\bSh_{x,y,r}\subset\barXt.$$
\end{dfn}

We observe the following relations between point and ball shadows:
$$ \bSh_{x,y,r}=\bigcup_{z\in B(y,r)} \pSh_{x,z} 
\quad\hbox{ and }\quad
\bSht_{x,y,r}=\bigcup_{z\in B(y,r)} \pSht_{x,z} $$
There are analogous relations between point and ball type basic subsets:
$$ \bO_{x,y,r}=\bigcup_{z\in B(y,r)} \pO_{x,z} 
\quad\hbox{ and }\quad
\bOt_{x,y,r}=\bigcup_{z\in B(y,r)} \pOt_{x,z} $$
We note that 
the $\taumod$-versions of the shadows and basic subsets 
are generalizations of these to arbitrary rank 
and agree with them in rank one. 

\medskip
We first recall the description of the visual topology 
on the visual compactification $\bar X$.
\begin{fact}
\label{fact:topviscomp}
(i) For every point $x\in X$,
the basic subsets $\bO_{x,\cdot,\cdot}$ form
together with the open subsets of $X$ a basis 
of the visual topology on $\bar X$.

(ii) For every ray $x\xi\subset X$,
every sequence $y_n\to\infty$ of points $y_n\in x\xi$
and every bounded sequence of radii $r_n>0$
the basic subsets $\bO_{x,y_n,r_n}$ form a neighborhood basis of $\xi$.
\end{fact}
This restricts to the following description of the visual topology on $\geo X$.
\begin{fact}
\label{fact:vistop}
(i) For every point $x\in X$,
the shadows $\bSh_{x,\cdot,\cdot}$ form a basis of the visual topology on $\geo X$.
If $X$ is a euclidean building,
then also the shadows $\pSh_{x,\cdot}$ form a basis.

(ii) For every ray $x\xi\subset X$,
every sequence $y_n\to\infty$ of points $y_n\in x\xi$ and every bounded sequence of radii $r_n>0$,
the shadows $\bSh_{x,y_n,r_n}$ form a neighborhood basis of $\xi$.
If $X$ is a euclidean building,
then also the shadows $\pSh_{x,y_n}$ form a neighborhood basis.
Moreover, 
if $X$ is a symmetric space, then for $x\neq y\in x\xi$
also the shadows $\bSh_{x,y,\cdot}$ form a neighborhood basis.
\end{fact}

Now we construct 
natural topologies on $\DtX$ and, at least partially, on $\barXt$.
\begin{lem}
The subsets $\bOt_{\cdot,\cdot,\cdot}$ are open in $X$. 
If $X$ is a euclidean building,
then also the subsets $\pOt_{\cdot,\cdot}$ are open.
\end{lem}
\proof
The openness of $\bOt_{\cdot,\cdot,\cdot}$
follows from the semicontinuity of diamonds, 
cf. Lemma~\ref{lem:semicontdiamo}.
If $X$ is a euclidean building,
then the openness of $\pOt_{\cdot,\cdot}$
is a consequence of Corollary~\ref{cor:leavecone}.
\qed

\begin{lem}
If $xy$ is $\taumod$-regular, 
then $y\in B(y,r)\subset\bOt_{x,y,r}$ 
for all sufficiently small $r>0$.
\end{lem}
\proof
If $r$ is sufficiently small,
then the segments $xz$ are $\taumod$-regular
for all $z\in B(y,r)$.
\qed

\begin{lem}
(i) If $z\in\bOt_{x,y,r}$, then there exists $s>0$ such that $\barbOt_{x,z,s}\subset\barbOt_{x,y,r}$.

(ii) If $X$ is a euclidean building 
and $z\in\pOt_{x,y}$, 
then there is $s>0$ with $\barbOt_{x,z,s}\subset\barpOt_{x,y}$.
\end{lem}
\proof 
(i) Due to the semicontinuity of diamonds, see Lemma~\ref{lem:semicontdiamo},
there exists $s>0$ such that for every $z'\in B(z,s)$
the segment $xz'$ is $\taumod$-regular and 
the diamond $\diamot(x,z')$ still intersects the ball $B(y,r)$.

(ii) The argument is the same,
but using Corollary~\ref{cor:leavecone}.
It implies that 
there exists $s>0$ such that for every $z'\in B(z,s)$
the segment $xz'$ is $\taumod$-regular and 
the diamond $\diamot(x,z')$ still contains $y$.
\qed

\begin{cor}
\label{cor:basfltop}
(i) The subsets $\barbOt_{x,\cdot,\cdot}$ form together with the open subsets of $X$
the basis of a topology ${\cal T}_x$ on $\barXt$.
If $X$ is a euclidean building,
then also the subsets $\barpOt_{x,\cdot}$ form a basis. 

(ii) For every simplex $\tau\in\DtX$, 
every asymptotically uniformly $\taumod$-regular sequence $y_n\to\infty$ in $V(x,\st(\tau))$
and every bounded sequence of radii $r_n>0$,
the basic subsets $\barbOt_{x,y_n,r_n}$
form a neighborhood basis for $\tau$ in $(\barXt,{\cal T}_x)$.
If $X$ is a euclidean building,
then also the subsets $\barpOt_{x,y_n}$ form a neighborhood basis. 
In particular, ${\cal T}_x$ is first-countable.
\end{cor}
\proof 
(i) Suppose that $\tau$ belongs to a finite intersection $\cap_i\bSht_{x,y_i,r_i}$. 
This means that $V(x,\st(\tau))$ intersects all balls $B(y_i,r_i)$.
Let $z\in V(x,\ost(\tau))-\{x\}$ be a point 
so that $\diamot(x,z)$ also intersects them.
Then $z\in\cap_i\bOt_{x,y_i,r_i}$.
With the lemma it follows that 
$\tau\in\barbOt_{x,z,s}\subset\cap_i\barbOt_{x,y_i,r_i}$ for all sufficiently small $s$.
Furthermore, $\cap_i\bOt_{x,y_i,r_i}$ is open in $X$.

The subsets $\barbOt_{x,\cdot,\cdot}$ are unions of subsets of the form $\barpOt_{x,\cdot}$.
If $X$ is a euclidean building,
then vice versa the subsets $\barpOt_{x,\cdot}$ 
are unions of subsets of the form $\barbOt_{x,\cdot,\cdot}$
by the last lemma.

(ii) Suppose that $\tau\in\bSh_{x,y,r}$
and that $\barbOt_{x,y_n,r_n}\not\subset\barbOt_{x,y,r}$ for all $n$.
Then there exist points $z_n\in B(y_n,r_n)$ such that 
$xz_n$ is $\taumod$-regular and 
$\diamot(x,z_n)\cap B(y,r)=\emptyset$.

If $X$ is locally compact, 
then after passing to a subsequence, 
$xz_n$ subconverges to a ray $x\zeta\subset V(x,\st(\tau))$ with $\zeta\in\ost(\tau)$.
Let $w\in x\zeta$ be a point such that $y\in\diamot(x,w)$,
and let $w_n\in xz_n$ be points converging to it, $w_n\to w$.
Then $\diamot(x,w_n)\cap B(y,r)\neq\emptyset$ for large $n$, 
due to the semicontinuity of diamonds, 
see Lemma~\ref{lem:semicontdiamo},
and hence also $\diamot(x,z_n)\cap B(y,r)\neq\emptyset$, a contradiction. 

If $X$ is a euclidean building, 
then it follows with Corollary~\ref{cor:leavecone}
that $y\in\diamot(x,z_n)$, for large $n$,
which is also a contradiction. 

Thus, the subsets $\barbOt_{x,y_n,r_n}$ form a neighborhood basis.
If $X$ is a euclidean building, 
it follows that also the smaller open subsets $\barpOt_{x,y_n}\subset\barbOt_{x,y_n,r_n}$ 
form a neighborhood basis.
\qed

\medskip
We compare now the topologies ${\cal T}_x$ for different base points $x$.

By construction, they all 
restrict to the given topology on $X$.

Regarding the comparison of the topologies ${\cal T}_x$
at infinity on $\DtX$ and on the entire bordification $\barXt$,
we use that if a topological space is first-countable,
then its topology is determined by the sequential convergence.
Namely, a subset is a neighborhood of a point,
if and only if it cannot be avoided by a sequence converging to this point. 
We therefore compare sequential convergence for the topologies ${\cal T}_x$.
We will do this only partially, 
namely for arbitrary sequences in $\DtX$,
but only for asymptotically uniformly $\taumod$-regular sequences in $X$.
This will be sufficient for the purposes of this paper. 

We first observe that ${\cal T}_x$-convergence translates into Hausdorff convergence of diamonds and Weyl cones.
\begin{lem}
\label{lem:txconv}
(i) The convergence $\tau_n\to\tau$ in $\DtX$ with respect to ${\cal T}_x$
is equivalent to the Hausdorff convergence 
$V(x,\st(\tau_n)))\cap\bar B(x,R)\to V(x,\st(\tau))\cap\bar B(x,R)$ of truncated Weyl cones
for all radii $R>0$.

(ii)
For an asymptotically uniformly $\taumod$-regular sequence $x_n\to\infty$ in $X$,
the convergence $x_n\to\tau$ in $(\barXt,{\cal T}_x)$
is equivalent to the Hausdorff convergence $\diamot(x,x_n)\cap\bar B(x,R)\to V(x,\st(\tau))\cap\bar B(x,R)$ 
of truncated diamonds for all radii $R>0$
\end{lem}
\proof
The first statement follows from the second one in view of Lemma~\ref{lem:diamoco}.

For the second statement, suppose that $x_n\to\tau$.
Then for every point $y\in V(x,\st(\tau))$ and radius $r>0$,
the diamonds $\diamot(x,x_n)$ intersect $B(y,r)$ for all sufficiently large $n$.
Hence,
$d(y,\diamot(x,x_n))\to 0$ as $n\to+\infty$,
and the continuity of diamonds 
(Proposition~\ref{prop:contdiamo})
implies that 
$\diamot(x,y)\subset N_{\eps_n}(\diamot(x,x_n))$ with a sequence $\eps_n\to0$.
Again in view of Lemma~\ref{lem:diamoco},
this yields the asserted Hausdorff convergence.
The converse direction is clear.
\qed

\begin{cor}
The topology ${\cal T}_x$ on $\barXt$ is Hausdorff.
\end{cor}
\proof
This is a direct consequence of first-countability and the last lemma,
because it implies that limits of sequences are unique.
\qed

\medskip
We now compare the topologies ${\cal T}_x$ on $\DtX$.
We do this by comparing them to the visual topology on $\geo X$.
For every type $\bar\xi\in\inte(\taumod)$,
there is the natural identification
\begin{equation}
\label{eq:canbijvis}
\theta^{-1}(\bar\xi)\buildrel1:1\over\lra\DtX
\end{equation}
with the subspace $\theta^{-1}(\bar\xi)\subset\geo X$,
assigning to a point $\xi\in\geo X$ with type $\theta(\xi)=\bar\xi$
the type $\taumod$ simplex $\tau$ spanned by it, $\xi\in\inte(\tau)$.
\begin{lem}
For every type $\bar\xi\in\inte(\taumod)$ and every point $x\in X$,
the bijection (\ref{eq:canbijvis}) 
is a homeomorphism with respect to the restrictions of 
the visual topology on $\bar X$ to $\theta^{-1}(\bar\xi)$ 
and the topology ${\cal T}_x$ on $\barXt$ to $\DtX$.
\end{lem}
\proof
Let $(\xi_n)$ and $\xi$ be a sequence and a point in $\theta^{-1}(\bar\xi)\subset\geo X$,
and let $(\tau_n)$ and $\tau$ be the corresponding sequence and point in $\DtX$,
i.e. $\xi_n\in\inte(\tau_n)$ and $\xi\in\inte(\tau)$.
We must show that $\xi_n\to\xi$ if and only if $\tau_n\to\tau$ 
with respect to the topologies in consideration. 

Suppose that $\tau_n\to\tau$ with respect to ${\cal T}_x$.
Then $V(x,\st(\tau_n))\to V(x,\st(\tau))$ by the previous lemma.
In particular, increasingly long subsegments $xy_n\subset x\xi_n\subset V(x,\st(\tau_n))$
become arbitrarily close to segments $x\bar y_n\subset V(x,\st(\tau))$.
We want to find Hausdorff close segments in $V(x,\st(\tau))$ of the same type $\bar\xi$.
By the triangle inequality for $\De$-lengths (\ref{eq:tridele}),
$\|d_\De(x, y_n) - d_\De(x,\bar y_n)\|\le  d(y_n,\bar y_n)\to0 $
and, in a euclidean Weyl chamber through $\bar y_n$ with tip $x$,
we find a point $z_n\in V(x,\st(\tau))$ 
with $d_\De(x, z_n)=d_\De(x, y_n)$.
Then $d(z_n,\bar y_n)=\|d_\De(z_n,\bar y_n)\|=\|d_\De(x, z_n)-d_\De(x,\bar y_n)\|\to0$,
and hence $d(z_n,y_n)\to0$ by the triangle inequality.
Moreover, $\theta(xz_n)=\theta(xy_n)=\bar\xi$ and therefore $z_n\in x\xi$,
because $\xi$ is the only point in $\st(\tau)$ with type $\bar\xi$.
It follows that $x\xi_n\to x\xi$,
i.e. $\xi_n\to\xi$.

Conversely, 
suppose that $\xi_n\to\xi$, i.e. $x\xi_n\to x\xi$.
Then any ball centered at $x\xi$ is intersected by $x\xi_n$ for all sufficiently large $n$.
Thus, $\tau_n\to\tau$ by our description of ${\cal T}_x$-neighborhood bases.
\qed

\begin{cor}
The restriction of the topology ${\cal T}_x$ to $\DtX$ does not depend on $x$.
\end{cor}
\begin{dfn}[Visual topology]
We call this topology on $\DtX$ the {\em visual topology}.
\end{dfn}

Now we show that the topologies ${\cal T}_x$ agree on the entire bordification $\barXt$ 
``in $\taumod$-regular directions''.
We reformulate the condition for ${\cal T}_x$-convergence 
for asymptotically uniformly $\taumod$-regular sequences $x_n\to\infty$ 
given in Lemma~\ref{lem:txconv} above,
in order to show its independence of $x$.
We do this separately in the symmetric space (locally compact) and euclidean building cases.

In the locally compact case,
we can express ${\cal T}_x$-convergence
in terms of accumulation at infinity (the limit set) in $\bar X$:
\begin{lem}
\label{lem:txconvlc}
Suppose that $X$ is locally compact.
Then $x_n\to\tau\in\DtX$ with respect to ${\cal T}_x$, 
if and only if 
the accumulation set of $(x_n)$ in $\bar X$ 
(with respect to the visual topology of $\bar X$)
is contained in $\ost(\tau)\subset\geo X$.
\end{lem}
\proof
Since $X$ is locally compact,
the sequence $(x_n)$ subconverges in both $\bar X$ and $\barXt$.
The latter holds, because the sequence of diamonds $\diamot(x,x_n)$ Hausdorff subconverges
and, in view of Lemma~\ref{lem:diamoco},
the Hausdorff sublimits must be type $\taumod$ Weyl cones. 
Note also that $(x_n)$ accumulates in $\bar X$
only at the $\taumod$-regular part $\theta^{-1}(\ost(\taumod))$ of $\geo X$,
as a consequence of asymptotic uniform $\taumod$-regularity.

Therefore,
if the assertion is wrong,
we may assume after passing to a subsequence, 
that $x_n\to\tau$ in $\barXt$ and $x_n\to\xi'\in\ost(\tau')$ in $\bar X$ 
for different simplices $\tau,\tau'\in\DtX$.
But then $\diamot(x,x_n)\to V(x,\st(\tau))$
according to Lemma~\ref{lem:txconv}.
Since $xx_n\to x\xi'$, it follows that $\xi'\in\st(\tau)$,
a contradiction. 
\qed

\medskip
In the euclidean building case, 
we can strengthen the condition of Hausdorff convergence
of Weyl cones to initial coincidence up to increasing radii.
\begin{lem}
\label{lem:txconveb}
Suppose that $X$ is a euclidean building.
Then $x_n\to\tau\in\DtX$ with respect to ${\cal T}_x$, 
if and only if 
for every $R>0$ it holds that 
$\diamot(x,x_n)\cap\bar B(x,R)=V(x,\st(\tau))\cap\bar B(x,R)$ for all sufficiently large $n$.
\end{lem}
\proof
This is a consequence of Lemmas~\ref{lem:txconv} and~\ref{lem:closediamos}.
\qed

\begin{cor}
Whether an asymptotically uniformly $\taumod$-regular sequence $x_n\to\infty$ in $X$
converges to a simplex $\tau\in\DtX$ in $(\barXt,{\cal T}_x)$,
does not depend on $x$.
\end{cor}
\proof
If $X$ is locally compact, this follows immediately from Lemma~\ref{lem:txconvlc}.
Assume therefore that $X$ is a euclidean building. 

Let $x,x'\in X$ and suppose that $x_n\to\tau\in\DtX$ with respect to ${\cal T}_x$.
By Lemma~\ref{lem:txconveb},
there exists a sequence $y_n\to\infty$ of points $y_n\in xx_n\cap V(x,\st(\tau))$.
Let $y'_n\in x'x_n$ be points uniformly close to the points $y_n$,
e.g. such that $d(y'_n,y_n)\leq d(x',x)$.
Then the sequence $(y'_n)$ is contained in a tubular neighborhood (of radius $d(x',x)$)
of $V(x,\st(\tau))$,
and hence also 
in a tubular neighborhood (of radius $2d(x',x)$)
of $V(x',\st(\tau))$,
because the two Weyl cones have finite Hausdorff distance ($\leq d(x',x)$).
The sequences $(y_n)$ and $(y'_n)$ inherit from $(x_n)$ asymptotically uniform $\taumod$-regularity.

Consider the subsegments $x'z'_n=x'y'_n\cap V(x',\st(\tau))$.
According to Corollary~\ref{cor:leavecone},
the distances $d(z'_n,y'_n)$ are uniformly bounded,
and therefore $z'_n\to\infty$.
Since $\diamot(x',z'_n)\subset\diamot(x',x_n)\cap V(x',\st(\tau))$,
it follows,
using again Lemma~\ref{lem:diamoco}, 
that $\diamot(x',x_n)\to V(x',\st(\tau))$.
Hence, $x_n\to\tau$ also with respect to ${\cal T}_{x'}$.
\qed

\medskip
The corollary justifies the following definition.
\begin{dfn}[Flag convergence]
We say that an asymptotically uniformly $\taumod$-regular sequence $x_n\to\infty$ in $X$
{\em flag converges} to a simplex $\tau\in\DtX$, if $x_n\to\tau$ in $(\barXt,{\cal T}_x)$ for some base point $x$.
\end{dfn}

Now we can also make precise the coincidence of the topologies ${\cal T}_x$ 
``in $\taumod$-regular directions''.
Suppose that $A\subset X$ is an asymptotically uniformly $\taumod$-regular subset,
and consider the subset 
$$\tiAt:=A\cup\DtX\subset\barXt.$$
\begin{cor}
The topology induced by ${\cal T}_x$ on $\tiAt$ does not depend on $x$.
\end{cor}
\begin{dfn}[Topology of flag convergence]
We call this topology on $\tiAt$ the {\em topology of flag convergence}.
\end{dfn}
As shown above,
the topologies ${\cal T}_x$
and hence the topology of flag convergence on $\tiAt$
are Hausdorff and first-countable.
Neighborhood bases at infinity have been described in Corollary~\ref{cor:basfltop}.

\medskip
We further discuss the flag convergence of sequences. 

A situation when an asymptotically uniformly regular sequence flag converges,
is when it stays close to a Weyl cone:
\begin{lem}
\label{lem:flagconv}
Suppose that the asymptotically uniformly $\taumod$-regular sequence $x_n\to\infty$
is contained in the tubular neighborhood of the type $\taumod$ Weyl cone $V(x,\st(\tau))$.
Then $x_n\to\tau$.
\end{lem}
\proof
If $X$ is locally compact, this follows from Lemma~\ref{lem:txconvlc}.

Suppose therefore that $X$ is a euclidean building. 
Consider the points $y_n$ where the segments $xx_n$ exit the Weyl cone $V(x,\st(\tau))$,
i.e. $xy_n=xx_n\cap V(x,\st(\tau))$.
Then Corollary~\ref{cor:leavecone} implies that $d(y_n,x_n)$ is bounded.
Hence $y_n\to\infty$ is an asymptotically uniformly $\taumod$-regular sequence in $V(x,\st(\tau))$,
and $x_n\in\barpOt_{x,y_n}$.
The basic subsets $\barpOt_{x,y_n}$ form a neighborhood basis of $\tau$.
Thus, $x_n\to\tau$ also in this case.
\qed

\medskip
We give a name to this stronger form of flag convergence:
\begin{dfn}[Conical convergence, cf. {\cite[Def. 6.1]{morse}}]
\label{dfn:conicaltau}
We say that an asymptotically uniformly $\taumod$-regular sequence $x_n\to\infty$ in $X$ 
{\em flag converges conically} to $\tau\in\DtX$
if it is contained in a tubular neighborhood of the Weyl cone $V(x,\st(\tau))$. 
\end{dfn}

\begin{cor}
\label{cor:closewco}
Let $V(x,\st(\tau))$ and $V(x',\st(\tau'))$ be type $\taumod$ Weyl cones.
Suppose that for some $D>0$ the intersection of their $D$-neighborhoods 
contains an asymptotically uniformly $\taumod$-regular sequence.
Then $\tau=\tau'$.
\end{cor}
\proof
If $(x_n)$ is such an asymptotically uniformly $\taumod$-regular sequence, 
then $x_n\to\tau$ and $x_n\to\tau'$.
The assertion follows from the Hausdorff property of the topologies ${\cal T}_x$.
\qed

\medskip
The following {\em convergence criterion} will be useful when $X$ is not locally compact.
\begin{lem}
\label{lem:taunconv}
Let $x_n\to\infty$ be an asymptotically uniformly $\taumod$-regular sequence in $X$,
and let $(\tau_n)$ be a sequence in $\DtX$ 
such that for some point $x\in X$ and some constant $D\geq0$ 
it holds that $x_n\in\bar N_D(V(x,\st(\tau_m)))$ for all $m\geq n$.
Then the sequence $(\tau_n)$ converges, $\tau_n\to\tau\in\DtX$, 
and $x_n\in\bar N_D(V(x,\st(\tau)))$ for all $n$. 
In particular, $x_n\to\tau$ conically.
\end{lem}
\proof
If $X$ is locally compact,
then there exists a convergent subsequence of simplices, $\tau_{n_k}\to\tau$.
It follows that $x_n\in\bar N_D(V(x,\st(\tau)))$ for all $n$, and the assertion holds in this case.

Suppose therefore that $X$ is a euclidean building (because otherwise $X$ is locally compact).
For suitable $\Theta$,
the segments $xx_n$ are $\Theta$-regular for large $n$.
Let $\tau'_n\in\pSht_{x,x_n}$.
Applying Lemma~\ref{lem:closediamos},
we obtain for any radius $r>0$ that 
$$V(x,\st(\tau'_n))\cap\bar B(x,r)=V(x,\st(\tau_m))\cap\bar B(x,r)$$
for $m\geq n\geq n_0(r)$.
Thus, both sides are independent of $m$ and $n$,
and isometric to 
$$V(x,\st(\tau_m))\cap\bar B(x,r)=C(x,r)$$
for $m\geq n_0(r)$. 
The union of the nested family of cones $C(x,r)$ as $r\to+\infty$ is a type $\taumod$ Weyl cone 
$V(x,\st(\tau))$.
It follows that $\tau_m\to\tau$ and $x_n\in\bar N_D(V(x,\st(\tau)))$.
\qed

\subsection{Ultralimits of parallel sets, Weyl cones and diamonds}
\label{sec:ulimps}

Let $X$ be a model space.

For a sequence of base points $\star_n\in X$ 
and a sequence of scale factors $\la_n>0$ with $\ulim\la_n=0$,
we consider the ultralimit 
$$(X_{\om},\star_{\om})=\ulim(\la_nX,\star_n)$$
of rescaled copies of $X$.
We will use the following result :
\begin{thm}
[B. Leeb, B. Kleiner, {\cite[ch.\ 5]{qirigid}}] 
$X_{\om}$ is a euclidean building of the same rank and type 
as the model space $X$.
\end{thm} 

We will need later that certain families of subsets are closed under taking ultralimits.

Sequences of maximal flats in $X$
ultraconverge to maximal flats in $X_{\om}$,
see also \cite[ch.\ 5]{qirigid}:
If $F_n\subset X$, $n\in\N$, are maximal flats 
such that $\ulim \la_n d(F_n,\star_n)<+\infty$,
then $$F_{\om}:=\ulim \la_nF_n\subset X_{\om}$$
is a maximal flat.
Furthermore,
if $\kappa_n:F_n\to F_{mod}$ are charts such that 
$\ulim \la_n d(\kappa_n^{-1}(0),\star_n)<+\infty$
for the base point $0\in\Fmod$,
then the ultralimit 
\begin{equation}
\label{eq:ulimchfl}
\kappa_{\om}^{-1}:=\ulim\kappa_n^{-1}:\Fmod\to X_{\om}
\end{equation}
of the isometric embeddings $\kappa_n^{-1}:\la_n\Fmod\to \la_nX_n$
is an isometric embedding,
and it is the inverse of a chart $\kappa_{\om}$ for $F_{\om}$.
(Note that 
$\ulim_n(\la_n\Fmod,0)\cong(\Fmod,0)$ 
canonically,
because $\Fmod$ is self similar, $(\la_n\Fmod,0)\cong(\Fmod,0)$ canonically, and locally compact.)

Euclidean Weyl sectors (chambers) ultraconverge to euclidean Weyl sectors (chambers),
if their tips ultraconverge:
Let $V(x_n,\tau_n)\subset X_n$ be Weyl sectors 
and suppose that $x_{\om}=(x_n)\in X_{\om}$ exists. 
Since sectors are contained in maximal flats,
we may assume that $V(x_n,\tau_n)\subset F_n$
and work with the charts $\kappa_n$ and $\kappa_{\om}$.
Then $(\geo\kappa_n)\tau_n\subset\geo\Fmod$ is one of finitely many faces,
and therefore $\om$-always the same face $\bar\tau$.
We put $\tau_{\om}:=(\geo\kappa_{\om}^{-1})\bar\tau \subset \geo F_{\om}$ and obtain that 
\begin{equation}
\label{eq:ulimsect}
\ulim \la_nV(x_n,\tau_n)=V(x_{\om},\tau_{\om})
\end{equation}
and, regarding types,
$\theta(\tau_{\om})=\theta(\tau_n)$ for $\om$-all $n$.
Applying (\ref{eq:ulimsect}) 
to the chambers $\bar\si\supset\bar\tau$ in $\geo\Fmod$ and taking the union,
one obtains in particular that 
\begin{equation}
\label{eq:ulimstinfl}
\ulim \la_nV(x_n,\st(\tau_n)\cap\geo F_n)=V(x_{\om},\st(\tau_{\om})\cap\geo F_{\om}),
\end{equation}
a fact, which will be useful below.

Generalizing the fact for maximal flats,
we will show next that parallel sets ultraconverge to parallel sets.
Consider a sequence of parallel sets $P_n=P(\tau_n^-,\tau_n^+)\subset X$ 
and assume that 
$\ulim \la_n d(\star_n,P_n)<+\infty$.
Let 
$$P_{\om}:=\ulim \la_nP_n \subset X_{\om} .$$
\begin{lem}[Ultralimits of parallel sets]
\label{lem:ulimps}
$P_{\om}$ is again a parallel set, i.e.\
$P_{\om}=P(\tau_{\om}^-,\tau_{\om}^+)$ 
with a pair of opposite 
simplices $\tau_{\om}^{\pm}\subset\geo X_{\om}$.
Moreover, 
$\theta(\tau_{\om}^{\pm})=\theta(\tau_n^{\pm})$ for $\om$-all $n$.
\end{lem}
\proof
We may assume without loss of generality that $\star_n\in P_n$
and that $\theta(\tau_n^{\pm})=\taumod^{\pm}$ for all $n$.

In order to represent the $P_n$ as parallel sets of geodesic lines,
we fix a type $\bar\xi\in\inte(\taumod)$
and denote by $\xi_n\in\inte(\tau_n)$
the ideal points of type $\theta(\xi_n)=\bar\xi$.
Then $P_n=P(l_n)$ with the oriented geodesic line $l_n\subset P_n$ extending the ray $\star_n\xi_n$.

The ultralimit of lines 
$$l_{\om}:=\ulim \la_nl_n  \subset P_{\om}$$ 
is again an oriented line of type $\bar\xi$.
Let $\xi_{\om}:=l_{\om}(+\infty)\in\geo X_{\om}$
denote its forward ideal endpoint,
and let $\tau_{\om}^{\pm}\subset\geo X_{\om}$ 
denote the type $\taumod^{\pm}$ simplices spanned by the ideal endpoints $l_{\om}(\pm\infty)$.

Since every point in $P_n$ is contained in a maximal flat $F_n\supset l_n$,
and since sequences of maximal flats $F_n\supset l_n$
ultraconverge to maximal flats $F_{\om}\supset l_{\om}$,
we have that 
$$ P_{\om}\subset P(l_{\om}).$$
We must show that $P_{\om}$ fills out $P(l_{\om})$.
Note that, as an ultralimit of subsets, $P_{\om}$ is closed.

Let $x_{\om}=(x_n)\in X_{\om}$.
The ray $x_{\om}\xi_{\om}\subset X_{\om}$ is the ultralimit 
of the rays $x_n\xi_n\subset X_n$.
We apply Proposition~\ref{prop:diveparsymm} and Proposition~\ref{prop:diveparbuil} 
to conclude that $x_{\om}\xi_{\om}$ dives within uniformly bounded time into $P_{\om}$.
Namely, fix a constant $d>0$ 
and choose $\Theta\ni\bar\xi$.
Let $C'':=\max(C(\Theta,d),(\sin\eps_0(\Theta))^{-1})$ with the constants 
appearing in these results.
Let $y_n\in x_n\xi_n$ be the point at distance 
$d(x_n,y_n)=C''\cdot d(x_n,P_n)$.
Then $y_n\in N_d(P_n)$. 
(If $X$ is a euclidean building, even $y_n\in P_n$.)
The ultralimit point $y_{\om}\in x_{\om}\xi_{\om}$ is defined
and has distance 
$d(x_{\om},y_{\om})=C'' d(x_{\om},P_{\om})<+\infty$ from $x_{\om}$.
Since $\la_nd\to0$, 
we have that 
$y_{\om}\in\ulim\la_n N_d(P_n)=\ulim N_{\la_nd}(\la_nP_n)= N_0(P_{\om})=P_{\om}$.
Thus, the ultralimit ray $x_{\om}\xi_{\om}$ enters $P_{\om}$ within uniformly bounded time.

As a consequence,
every geodesic line $l'_{\om}\subset X_{\om}$ parallel to $l_{\om}$
must already be contained in $P_{\om}$.
This means that $P_{\om}=P(l_{\om})=P(\tau_{\om}^-,\tau_{\om}^+)$.
\qed

\medskip
Using the result on parallel sets,
we will deduce the ultraconvergence of Weyl cones 
from the ultraconvergence of Weyl sectors.
Consider a sequence of Weyl sectors $V(x_n,\tau_n)\subset X$
which ultraconverge as in (\ref{eq:ulimsect}).
Then the corresponding Weyl cones ultraconverge, too:
\begin{lem}[Ultralimits of Weyl cones]
\label{lem:uconvcon}
$\ulim \la_nV(x_n,\st(\tau_n))=V(x_{\om},\st(\tau_{\om}))$.
\end{lem}
\proof
The left-hand side is the union of the ultralimits 
$\ulim \la_nV(x_n,\si_n)$ 
for all sequences of chambers $\si_n\supset\tau_n$ in $\geo X$.
By (\ref{eq:ulimstinfl}),
these ultralimits are euclidean Weyl chambers $V(x_{\om},\si_{\om})$
with chambers $\si_{\om}\supset\tau_{\om}$ in $\geo X_{\om}$,
i.e.\ $\si_{\om}\subset\st(\tau_{\om})$. 
This shows that 
$$\ulim \la_nV(x_n,\st(\tau_n))\subset V(x_{\om},\st(\tau_{\om})).$$
To verify the reverse inclusion, 
we work inside parallel sets containing the Weyl cones.
Let $\hat\tau_n\subset\geo X$ be faces $x_n$-opposite to the faces $\tau_n$.
Then $V(x_n,\st(\tau_n))\subset P_n=P(\hat\tau_n,\tau_n)$. 

Consider a sequence of maximal flats $F_n\subset P_n$ containing the points $x_n$, 
and the ultralimit flat $F_{\om}=\ulim\la_nF_n$.
Then $\tau_n\subset\geo F_n$ and $\tau_{\om}\subset\geo F_{\om}$.
Applying (\ref{eq:ulimstinfl}) yields that 
$$V(x_{\om},\st(\tau_{\om}))\cap F_{\om}=
V(x_{\om},\st(\tau_{\om})\cap\geo F_{\om})
=\ulim \la_nV(x_n,\st(\tau_n)\cap\geo F_n).$$
The union of all flats $F_{\om}$ arising in this way as ultralimits is precisely $P_{\om}:=\ulim\la_nP_n$.
Hence 
$$V(x_{\om},\st(\tau_{\om}))\cap P_{\om}\subset\ulim \la_nV(x_n,\st(\tau_n)).$$
Now we use that parallel sets ultraconverge to parallel sets.
By Lemma~\ref{lem:ulimps},
$P_{\om}=P(\hat\tau_{\om},\tau_{\om})$
with a face $\hat\tau_{\om}$ which is $x_{\om}$-opposite to $\tau_{\om}$. 
It follows that $V(x_{\om},\st(\tau_{\om}))\subset P_{\om}$
and 
$$V(x_{\om},\st(\tau_{\om}))\subset \ulim \la_nV(x_n,\st(\tau_n)),$$
which finishes the proof.
\qed

\medskip
Finally, we describe ultralimits of sequences of diamonds.
Consider a sequence of $\Theta$-regular segments $x_n^-x_n^+\subset X$
and the $\taumod$-diamonds $\diamo_n:=\diamot(x_n^-,x_n^+)$ spanned by them. 
Let 
$$ \diamoom:=\ulim\la_n\diamo_n .$$
\begin{lem}[Ultralimits of diamonds]
\label{lem:uconvdiamo}
If the sequence of segments $x_n^-x_n^+$ ultraconverges to a segment 
$x_{\om}^-x_{\om}^+\subset X_{\om}$, then 
$\diamoom=\diamot(x_{\om}^-,x_{\om}^+)$.
\end{lem}
\proof
We recall that diamonds are forever.
In order to work inside sequences of parallel sets,
let $(\tau_n^-,\tau_n^+)$ be 
pairs of $(x_n^-,x_n^+)$-op\-po\-site type $\taumod^{\pm}$ simplices in $\geo X$.
Putting $P_n:=P(\tau_n^-,\tau_n^+)$ and 
$V_n^{\pm}:=V(x_n^{\mp},\st(\tau_n^{\pm}))\subset P_n$,
we have that 
$$\diamo_n=V_n^-\cap V_n^+\subset P_n .$$
Lemma~\ref{lem:ulimps} implies that 
$$P_{\om}:=\ulim\la_nP_n=P(\tau_{\om}^-,\tau_{\om}^+)$$
with a pair $(\tau_{\om}^-,\tau_{\om}^+)$ 
of $(x_{\om}^-,x_{\om}^+)$-opposite type $\taumod^{\pm}$ simplices in $\geo X_{\om}$.
Moreover, \newline by Lemma~\ref{lem:uconvcon},
$$V_{\om}^{\pm}:=\ulim\la_nV_n^{\pm}=V(x_{\om}^{\mp},\st(\tau_{\om}^{\pm}))\subset P_{\om}.$$
Clearly,
$$\diamoom\subset V_{\om}^-\cap V_{\om}^+\subset P_{\om} ,$$
and we must prove that 
$\diamoom=V_{\om}^-\cap V_{\om}^+$.

Since the segment $x_{\om}^-x_{\om}^+$ is $\Theta$-regular, 
and hence, in particular, is $\taumod$-regular, 
the intersection of interiors $\inte(V_{\om}^-)\cap\inte(V_{\om}^+)$
is dense in $V_{\om}^-\cap V_{\om}^+$.
Since $\diamoom$ is closed (being an ultralimit of subsets),
it therefore suffices to show that 
$\inte(V_{\om}^-)\cap\inte(V_{\om}^+)\subset\diamoom$.

Let $z_{\om}\in\inte(V_{\om}^-)\cap\inte(V_{\om}^+)$.
We may assume that $z_{\om}=(z_n)$ with $z_n\in P_n$.
Since $z_{\om}\in\inte(V_{\om}^{\pm})$,
the segments $x_{\om}^-z_{\om}$ and $z_{\om}x_{\om}^+$ are longitudinal.
It follows that the segments $x_n^-z_n$ and $z_nx_n^+$ are longitudinal
for $\om$-all $n$.
So, $z_n\in V_n^-\cap V_n^+=\diamo_n$
and $z_{\om}\in\diamoom$.
\qed

\section{Modified Carnot--Finsler metric and its contraction\\ property }\label{sec:Carnot-Finsler}

\subsection{A modified Carnot-Finsler type metric on diamonds}

Suppose that $X$ is a model space.

\begin{defn}
A broken path $x_0x_2\ldots x_k$ in a diamond $\diamo$ is called 
{\em non-longitidinal}, if each segment 
$x_i x_{i+1}$ of this path is nonlongitudinal. 
\end{defn}

On the diamond $\diamo=\diamot(x_-,x_+)$ 
we introduce the {\em pseudo-metric} $d_{\diamo}$
which is obtained by infimizing the length of 
broken non-longitudinal paths $x_0x_1\dots x_k$ in $\diamo$ with $x_0=x_-, x_k=x_+$. 
The triangle inequality and symmetry are clearly satisfied by $d_{\diamo}$, but, 
in general, $d_{\diamo}$ is only a pseudo-metric, because points may have infinite distance. 
 
The modified metric $d_{\diamo}$ is larger than the original metric,
$$d_{\diamo}\geq d|_{\diamo} .$$
It obviously agrees with $d$ in the non-longitudinal directions,
i.e.\ for a non-longitudinal segment $xy\subset\diamo$
we have 
$d_{\diamo}(x,y)=d(x,y)$.
However, $d_{\diamo}$ is strictly larger than $d$ 
in the longitudinal directions:

\begin{lem}
\label{lem:horlong}
For a longitudinal segment $xy$ in $\diamo$, 
we have
\begin{equation*}
d_{\diamo}(x,y)\geq C\cdot d(x,y)
\end{equation*}
with a constant $C=C(\theta(\oa{xy}))>1$ depending continuously on the direction type $\theta(\oa{xy})$.
\end{lem}
\proof
We choose an $(x_-,x_+)$-opposite type $(\iota\taumod,\taumod)$ pair of simplices $(\tau_-,\tau_+)$.
Then $\diamo\subset P(\tau_-,\tau_+)$. 
The longitudinal segment $xy$ 
can be extended inside $P(\tau_-,\tau_+)$ to a longitudinal ray $x\xi$,
i.e.\ $\xi\in\ost(\tau_+)\cup\ost(\tau_-)$, say $\xi\in\ost(\tau_+)$. 
Along $xy$, the Busemann function $b_{\xi}$ decays with minimal possible slope $\equiv-1$. 
On the other hand,
along any {\em non}-longitudinal segment in $P(\tau_-,\tau_+)$, 
and hence along any piecewise non-longitudinal geodesic path in $\diamo$ connecting $x$ to $y$,
it has slope $\geq-1+\eps$ with a constant $\eps>0$
depending continuously on the (Tits) distance of $\xi$ from $\D\st(\tau_+)$,
which in turn depends only on $\theta(\oa{xy})$.
It follows that 
$(1-\eps)\cdot d_{\diamo}(x,y)\geq d(x,y)$,
whence the assertion. 
\qed

\medskip
Whether the modified metric $d_{\diamo}$ can be bounded above 
in terms of the original metric,
depends on the geometry of the face type $\taumod\subset\simod$.
Note that\footnote{$\st(\taumod)$ here 
refers to the star {\em within $\amod$}.}
$$\st(\taumod)=W_{\taumod}\simod\subset\amod$$
is a proper convex subcomplex,
because $X$ has no euclidean factor,
and hence contained in a closed hemisphere.
In fact, it is contained in all closed hemispheres with center in $\taumod$.
(Recall that chambers have diameter $\leq\pihalf$.)
Thus,
$\st(\taumod)$ is itself a hemisphere 
if and only if $\taumod$ is a root type vertex, $\taumod=\{\bar\zeta\}$,
and the spherical Coxeter complex $(\amod,W)$ is reducible 
with the 0-sphere $\{\pm\bar\zeta\}$ as a join factor.
\begin{lem}
\label{lem:modifmetequiv}
If $\st(\taumod)$ is not a closed hemisphere,
then $$d_{\diamo}\leq C\cdot d|_{\diamo}$$
with a constant $C=C(\simod)>1$.
\end{lem}
\proof
We must bound above the $\diamo$-length of longitudinal segments.

Let $P(\tau_-,\tau_+)\supset\diamo$ be an ambient parallel set 
for an $(x_-,x_+)$-opposite pair of simplices $(\tau_-,\tau_+)$ of 
type $(\taumod^-,\taumod^+)$.
Then a longitudinal segment $x'_-x'_+\subset\diamo$ 
is contained in a maximal flat $F\subset P(\tau_-,\tau_+)$. 
Assuming that the segment is oriented so that $x'_{\pm}\in V(x'_{\mp},\ost(\tau_{\pm}))$,
we have 
$\diamo':=\diamot(x'_-,x'_+)\subset\diamo$,
and the ``flat diamond" $\diamo'\cap F$ is the intersection of the two flat sectors 
$V(x'_{\mp},\st(\tau_{\pm}))\cap F=V(x'_{\mp},\st(\tau_{\pm})\cap\geo F)$.
Note that $\tau_{\pm}\subset\geo F$,
so $\st(\tau_{\pm})\cap\geo F\cong\st(\taumod)$.

We will bound $d_{\diamo}(x'_-,x'_+)$ above by connecting the points $x'_{\pm}$
inside $\diamo'\cap F$
by a piecewise non-longitudinal path 
with controlled length.
This can be done by a path $x'_- y x'_+$ in the boundary of $\diamo'\cap F$.
To see this,
choose a pair of antipodes $\zeta_{\pm}\in\tau_{\pm}$
and note that $\st(\tau_{\pm})\subset\bar B(\zeta_{\pm},\pihalf)$ is a proper convex subset.
Accordingly,
the convex subcomplexes $\st(\tau_{\pm})\cap\geo F$ of the apartment $\geo F$
are proper subsets of the ``complementary" closed hemispheres $\bar B(\zeta_{\pm},\pihalf)\cap\geo F$.
Since the open hemispheres are disjoint,
a ray $x'_-\eta_+$ in the boundary of the flat sector 
$V(x'_-,\st(\tau_+)\cap\geo F)$ 
with $\eta_+\in\st(\tau_+)\cap B(\zeta_+,\pihalf)\cap\geo F$
intersects the boundary of the other flat sector 
$V(x'_+,\st(\tau_-)\cap\geo F)$,
and we take $y$ to be the (unique) intersection point.
The path $x_-yx_+$ in the boundary of $\diamo'\cap F$ 
then consists of two non-longitudinal segments
contained in boundaries of $\taumod$-Weyl cones.

To control the length of the path $x'_- y x'_+$,
we note that, 
since there are only finitely many face types $\taumod\subset\simod$,
and hence only finitely many possible isometry types of subcomplexes $\st(\tau_{\pm})\cap\geo F$,
the ray $x'_-\eta_+$ in the boundary of $V(x'_-,\st(\tau_+)\cap\geo F)$ 
can be chosen so that $\tangle(\eta_+,\zeta_+)\leq\pihalf-\de$ 
for a uniform $\de=\de(\simod)>0$.
Then $\angle_y(x_-,x_+)\geq\de$.
The triangle $\De(x'_-,y,x'_+)$ lies in the flat $F$,
and elementary euclidean geometry yields an estimate of the form
$d(x'_-,x'_+)\geq c\cdot(d(x'_-,y)+d(y,x'_+))$
with a constant $c=c(\de)>0$,
and hence 
$d_{\diamo}(x'_-,x'_+)\leq c^{-1}\cdot d(x'_-,x'_+)$.
\qed

\medskip
Thus, under the assumption of the lemma, 
the modified metric $d_{\diamo}$ is uniformly equivalent to the original metric $d$ on $X$,
and in particular it is an honest {\em metric}.
Furthermore, the distortion is small in almost non-longitudinal directions:
\begin{lem}
\label{lem:almlongalmundist}
If $\st(\taumod)$ is not a closed hemisphere 
and if 
$xy\subset\diamo$
is a longitudinal segment 
with direction $\eps$-close to a non-longitudinal direction,
then 
$$d_{\diamo}(x,y)\leq(1+C\eps)\cdot d(x,y)$$
with the constant $C$ from the previous lemma.
\end{lem}
\proof
For the proof, we switch notation (replacing $x,y$ by $x'_{\pm}$) 
and use some of the notation in the proof of the previous lemma.

Suppose that the direction of the longitudinal segment $x'_-x'_+$
is $\eps$-close to a non-longitudinal direction. 
Then a ray $x'_-\eta_+$ in the boundary of $V(x'_-,\st(\tau_+)\cap\geo F)$ 
can be chosen so that $\angle_{x'_-}(x'_+,\eta_+)\leq\eps$.
Let $y'_+\in x'_-\eta_+$ be the point 
with $d(x'_-,y'_+)=d(x'_-,x'_+)$.
Note that the triangle $\De(x'_-,y'_+,x'_+)$ lies in the flat $F$
and that its side $x'_-y'_+$ is non-longitudinal.

If $y'_+\in\diamo$, then we can estimate using the previous lemma:
$$d_{\diamo}(x'_-,x'_+)\leq d_{\diamo}(x'_-,y'_+)+d_{\diamo}(y'_+,x'_+)
\leq d(x'_-,x'_+)+C\cdot d(y'_+,x'_+)
\leq(1+C\eps)\cdot d(x'_-,x'_+)$$
Otherwise, the segment $x'_-y'_+$ leaves $\diamo'$, equivalently, the Weyl cone $V(x'_+,\st(\tau_-))$
in a point $y'$.
Then the path $x'_-y'x'_+$
consists of two non-longitudinal segments 
and, according to the triangle inequality, has length 
$\leq d(x'_-,y'_+)+d(y'_+,x'_+)\leq (1+\eps)\cdot d(x'_-,x'_+)$.
\qed

\medskip
Suppose now 
that $\st(\taumod)$ is a hemisphere,
which means as mentioned above
that $\taumod=\{\bar\zeta\}$ is a root type vertex and $\st(\taumod)=\bar B(\bar\zeta,\pihalf)$,
i.e.\ the spherical Coxeter complex $(\amod,W)$ is reducible 
and splits off the 0-sphere $\{\pm\bar\zeta\}$ as a join factor.

Accordingly,
$\tau_{\pm}=\{\zeta_{\pm}\}$ are antipodal root type vertices,
$\st(\tau_{\pm})=\bar B(\zeta_{\pm},\pihalf)$
and the model space splits off a rank one factor, 
i.e.\ it splits metrically as the product 
$$X\cong T\times X'$$
of a rank one symmetric space or a metric tree $T$ 
and a model space $X'$ of corank one.
The ideal vertices in $\tits T\subset\tits X$ are the type $\bar\zeta$ ideal points.
The type $\taumod$ parallel sets are of the form $l\times X'$ 
for a geodesic line $l\subset T$,
the $\taumod$-Weyl cones of the form $r\times X'$ 
for a geodesic ray $r\subset T$,
and the $\taumod$-diamonds are of the form $s\times X'$ 
for a geodesic segment $s\subset T$.
A segment in $X$ is $\taumod$-regular 
if and only if it is not contained in a cross section
$pt\times X'$.
The longitudinal segments in $\taumod$-diamonds
are precisely the $\taumod$-regular ones. 
Thus, the non-longitudinal paths are precisely the paths contained in cross sections.

It follows that two points in $\diamo$ have finite $\diamo$-distance 
if and only if they lie in the same cross section,
and on cross sections $d_{\diamo}$ coincides with $d$.
In particular,
$d_{\diamo}$ is not an honest metric in this case.

\begin{rem}
Our discussion shows that any two points $x,y\in\diamo$ 
with finite $\diamo$-distance 
can be connected by a polygonal path in $\diamo$ with $d$-length $d_{\diamo}(x,y)$ 
consisting of at most two non-longitudinal segments. 
\end{rem}

\subsection{Contraction properties of nearest point projections in euclidean buildings}
\label{sec:proj}

In this section, let $X$ be a euclidean building without flat factor.

Recall that 
for a closed convex subset $C\subset X$ 
the nearest point projection $\pi_C:X\to C$ is 1-Lipschitz.
In this section, 
we give sharper contraction estimates 
for projections to diamonds.
This is based on the following general observation.
Here, $\inte(\Si_{\bar x}C)$ denotes the interior of $\Si_{\bar x}C$
as a subset of $\Si_{\bar x}X$.
\begin{lem}
\label{lem:arect}
For $\eps>0$ and $A>1$ there exists $R=R(\eps,A)>0$ such that the following holds:

Let $C\subset X$ be closed convex, 
and let $x,y\in X$ be points with projections $\bar x=\pi_C(x)$ and $\bar y=\pi_C(y)$.
Suppose that $d(x,y)<A\cdot d(\bar x,\bar y)$ and $d(x,\bar x)>R\cdot d(\bar x,\bar y)$.
Then the direction $\oa{\bar x\bar y}$ is $\eps$-close to a direction in 
$\Si_{\bar x}C-\inte(\Si_{\bar x}C)$.
\end{lem}
\proof
The assertion is scale invariant and we may therefore assume that 
$d(\bar x,\bar y)=1$.

Since $\bar x=\pi_C(x)$,
we have that 
$\angle_{\bar x}(\oa{\bar xx},\Si_{\bar x}C)\geq\pihalf$.
In particular, 
$\angle_{\bar x}(x,\bar y)\geq\pihalf$ and, analogously, $\angle_{\bar y}(y,\bar x)\geq\pihalf$.
To see that 
these angles can exceed $\pihalf$ 
only by arbitrarily little if $R$ is sufficiently large, 
we proceed as follows 
using triangle comparison.

In order to bound $\cos\angle_{\bar x}(x,\bar y)$ from below,
we divide the quadrilateral $\boxvoid(x,y,\bar y,\bar x)$
into the triangles $\De(\bar x,\bar y,y)$ and $\De(\bar x,x,y)$.
Let $D=d(\bar x,y)$.
Applying comparison to $\De(\bar x,\bar y,y)$ yields 
for the angle $\al=\angle_{\bar x}(\bar y,y)$ that 
$$ \cos\al \geq \frac{1}{D} $$
because $\angle_{\bar y}(\bar x,y)\geq\pihalf$.
And applying comparison to $\De(\bar x,x,y)$ yields 
for the angle $\beta=\angle_{\bar x}(x,y)$ that 
$$ \sin\beta \leq \frac{A}{D} $$
because $d(x,y)\leq A$.
It follows for $\angle_{\bar x}(x,\bar y)\leq\al+\beta$ that 
$$\cos\angle_{\bar x}(x,\bar y) \geq
\cos\al\cos\beta-\sin\al\sin\beta \geq 
\frac{1}{D}\sqrt{1-\frac{A^2}{D^2}}-\frac{A}{D}\sqrt{1-\frac{1}{D^2}}$$
The right-hand side tends $\to0$ as $R\to+\infty$.
Thus, 
$\angle_{\bar x}(x,\bar y)<\pihalf+\eps$ for suitable $R\geq R(\eps,A)$.

Now let $v\in\Si_{\bar x}C$ be the direction 
where the shortest arc in $\Si_{\bar x}X$ connecting $\oa{\bar xx}$ to $\oa{\bar x\bar y}$ 
enters $\Si_{\bar x}C$. 
Since $\angle_{\bar x}(\oa{\bar xx},v)\geq\pihalf$,
it follows that $\angle_{\bar x}(v,\oa{\bar x\bar y})<\eps$.
By its definition, $v\not\in\inte(\Si_{\bar x}C)$.
\qed

\medskip
In the special case of $\taumod$-diamonds $\diamo=\diamot(x_-,x_+)$, the lemma yields,
cf. Lemma~\ref{lem:spdirdiamo}:
\begin{cor}
If $C=\diamo$,
then the direction 
$\oa{\bar x\bar y}$ (in the above lemma) is $\eps$-close to a non-longitudinal direction. 
\end{cor}
We apply this observation to estimate the local contraction of projections. 
The main result of this section is the following estimate 
which is interesting in itself:
\begin{thm}[Contraction estimate]
\label{thm:contrprojdiamo}
For a $\taumod$-diamond $\diamo=\diamot(x_-,x_+)$, 
the map
\begin{equation*}
(X,d)\buildrel \pi_{\diamo}\over\lra(\diamo,d_{\diamo})
\end{equation*}
is locally 1-Lipschitz outside $\diamo$.
\end{thm}
\proof
Suppose that $xy$ is a segment disjoint from $\diamo$,
and let $r>0$ be so small that $xy$ stays outside the $r$-neighborhood of $\diamo$.
We fix constants $\eps,A^{-1}\simeq0$
and subdivide $xy$ into subsegments of length $<R^{-1}r$ with the constant $R=R(\eps,A)$ from 
Lemma~\ref{lem:arect}.
We project the subdivision points to $\diamo$.
According to the corollary, the segments connecting the projections of subsequent subdivision points 
have directions $\eps$-close to non-longitudinal directions
or have length $\leq A^{-1}$ times the length of the corresponding subdivision segment. 

If $\st(\taumod)\subset\amod$ is not a hemisphere, 
then the modified metric $d_{\diamo}$
is uniformly equivalent to the original metric $d|_{\diamo}$ 
and almost undistorted in almost non-longitudinal directions,
cf.\ Lemmas~\ref{lem:modifmetequiv} and~\ref{lem:almlongalmundist}.
Denoting by $\bar c$ the polygonal path in $\diamo$
connecting the projections of the subdivision points,
and by $L$ and $L_{\diamo}$ the lengths measured with respect to the metrics $d$ and $d_{\diamo}$,
we obtain
$$d_{\diamo}(\pi_{\diamo}(x),\pi_{\diamo}(y))\leq
L_{\diamo}(\bar c) \leq(1+C\eps)\cdot L(\bar c)+\frac{C}{A}\cdot d(x,y)
\leq \Bigl(1+C(\eps+\frac{1}{A})\Bigr)\cdot d(x,y)$$
with the constant $C$ of Lemma~\ref{lem:modifmetequiv}.
The assertion follows in this case by letting $\eps,A^{-1}\to0$.

Otherwise,
if $\st(\taumod)$ is a hemisphere,
the assertion becomes trivial:
The euclidean building splits as the product $X\cong T\times X'$
of a metric tree $T$ and a euclidean building $X'$, 
and the $\taumod$-diamonds are of the form $\diamo=s\times X'$
for segments $s\subset T$
(cf.\ above).
The projection has the form 
$\pi_{\diamo}=\pi_s\times\id_{X'}$
with the nearest point projection $\pi_s:T\to s$.
Outside $\diamo$,
the $\pi_s$-component is locally constant,
so $\pi_{\diamo}$ locally maps into one cross section.
The assertion holds, because $d_{\diamo}$ and $d$ agree on cross sections.
\qed

\section{Regular implies Morse}\label{sec:Regular and coarse regular}

In this chapter we prove the main result of this paper, 
the {\em Morse Lemma} for $\taumod$-regular quasigeodesics in model spaces, Theorem \ref{main-tau} in the introduction.

\subsection{Rectifiable paths in euclidean buildings}
\label{sec:rectpaths}

Let $X$ be a euclidean building. 

It is natural to ask
(compare the definitions in sections~\ref{sec:regcreg} and \ref{sec:longi}):
\begin{ques}
\label{ques:regpathbuil}
Are $\taumod$-regular paths in euclidean buildings 
contained in type $\taumod$ parallel sets
as longitudinal paths?
\end{ques}

The goal of this section is to answer the question affirmatively for 
locally rectifiable paths, 
cf.\ Theorem~\ref{thm:rrpathsinflats} below.

\medskip
We begin by discussing basic properties of $\taumod$-regular paths $c:I\to X$.

The segment $c(t)c(s)$ is $\taumod^{\pm}$-regular for $\pm s>\pm t$
and its $\taumod^{\pm}$-direction $\tau_{\pm}(\oa{c(t)c(s)})\subset\Si_{c(t)}X$ at $c(t)$ 
is therefore well-defined.
\begin{lem}
\label{lem:tmdir}
Let $t\in I$ such that $\pm t$ is not maximal in $\pm I$.
Then the $\taumod^{\pm}$-direction $\tau_{\pm}(\oa{c(t)c(s)})$ for $\pm s>\pm t$ 
does not depend on $s$,
i.e.\ there is a well-defined type $\taumod^{\pm}$ simplex $\tau_{\pm}(t)\subset\Si_{c(t)}X$
such that 
$$\tau_{\pm}(\oa{c(t)c(s)})=\tau_{\pm}(t)
\qquad\hbox{for $\pm s>\pm t$}$$ 
\end{lem}
\proof
The direction $\oa{c(t)c(s)}\in\Si_{c(t)}X$ varies continuously with $s$, 
and the type $\taumod^{\pm}$ open stars are the connected components 
of the $\taumod^{\pm}$-regular part of $\Si_{c(t)}X$. 
The direction must therefore remain in the same open star.
\qed

\begin{lem}
\label{lem:opptdir}
For $[a,b]\subset I$ and $t\in(a,b)$, 
the $\taumod^{\pm}$-directions $\tau_{\pm}(t)$ are opposite to each other
if and only if $c(t)\in\diamot(c(a),c(b))$.
\end{lem}
\proof
The segment $c(t)c(a)$ is $\taumod^-$-regular and $c(t)c(b)$ is $\taumod^+$-regular. 
Hence $c(t)$ can only lie in the interior of $\diamot(c(a),c(b))$ and,
according to the description (\ref{eq:intdiamo}) of the interior of diamonds, 
it does so if and only if the $\taumod^{\pm}$-directions 
$\tau_-(\oa{c(t)c(a)})=\tau_-(t)$ and $\tau_+(\oa{c(t)c(b)})=\tau_+(t)$ 
of these segments at $c(t)$ are opposite to each other.
\qed

\medskip
As a consequence of these two lemmas,
we obtain a ``local-global" equivalence for the straightness of triples on the path:
\begin{cor}
\label{cor:diamcontequiv}
For $[a',b']\subset[a,b]\subset I$ and $t\in(a',b')$,
we have:
$$c(t)\in\diamot(c(a),c(b))\Leftrightarrow c(t)\in\diamot(c(a'),c(b'))$$
\end{cor}
Another useful consequence is:
\begin{cor}
\label{cor:intpairlongi}
If $[a',b']\subset[a,b]\subset I$ 
and if $c(a'),c(b')\in\diamo=\diamot(c(a),c(b))$,
then the segment $c(a')c(b')$ is longitudinal in $\diamo$.
\end{cor}
\proof
With the last lemma,
we see that 
the pair of $\taumod^{\pm}$-directions $\tau_{\pm}(a')$ is opposite,
as well as the pair $\tau_{\pm}(b')$.
The assertion then follows e.g.\ from Corollary~\ref{cor:longiconvdiamo}.
\qed

\medskip
The main result of this section is a positive answer to Question~\ref{ques:regpathbuil} 
for arbitrary (i.e.\ possibly non-discrete) euclidean buildings 
in the rectifiable uniform case:
\begin{thm}
\label{thm:rrpathsinflats}
Let $X$ be a euclidean building.
Then every rectifiable $\Theta$-regular path $c:[a,b]\to X$
is contained in the $\Theta$-diamond $\diamo=\diamotTh(c(a),c(b))$ spanned by its endpoints
and is longitudinal in $\diamo$.
\end{thm}
We break the proof up into several steps.

We observe first that, arguing by contradiction,
we may assume that the path does not touch the diamond at all except at its endpoints:
\begin{lem}
\label{lem:subpouts}
If the path 
$c:[a,b]\to X$ is $\taumod$-regular 
and not contained in $\diamo=\diamot(c(a),c(b))$,
then there exists a nondegenerate subinterval $[a',b']\subset[a,b]$
such that $c(t)\not\in\diamot(c(a'),c(b'))$ for all $t\in(a',b')$. 
\end{lem}
\proof
Denote $\diamo=\diamot(c(a),c(b))$,
and let $(a',b')$ be a connected component of the nonempty open subset 
$\{t\in I:c(t)\not\in\diamo\}$.
Then $c(a'),c(b')\in\diamo$ 
and, invoking Corollary~\ref{cor:intpairlongi},
we know that the segment $c(a')c(b')$ is longitudinal,
and hence $\diamot(c(a'),c(b'))\subseteq\diamo$.
\qed

\medskip
The intuition behind the proof of the theorem is that 
it ``costs length'' for a $\taumod$-regular path 
to move outside the diamond of its endpoints,
due to the contraction properties of projections as described in section~\ref{sec:proj}.
The key step in the proof is:
\begin{lem}
Let $c:[a,b]\to X$ be a 
path such that the oriented segment $c(a)c(b)$ connecting its endpoints is $\taumod$-regular 
and such that 
$c(t)\not\in\diamo=\diamot(c(a),c(b))$ for all $t\in(a,b)$. 
Then
$$L(c)\geq d_{\diamo}(c(a),c(b)).$$
\end{lem}
\proof
Let $\pi_{\diamo}:X\to\diamo$ denote the nearest point projection.
We consider the projected curve $\bar c:=\pi_{\diamo}\circ c:I\to\diamo$. 
Since $c$ lies outside $\diamo$ except for its endpoints,
Theorem~\ref{thm:contrprojdiamo} yields that 
$$L(c)\geq L_{\diamo}(\bar c)$$
where $L_{\diamo}$ is the length measured with respect to the modified metric $d_{\diamo}$. 
The assertion follows, because $L_{\diamo}(\bar c)\geq d_{\diamo}(c(a),c(b))$.
\qed

\medskip
As a consequence, 
almost distance minimizing 
(cf.\ Definition ~\ref{dfn:almdmin})
uniformly $\taumod$-regular paths 
must touch the diamond of their endpoints:
\begin{lem}
\label{lem:dmintdiamo}
There exists $\eps=\eps(\Theta)>0$ such that the following holds:
If the path $c:[a,b]\to X$ is $\Theta$-regular and $\eps$-distance minimizing, 
then $c(t)\in\diamo$ for some $t\in(a,b)$. 
\end{lem}
\proof
Due to the compactness of $\Theta$,
there exists an $\eps=\eps(\Theta)>0$ such that $d_{\diamo}>(1+\eps)\cdot d$ for $\Theta$-longitudinal 
pairs of points in $\diamo$,
cf.\ Lemma~\ref{lem:horlong}.
The assertion then follows from the previous lemma.
\qed

\medskip
Based on Corollary~\ref{cor:diamcontequiv},
we can extend the last result to rectifiable paths of arbitrary length,
because these contain arbitrarily distance minimizing subpaths:
\begin{lem}
\label{lem:intdiamo}
If $c:[a,b]\to X$ is $\Theta$-regular and rectifiable, 
then $c(t)\in\diamo$ for some $t\in(a,b)$. 
\end{lem}
\proof
Let $\eps=\eps(\Theta)$ be the constant from Lemma~\ref{lem:dmintdiamo}.
There exists a nondegenerate subinterval $[a',b']\subset[a,b]$
such that the subpath $c|_{[a',b']}$ is 
$\eps$-distance minimizing, 
cf.\ Lemma~\ref{lem:distminsubp}.
Then Lemma~\ref{lem:dmintdiamo} implies that 
$c(t)\in\diamot(c(a'),c(b'))$ for some $t\in(a',b')$.
With Corollary~\ref{cor:diamcontequiv},
it follows that also $c(t)\in\diamo$.
\qed

\medskip
We are ready to conclude the proof of the theorem. 

\noindent 
{\em Proof of Theorem~\ref{thm:rrpathsinflats}.}
Suppose that $c$ is not contained in $\diamotTh(c(a),c(b))$.
Then, by $\Theta$-regularity, 
it is also not contained in $\diamo=\diamot(c(a),c(b))$.
According to Lemma~\ref{lem:subpouts},
after replacing $c$ by a subpath,
we may assume that $c(t)\not\in\diamo$ for all $t\in(a,b)$. 
But this contradicts Lemma~\ref{lem:intdiamo},
so $c$ is contained in $\diamotTh(c(a),c(b))$.
It is longitudinal by Corollary~\ref{cor:intpairlongi}.
\qed

\medskip
As a consequence of the theorem, we obtain with Lemma~\ref{lem:bdddet}
that rectifiable uniformly $\taumod$-regular paths are, 
up to reparametrization, bilipschitz;
they become bilipschitz when parametrized by arc length:
\begin{cor}[Bounded detours]
\label{cor:regpundist}
$L(c)\leq L(\Theta)\cdot d(c(a),c(b))$
\end{cor}

We have the following implications of the theorem for infinite paths:
\begin{cor}
\label{cor:rectlocunifregp}
(i) Every locally rectifiable and uniformly $\taumod$-regular path $c:I\to X$ 
is contained in a type $\taumod$ parallel set as a longitudinal path.

(ii) Every locally rectifiable and 
uniformly $\taumod$-regular path $c:[0,+\infty)\to X$ with infinite length
is contained as a longitudinal path in a Weyl cone $V(c(0),\st(c(+\infty))$ 
for a unique simplex $c(+\infty)\in\Flagt(\geo X)$. 

(iii) Every locally rectifiable and 
uniformly $\taumod$-regular path $c:\R\to X$,
both of whose ends have infinite length,
is contained as a longitudinal path in a parallel set 
$P(c(-\infty),c(+\infty))$ 
for a unique pair of opposite simplices $c(\pm\infty)\in\Flag_{\taumod^{\pm}}(\geo X)$. 
\end{cor}
\proof
(i) By Theorem~\ref{thm:rrpathsinflats}, 
for every compact subinterval $[a,b]\subset I$, 
the corresponding part of the path is contained in the diamond $\diamot(c(a),c(b))$.
These diamonds are nested, i.e. for 
$[a',b']\subset[a,b]\subset I$,
it holds that $\diamot(c(a'),c(b'))\subset\diamot(c(a),c(b))$.
Since the path $c$ is uniformly $\taumod$-regular,
the closure of the union of these diamonds over all compact subintervals of $I$ 
is either a type $\taumod$ diamond, a Weyl cone or a parallel set. 
The longitudinality follows from the longitudinality part of the theorem. 

(ii)
The sequence $(c(n))_{n\in\N_0}$ is $\Theta$-regular 
and diverges to infinity in view of Corollary~\ref{cor:regpundist}.
By the theorem,
$c([0,n])\subset\diamoTh(c(0),c(n))$ for all $n$.
In this situation, Lemma~\ref{lem:taunconv} applies,
after enclosing the diamonds $\diamot(c(0),c(n))$ into auxiliary Weyl cones $V(c(0),\st(\tau_n))$,
and yields that the sequence $(c(n))$ is contained in a Weyl cone $V(c(0),\st(c(+\infty))$
for a simplex $c(+\infty)\in\Flagt(\geo X)$,
which is unique according to Lemma~\ref{lem:flagconv}.
Since $V(c(0),\st(c(+\infty))$ then also contains the diamonds $\diamot(c(0),c(n))$,
it contains the entire path $c$.

(iii) By part (ii), there exist unique simplices $c(\pm\infty)\in\Flag_{\taumod^{\pm}}(\geo X)$
so that $c(\pm[t,+\infty))\subset V(c(\pm t),\st(c(\pm\infty)))$.
It follows that for any $t_1<t_2$ and any ideal points $\xi_{\pm}\in\ost(c(\pm\infty))$, 
the biinfinite broken path $\xi_-c(t_1)c(t_2)\xi_+$ is $\taumod$-straight.
Proposition~\ref{prop:longiconvpar} then implies that 
the simplices $c(\pm\infty)$ are opposite and the path $c$ is contained in the parallel set 
$P(c(-\infty),c(+\infty))$ as a longitudinal path.
\qed

\begin{dfn}[Endpoint at infinity]
For a locally rectifiable and uniformly $\taumod$-regular path $c:[0,+\infty)\to X$ with infinite length,
we call the simplex $c(+\infty)\in\Flagt(\geo X)$
its {\em $\taumod$-endpoint at infinity} or {\em ideal $\taumod$-endpoint}.
\end{dfn}

\medskip
We apply our results to paths of infinite length which remain close to a Weyl cone or a parallel set.
In some situations one can show that they must be contained in it.
\begin{cor}
\label{cor:blrclwco}
Let $c:[0,+\infty)\to X$ be a $\Theta$-regular $L$-bilipschitz ray 
which is contained in the tubular $D$-neighborhood of a type $\taumod$ Weyl cone $V$
with tip at $c(0)$.
Then $V=V(c(0),\st(c(+\infty)))$
and $c$ is contained in $V$ as a longitudinal path.
\end{cor}
\proof
Suppose that 
$c([0,+\infty))\subset\bar N_D(V(c(0),\st(\tau_+)))$ for a simplex $\tau_+\in\Flagt(\geo X)$ 
and some $D>0$.
According to Corollary~\ref{cor:rectlocunifregp},
$c$ is contained in the Weyl cone $V(c(0),\st(c(+\infty)))$ as a longitudinal path.
The sequence $(c(n))_{n\in\N}$ is $\Theta$-regular. 
Corollary~\ref{cor:closewco} therefore implies that $\tau_+=c(+\infty)$.
\qed

\medskip
If a biinfinite $\taumod$-regular bilipschitz path is close to a type $\taumod$ parallel set,
we need a longitudinality property for its projection
to be able to conclude that it must be contained in the parallel set.
We denote by $\bar c=\pi_P\circ c$ 
the projection of the path $c$ to the parallel set $P$.

\begin{cor}
\label{cor:bllclps}
There exists a constant $l=l(L,\Theta,D)>0$ such that the following holds:

Suppose that $c:\R\to X$ is a $\Theta$-regular $L$-bilipschitz line
which is contained in the tubular $D$-neighborhood 
of a type $\taumod$ parallel set $P$.
If for some interval $[a',b']\subset\R$ of length $\geq l$
the segment $\bar c(a')\bar c(b')\subset P$ is longitudinal,
then $P=P(c(-\infty),c(+\infty))$
and $c$ is contained in $P$ as a longitudinal path.
\end{cor}
\proof
Suppose that $P=P(\tau_-,\tau_+)$
with opposite simplices $\tau_{\pm}\in\Flag_{\taumod^{\pm}}(\geo X)$, 
and $c(\R)\subset\bar N_D(P)$ for some $D>0$.
The projection $\bar c=\pi_P\circ c$ is coarsely longitudinal by Lemma~\ref{lem:clqg}.
More precisely, 
we choose $\Theta'$ depending on $\Theta$ 
and put $l=LcD$
with the constant $c=c(\Theta,\Theta')>0$ from Lemma~\ref{lem:lcoaregureg}.
Then Lemma~\ref{lem:clqg} yields 
that for {\em all} subintervals $[a',b']\subset\R$ of length $\geq l$
the segment $\bar c(a')\bar c(b')\subset P$ is longitudinal.
It follows for the bilipschitz rays $r_+=c|_{[0,+\infty)}$ and $r_-=c|_{(-\infty,0]}$ 
that $\pi_P\circ r_{\pm}$ 
is contained in a tubular neighborhood of the Weyl cone 
$V(\bar c(0),\st(\tau_{\pm}))\subset P$,
and hence $r_{\pm}$ in a tubular neighborhood of $V(\bar c(0),\st(\tau_{\pm}))$.
By Corollary~\ref{cor:blrclwco}, $\tau_{\pm}=c(\pm\infty)$.
According to Corollary~\ref{cor:rectlocunifregp}, 
$c$ is a longitudinal path in $P(c(-\infty),c(+\infty))=P$.
\qed

\begin{rem}
For {\em discrete} buildings,
the answer to Question~\ref{ques:regpathbuil} is affirmative 
without restriction on the paths.
Using that discrete euclidean buildings are locally conical, 
it is not hard to show that 
every $\taumod$-regular path $c:[a,b]\to X$ in a discrete euclidean building $X$ 
is contained in $\diamot(c(a),c(b))$.
\end{rem}

\subsection{The Morse Lemma for quasigeodesics in CAT(0) model spaces}
\label{sec:morselem}

We recall that 
the Morse Lemma for quasigeodesics in Gromov hyperbolic spaces asserts 
that uniform quasigeodesics are uniformly close to geodesics.
The main result of this paper is the following generalization to model spaces of arbitrary rank,
where geodesic lines (rays, segments) are replaced by parallel sets (cones, diamonds):
\begin{thm}[Morse Lemma]
\label{thm:mlemt}
Let $X$ be a model space.
Suppose that $q:[a_-,a_+]\to X$ is a $(\Theta,B)$-regular $(L,A)$-qua\-si\-geo\-de\-sic 
and that $x_-x_+$ is a $\Theta$-regular segment oriented $B$-Hausdorff close to $q(a_-)q(a_+)$.
Then the image of $q$ is contained in the $D$-neighborhood 
of the diamond $\diamot(x_-,x_+)$,
with a constant $D=D(L,A,\Theta,B,X)>0$.
\end{thm}
\proof
We will deduce the theorem 
from the corresponding result for bilipschitz paths in euclidean buildings,
cf.\ Theorem~\ref{thm:rrpathsinflats},
by passing to ultralimits.
We may work without loss of generality with {\em continuous} quasigeodesics.

We argue by contradiction. 
Suppose that a uniform constant $D$ does not exist
and consider,
for a fixed model space $X$ and fixed data $(L,A,\Theta,B)$,
sequences of $(\Theta,B)$-regular $(L,A)$-quasigeodesics 
$$q_n:I_n=[a_n^-,a_n^+]\to X,$$
of $\Theta$-regular segments $x_n^-x_n^+$
oriented $B$-Hausdorff close to the segments $q_n(a_n^-)q_n(a_n^+)$,
and of positive numbers $D_n\to+\infty$,
such that the image of $q_n$ is contained in the $D_n$-neigh\-bor\-hood of 
$$\diamon:=\diamot(x_n^-,x_n^+) ,$$
but not in its $\frac{2013}{2014}D_n$-neighborhood.
We may assume that 
$a_n^-\leq0\leq a_n^+$
and that $q_n(0)$ has almost maximal distance $>\frac{2013}{2014}D_n$ from $\diamon$.
Note that $\liminf_nD_n^{-1}|a_n^{\pm}|>0$.

Let $P_n=P(\tau_n^-,\tau_n^+)\subset X$ be a type $\taumod$ parallel set 
through the points $x_n^{\pm}$ 
such that the segment $x_n^-x_n^+$ is longitudinal. 
Then 
$$\diamon=V(x_n^-,\st(\tau_n^+))\cap V(x_n^+,\st(\tau_n^-))\subset P_n$$
and the image of $q_n$ is contained in the $D_n$-neighborhood of $P_n$.

The next result provides important information 
on the position of the quasigeodesics $q_n$ relative to the parallel sets $P_n$
if the length of $q_n$ grows faster than the scale $D_n$.

Let $\bar q_n=\pi_{P_n}\circ q_n$ denote 
the nearest point projection of $q_n$ to $P_n$ .
We fix some $\Theta'$.
(As usual, $\Theta'$ is supposed to contain $\Theta$ in its interior.)
\begin{lem}[Coarsely longitudinal on scale $D_n$]
\label{lem:qnclongi}
For every subinterval $[b_n^-,b_n^+]\subset I_n$ of length 
$\geq L(A+c(B+D_n))$
the segment $\bar q_n(b_n^-)\bar q_n(b_n^+)\subset P_n$ is $\Theta'$-longitudinal,
where $c=c(\Theta,\Theta')>0$ is the constant from Lemma~\ref{lem:lcoaregureg}.
\end{lem}
\proof
This is a direct consequence of Lemmas~\ref{lem:cregqgeolongsegmregproj} and~\ref{lem:clqg},
because the segment $\bar q_n(a_n^-)\bar q_n(a_n^+)$ 
is longitudinal by the choice of $P_n$.
\qed

\medskip
Now we pass to the ultralimit.

We choose base points $\star_n\in\diamon$ with $d(q_n(0),\star_n)\leq D_n$,
rescale (copies of) the space 
$X$ with the scale factors $D_n^{-1}\to0$
and then take the ultralimit (with respect to some nonprincipal ultrafilter $\om$).
As proven in \cite[ch.\ 5]{qirigid},
the ultralimit of rescaled model spaces
$$(X_{\om},\star_{\om})=\ulim_n(D_n^{-1}X,\star_n)$$
is a euclidean building of the same type $\simod$,
cf.\ section~\ref{sec:ulimps}.
The ultralimit of parallel sets 
$$P_{\om}:=\ulim_n D_n^{-1}P_n\subset X_{\om}$$
is again a type $\taumod$ parallel set,
$$P_{\om}=P(\tau_{\om}^-,\tau_{\om}^+)$$
for a pair of opposite type $\taumod^{\pm}$ simplices $\tau_{\om}^{\pm}\subset\geo X_{\om}$, 
cf.\ Lemma~\ref{lem:ulimps}.
The ultralimit of diamonds
$$\diamoom:=\ulim D_n^{-1}\diamon\subset P_{\om}$$
is a closed convex subset which contains the base point $\star_{\om}$.
It is in general not a diamond,
but it inherits the following geometric property from the diamonds $\diamo_n$:
\begin{lem}
\label{lem:whlongindi}
If the segment $y_{\om}^-y_{\om}^+\subset\diamoom$ is longitudinal, 
then $\diamot(y_{\om}^-,y_{\om}^+)\subset\diamoom$.
\end{lem}
\proof
The segment $y_{\om}^-y_{\om}^+$ is the ultralimit of segments $y_n^-y_n^+\subset\diamo_n$,
and these segments are longitudinal for $\om$-all $n$. 
Hence $\diamot(y_n^-,y_n^+)\subset\diamo_n$
due to Lemma~\ref{lem:nestdiamo}.
With Lemma~\ref{lem:uconvdiamo} 
it follows that
 $\diamot(y_{\om}^-,y_{\om}^+)=\ulim D_n^{-1}\diamot(y_n^-,y_n^+) \subset\diamoom$.
\qed

\medskip
The rescaled paths 
$D_n^{-1}q_n:D_n^{-1}I_n\to D_n^{-1}X$ given by $$(D_n^{-1}q_n)(t_n)=q_n(D_nt_n)$$
are $(\Theta,D_n^{-1}B)$-regular $(L,D_n^{-1}A)$-quasigeodesics.
Their ultralimit 
$q_{\om}=\ulim_nD_n^{-1}q_n:I_{\om}\to X_{\om}$ given by 
$$q_{\om}(t_{\om})=(q_n(D_nt_n))$$
is a well-defined {\em $L$-bilipschitz} path because $D_n^{-1}A\to0$,
cf.\ Lemma~\ref{lem:ulimqiebil},
and {\em $\Theta$-regular} because $D_n^{-1}B\to0$.
Its domain is the interval $I_{\om}=[a_{\om}^-,a_{\om}^+]\cap\R$,
where $a_{\om}^{\pm}=\ulim D_n^{-1}a_n^{\pm}$ 
and $\pm a_{\om}^{\pm}\in(0,+\infty]$.
If $|a_{\om}^{\pm}|<+\infty$,
then the endpoint 
$$q_{\om}(a_{\om}^{\pm})=\ulim q_n(a_n^{\pm})=\ulim x_n^{\pm}=:x_{\om}^{\pm}$$
exists
and lies in $\diamoom$.
Otherwise, the corresponding end of $q_{\om}$ has infinite length and diverges to infinity.

By construction, 
$$q_{\om}(I_{\om})\subset\bar N_1(\diamoom),$$
but
$$q_{\om}(I_{\om})\not\subset\diamoom.$$
In particular,
$q_{\om}(I_{\om})\subset\bar N_1(P_{\om})$,
and we denote by $\bar q_{\om}=\pi_{P_{\om}}\circ q_{\om}$ the nearest point projection.
Then $d(\bar q_{\om},q_{\om})\leq1$.

Regarding the position of $q_{\om}$ relative to the parallel set, 
it inherits from the $q_n$
uniform longitudinality beyond a certain scale:
\begin{lem}[Coarsely longitudinal ultralimit]
\label{lem:qomclongimod}
For every subinterval $[b_-,b_+]\subset I_{\om}$ of length 
$\geq\frac{2014}{2013}cL$
the segment $\bar q_{\om}(b_-)\bar q_{\om}(b_+)\subset P_{\om}$ is $\Theta'$-longitudinal.
\end{lem}
\proof
Apply Lemma~\ref{lem:qnclongi} taking into account that $D_n^{-1}A,D_n^{-1}B\to0$.
\qed

\medskip
If $q_{\om}$ has infinite length,
then the coarse longitudinality restricts the asymptotics of its end(s);
they must flag converge to the simplices $\tau_{\om}^{\pm}$.
We get the following information on $\diamoom$:
\begin{lem}
If $|a_{\om}^{\pm}|=+\infty$, then $V(\star_{\om},\st(\tau_{\om}^{\pm}))\subset\diamoom$.
\end{lem}
\proof
Suppose that $a_{\om}^+=+\infty$.
We have that $d(\pi_{\diamoom}\circ q_{\om},\bar q_{\om})\leq2$.
Lemma~\ref{lem:qomclongimod} therefore implies that 
the segment connecting $\star_{\om}$ 
to the point $\pi_{\diamoom}(q_{\om}(t))\in\diamoom$
is longitudinal for all sufficiently large $t>0$,
i.e. $\pi_{\diamoom}(q_{\om}(t))\in V(\star_{\om},\st(\tau_{\om}^+))$.
In particular,
for any sequence $t_k\to+\infty$
it holds that $\pi_{\diamoom}(q_{\om}(t_k))\to\tau_{\om}^+$
(equivalently, $q_{\om}(t_k)\to\tau_{\om}^+$) as $k\to+\infty$,
even conically, 
cf. Lemma~\ref{lem:flagconv}.
The longitudinality of the segments implies furthermore that 
$$\diamot(\star_{\om},\pi_{\diamoom}(q_{\om}(t)))\subset\diamoom,$$
cf. Lemma~\ref{lem:whlongindi}.
The assertion follows now with the description of flag convergence in euclidean buildings
given in Lemma~\ref{lem:txconveb}.
The case $a_{\om}^-=-\infty$ is analogous.
\qed

\medskip
This allows us to classify the possibilities for $\diamoom$:
\begin{cor}
$\diamoom$ either equals the diamond $\diamot(x_{\om}^-,x_{\om}^+)$, 
or one of the two Weyl cones $V(x_{\om}^{\mp},\st(\tau_{\om}^{\pm}))$,  
or the full parallel set $P_{\om}$, 
depending on whether both, one or none of the points $x_{\om}^{\mp}$ are defined.
\end{cor}
\proof
If the endpoint $x_{\om}^{\pm}$ of $q_{\om}$ exists,
then clearly 
$\diamoom\subset V(x_{\om}^{\pm},\st(\tau_{\om}^{\mp}))$,
because $\diamo_n\subset V(x_n^{\pm},\st(\tau_n^{\mp}))$.
If it does not exist,
then $V(\star_{\om},\st(\tau_{\om}^{\pm}))\subset\diamoom$
by the previous lemma.

Thus, if none of the endpoints exists, then $\diamoom=P_{\om}$ by convexity.
And, if exactly one endpoint $x_{\om}^{\pm}$ exists, 
then $\diamoom=V(x_{\om}^{\pm},\st(\tau_{\om}^{\mp}))$,
also by convexity.
If both endpoints exist, then 
$\diamoom\subset\diamot(x_{\om}^-,x_{\om}^+)$,
and equality follows from Lemma~\ref{lem:uconvdiamo},
cf. also Lemma~\ref{lem:whlongindi}.
\qed

\medskip
Now we apply our results on rectifiable regular paths from section~\ref{sec:rectpaths} to $q_{\om}$
in order to control its position also on the small scale:
\begin{lem}
\label{lem:qomindiamo}
$q_{\om}(I_{\om})\subset\diamoom$.
\end{lem}
\proof
If both endpoints $x_{\om}^{\pm}$ of $q_{\om}$ exist,
then $\diamoom=\diamot(x_{\om}^-,x_{\om}^+)$
and Theorem~\ref{thm:rrpathsinflats} implies the assertion. 

In the other cases, we use that $q_{\om}(I_{\om})\subset\bar N_1(\diamoom)$.

If $q_{\om}$ has exactly one endpoint, say $x_{\om}^-$, 
and thus is a bilipschitz {\em ray},
then $\diamoom=V(x_{\om}^-,\st(\tau_{\om}^+))$
and Corollary~\ref{cor:blrclwco} implies the assertion. 

If $q_{\om}$ has no endpoints at all and thus is a bilipschitz {\em line},
then $\diamoom=P_{\om}$.
We use that the projection $\bar q_{\om}$ to $P_{\om}$ is coarsely longitudinal,
cf.\ Lemma~\ref{lem:qomclongimod}.
We therefore can apply Corollary~\ref{cor:bllclps} which yields the assertion in this case.
\qed

\medskip
The last lemma contradicts that $q_{\om}(I_{\om})\not\subset\diamoom$.
This concludes the proof of Theorem~\ref{thm:mlemt}.
\qed

\begin{rem}
\label{rem:aftmthm}
(i)
It follows moreover that the projection $\bar q=\pi_{\diamo}\circ q$ of $q$ 
to the diamond $\diamo=\diamot(x_-,x_+)$
is {\em coarsely longitudinal},
by which we mean that for every subinterval $[b_-,b_+]\subset[a_-,a_+]$ of length $\geq L(A+c(B+D))$
the segment $\bar q(b_-)\bar q(b_+)\subset\diamo$ is $\Theta'$-regular and longitudinal,
with a constant $c=c(\Theta,\Theta')>0$.
This is a consequence of Lemma~\ref{lem:clqg} and was used in the proof of the theorem, 
compare Lemma~\ref{lem:qnclongi}.

(ii)
In the building case,
the argument works equally well
if we replace $X$ by a sequence of euclidean buildings $X_n$ of fixed type $\simod$.
Hence, the bound for the size of the tubular neighborhood depends only on the {\em rank} of the euclidean building 
and not on further geometric properties of it,
$D=D(L,A,\taumod,\Theta,B,\rank(X))$.

(iii)
The theorem remains valid 
if one allows the model spaces to have flat factors, 
because the case with flat factors immediately reduces to the case without.
\end{rem}

We have the following implications of the theorem for infinite quasigeodesics:
\begin{cor}
\label{cor:infqgclcp}
(i) Suppose that $q:[0,+\infty)\to X$ is a $(\Theta,B)$-regular $(L,A)$-quasiray.
Then the image of $q$ is contained in the $(D+B)$-neighborhood 
of the Weyl cone $V(q(0),\st(q(+\infty))$ 
for a unique simplex $q(+\infty)\in\Flagt(\geo X)$.

(ii) Suppose that $q:\R\to X$ is a $(\Theta,B)$-regular $(L,A)$-quasiline.
Then the image of $q$ is contained in the $(D+B)$-neighborhood 
of the parallel set $P(q(-\infty),q(+\infty))$
for a unique pair of opposite simplices $q(\pm\infty)\in\Flag_{\taumod^{\pm}}(\geo X)$.

In both cases, 
$q$ is coarsely longitudinal in the sense of Remark~\ref{rem:aftmthm}. 
\end{cor}
\proof
(i) Let $p_nx_n$ be $\Theta$-regular segments 
oriented $B$-Hausdorff close to the segments $q(0)q(n)$.
According to the theorem,
$q([0,n])\subset\bar N_D(\diamot(p_n,x_n))$.
We extend the diamonds to cones,
i.e. we let $\tau_n\in\Flagt(\geo X)$ be simplices 
such that 
$\diamot(p_n,x_n)\subset V(p_n,\st(\tau_n))$.
Then $q([0,n])\subset\bar N_D(V(p_n,\st(\tau_n)))$, 
and hence $q([0,n])\subset\bar N_{D+B}(V(q(0),\st(\tau_n)))$.
The sequence $q(n)\to\infty$ is asymptotically uniformly $\taumod$-regular,
because the quasiray $q$ is $(\Theta,B)$-regular.
Applying the convergence criterion in Lemma~\ref{lem:taunconv},
it follows that the sequence $(\tau_n)$ converges, 
$\tau_n\to\tau_{\infty}\in\Flagt(\geo X)$,
and the image of $q$ is contained in the $(D+B)$-neighborhood of $V(q(0),\st(\tau_{\infty}))$.
We put $q(+\infty)=\tau_{\infty}$.
The uniqueness follows from Corollary~\ref{cor:closewco}.

(ii) According to part (i),
there exist unique simplices $\tau_{\pm\infty}=q(\pm\infty)\in\Flag_{\taumod^{\pm}}(\geo X)$
such that 
$q(\pm n)\to\tau_{\pm\infty}$ conically as $n\to+\infty$.
More precisely,
\begin{equation}
\label{eq:cltwco}
q(\pm[-n,+\infty))\subset\bar N_{D+B}(V(q(\mp n),\st(\tau_{\pm\infty}))).
\end{equation}
The segment $q(-n)q(n)$ is $\Theta'$-regular and arbitrarily long for large $n$.
Let $y_{-n}y_n$ be a subsegment of it
at distance $>(D+B)\cdot(\sin\eps_0(\Theta'))^{-1}$ from the endpoints $q(\pm n)$.
By Corollary~\ref{cor:leavecone},
$$y_{-n}y_n\subset V(q(-n),\st(\tau_{+\infty}))\cap V(q(n),\st(\tau_{-\infty})).$$
Then for any interior point $z_n$ of this segment, 
it holds that 
$\log_{z_n}\tau_{\pm\infty}=\tau_{\pm}(z_nq(\pm n))$,
and it follows that the simplices $\tau_{\pm\infty}$ are $z_n$-opposite. 
Furthermore, 
the $\Theta'$-cones $V(q(t),\st(\Theta'))$ enter the parallel set $P=P(\tau_{-\infty},\tau_{+\infty})$
within uniformly bounded time, cf. Proposition~\ref{prop:diveparbuil}, 
and in view of (\ref{eq:cltwco}) it follows that $\bar N_{D+B}(P)$ contains the image of $q$.

The coarse longitudinality is a consequence of Lemma~\ref{lem:clqg}.
\qed

\begin{rem}
As a consequence of the corollary,
the image of every uniformly coarsely $\taumod$-regular uniform quasiline
is contained in a uniform neighborhood of a union
of two opposite Weyl cones in a parallel set.
The cones have a common tip
which can be chosen uniformly close to any point on the quasiline. 
\end{rem}

\begin{dfn}[Endpoint at infinity]
\label{dfn:endqr}
For a $(\Theta,B)$-regular quasiray $q:[0,+\infty)\to X$ 
we call the simplex $q(+\infty)\in\Flagt(\geo X)$
its {\em $\taumod$-endpoint at infinity} or {\em ideal $\taumod$-endpoint}.
\end{dfn}

We apply our results to infinite quasigeodesics which remain close to a Weyl cone or a parallel set
and show that they must be uniformly close.
\begin{cor}
\label{cor:qrclwco}
Let $q:[0,+\infty)\to X$ be a $(\Theta,B)$-regular $(L,A)$-quasiray
which is contained in a tubular neighborhood of a type $\taumod$ Weyl cone $V$
with tip at $q(0)$.
Then $V=V(q(0),\st(q(+\infty)))$
and $q$ is contained in the tubular $(D+B)$-neighborhood of $V$ 
as a coarsely longitudinal path.
\end{cor}
\proof
Suppose that 
$q([0,+\infty))\subset N_r(V(q(0),\st(\tau_+)))$ for a simplex $\tau_+\in\Flagt(\geo X)$
and some $r>0$.
According to Corollary~\ref{cor:infqgclcp},
$q$ is contained in the $(D+B)$-neighborhood of the Weyl cone $V(q(0),\st(q(+\infty)))$ 
as a coarsely longitudinal path.
The sequence $(q(n))_{n\in\N}$ is asymptotically uniformly $\taumod$-regular. 
Corollary~\ref{cor:closewco} therefore implies that $\tau_+=q(+\infty)$.
\qed

\medskip
If a coarsely $\taumod$-regular quasline is close to a type $\taumod$ parallel set,
we need a coarse longitudinality property 
to be able to conclude that it must be uniformly close to the parallel set.
We denote by $\bar q=\pi_P\circ q$ 
the projection of $q$ to the parallel set $P$.
\begin{cor}
\label{cor:bllclps2}
There exists a constant $l=l(L,A,\Theta,B,r)>0$ such that the following holds:

Suppose that $q:\R\to X$ is a $(\Theta,B)$-regular $(L,A)$-quasiline
which is contained in a tubular neighborhood 
of a type $\taumod$ parallel set $P$. 
If for some interval $[a,b]\subset\R$ of length $\geq l$
the segment $\bar q(a)\bar q(b)\subset P$ is longitudinal,
then $P=P(q(-\infty),q(+\infty))$
and $q$ is contained in $\bar N_{D+B}(P)$
as a coarsely longitudinal path.
\end{cor}
\proof
Suppose that $P=P(\tau_-,\tau_+)$ with opposite simplices $\tau_{\pm}\in\Flag_{\taumod^{\pm}}(\geo X)$,
and $q(\R)\subset\bar N_r(P)$ for some $r>0$.
The projection $\bar q=\pi_P\circ q$ is coarsely longitudinal along $P$ by Lemma~\ref{lem:clqg}.
It follows with Corollary~\ref{cor:qrclwco}
that $\tau_{\pm}=q(\pm\infty)$.
According to Corollary~\ref{cor:infqgclcp},
$q$ is contained as a coarsely longitudinal path in the $(D+B)$-neighborhood 
of $P(q(-\infty),q(+\infty))=P$.
\qed

\subsection{Regular implies Morse for undistorted maps and actions}

We relate the Morse Lemma (Theorem~\ref{thm:mlemt})
to terminology used in our paper \cite{morse}.

There we defined (in the setting of symmetric spaces)
a Morse quasigeodesic as a quasigeodesic satisfying 
the conlusion of the Morse Lemma with $\taumod$-diamonds replaced by $\Theta$-diamonds,
i.e. every finite subpath of the quasigeodesic is uniformly close to a diamond
whose tips are uniformly close to the endpoints of the subpath.
More precisely:
\begin{dfn}[Morse quasigeodesic, cf. {\cite[Definition 7.14]{morse}}]
\label{dfn:mqg}
An {\em $(L,A,\Theta,D)$-Morse quasigeodesic} in $X$ 
is an $(L,A)$-quasigeodesic $q:I\to X$ 
such that for all subintervals $[t_1,t_2]\subset I$ 
the subpath $q|_{[t_1,t_2]}$ is contained in the tubular $D$-neighborhood of a $\Theta$-diamond 
$\diamo_{\Theta}(x_1,x_2)$
with $d(x_i,q(t_i))\leq D$.

We call a quasigeodesic {\em $\taumod$-Morse}
if it is $(L,A,\Theta,D)$-Morse for some data $(L,A,\Theta,D)$.
\end{dfn}
In particular,
$\taumod$-Morse quasigeodesics are uniformly coarsely $\taumod$-regular.

Our Morse Lemma yields the converse:
\begin{cor}[Regular implies Morse for quasigeodesics]
\label{cor:regimpmorqg}
Uniformly coarsely $\taumod$-regular quasigeodesics in model spaces
are uniform $\taumod$-Morse quasigeodesics.
\end{cor}
\proof
Let $q:I\to X$ be a $(\Theta,B)$-regular $(L,A)$-quasigeodesic.
From the conclusion of Theorem~\ref{thm:mlemt} and the $(\Theta,B)$-regularity of $q$,
it follows for any $\Theta'$ (whose interior contains $\Theta$) 
that $q$ is uniformly close also to the $\Theta'$-diamond $\diamo_{\Theta'}(x_-,x_+)$,
and analogously for the subsegments $q|_{[b_-,b_+]}$ for all subintervals $[b_-,b_+]\subset[a_-,a_+]$.
This means that $q$ is a $(L,A,\Theta',D')$-Morse quasigeodesic
for some uniform constant $D'=D'(L,A,\Theta,\Theta',B,X)$.
\qed

\medskip
Based on the notion of Morse quasigeodesic,
we defined in \cite{morse} Morse embeddings and Morse actions.
The definitions apply verbatim to all model spaces.

We first consider maps into model spaces.
Suppose that $\taumod$ and $\Theta$ are $\iota$-invariant.
\begin{dfn}[Morse embedding, cf. {\cite[Definition 7.23]{morse}}]
\label{dfn:morseemb}
A {\em $\taumod$-Morse embedding} 
from a quasigeodesic space $Z$ into $X$ is a map 
$f:Z\to X$ 
which sends uniform quasigeodesics in $Z$ 
to uniform Morse quasigeodesics in $X$. 
We call it a {\em $\Theta$-Morse embedding} 
if it sends uniform quasigeodesics 
to uniform $\Theta$-Morse quasigeodesics. 
\end{dfn}
Thus, the map is a $\taumod$-Morse embedding if for any parameters $l,a$ 
the $(l,a)$-quasi\-geo\-de\-sics in $Z$ are mapped to 
$(L,A,\Theta,D)$-Morse quasigeodesics in $X$ 
with the parameters $L,A,\Theta,D$ depending on $l,a$. 
It is a $\Theta$-Morse embedding,
if $\Theta$ is fixed and only $L,A,D$ depend on $l,a$. 

In particular,
$\taumod$-Morse embeddings 
are coarsely uniformly $\taumod$-regular quasiisometric embeddings.
(Note that they are embeddings only in a coarse sense.)
We obtain the converse:
\begin{cor}[Regular implies Morse for quasiisometric embeddings]
\label{cor:regimpmormap}
Coarsely uniformly $\taumod$-regular quasiisometric embeddings 
from quasigeodesic metric spaces into model spaces
are uniform $\taumod$-Morse embeddings.
\end{cor}
\proof
Let $f:Z\to X$ be a $(\Theta,B)$-regular $(L,A)$-quasiisometric embedding
from a quasigeodesic space $Z$. 
Then $q$ maps uniform quasigeodesics in $Z$ to 
$(\Theta,B)$-regular uniform quasigeodesics in $X$.
These are uniform $\taumod$-Morse quasigeodesics.
\qed

\medskip
Now we consider isometric group actions on model spaces.
We recall that, since $X$ has no flat factor, 
every such action becomes type preserving 
after restricting it to a suitable finite index subgroup.

We call an action Morse, if its orbit maps are Morse.
More precisely:
\begin{dfn}[Morse action, cf. {\cite[Definition 7.30]{morse}}]
\label{dfn:morseact}
We say that an isometric action $\Ga\acts X$
of a finitely generated group $\Ga$ 
is {\em $\Theta$-Morse} 
if one (any) orbit map $\Ga\to\Ga x\subset X$ 
is a $\Theta$-Morse embedding 
with respect to a(ny) word metric on $\Ga$. 
We call the action {\em $\tau_{mod}$-Morse} if it is $\Theta$-Morse 
for some $\Theta$.
\end{dfn}
Morse actions are {\em undistorted} in the sense that 
the orbit maps are quasiisometric embeddings. 
In particular, they are properly discontinuous.
Furthermore, $\Theta$-Morse actions are coarsely $\Theta$-regular.
Again, we obtain a converse:
\begin{cor}[Regular implies Morse for undistorted actions]
\label{cor:regimpmoract}
Coarsely uniformly $\taumod$-regular undistorted isometric actions
by finitely generated groups on model spaces
are uniformly $\taumod$-Morse.
\end{cor}
\proof
The orbit maps are coarsely uniformly $\taumod$-regular quasiisometric embeddings. 
\qed

\subsection{Examples of regular bilipschitz paths}\label{sec:examples}

We construct examples of regular bilipschitz paths in model spaces,
which are not close to geodesics. 
One finds such paths already in the euclidean model Weyl chamber $\De=\De_{euc}=V(0, \simod)$
and, accordingly, inside every euclidean Weyl chamber of a model space.

Let $\taumod=\simod$, and let $\Theta\subset\inte(\simod)$ be $\iota$-invariant. 
Pick a sequence of numbers $s_n\ge 1$
and a non-converging sequence of unit vectors $v_n\in \Theta$.
(Here, we identify $\simod$ with the ``unit sphere'' $V(0, \simod)\cap\D B(0,1)$.)
For instance, we can take $(v_n)$ to be an alternating sequence taking exactly two values $v_1=v_{2k-1}$ 
and $v_2=v_{2k}$, $k\in \N$, where $v_1, v_2\in\Theta$.
Now, define the 1-Lipschitz path $p: [0,+\infty)\to V(0, \simod)$ 
by concatenating the segments $x_n x_{n+1}$, where each vector $\overrightarrow{x_n x_{n+1}}$ equals $s_n v_n$. 
We claim that the path $p$ is a $\Theta$-regular quasigeodesic ray in $V(0, \simod)$. 
First of all, it follows from the definition of $p$ that it is $\Theta$-straight. Therefore, by Proposition \ref{prop:longiconvpar} the sequence $(x_n)$ is $\Theta$-longitudinal: All segments $x_m x_{n}$ for $n>m$ are $\Theta$-longitudinal. 
By the same reason (inserting additional subdivision points),
all segments $p(s)p(t)$ for $s<t$ are $\Theta$-longitudinal as well,
i.e. the path $p$ is $\Theta$-longitudinal.
In particular, $p$ is $\Theta$-regular. 
That $p$ is quasigeodesic, follows from the fact that the distance $d(\cdot,0)$ from the origin 
grows along $p$ with uniformly positive slope $\geq\eps(\Theta)>0$.

Our next goal is to ensure that $p$ is not close to a geodesic ray. In order to accomplish this, 
we choose the sequence $(s_n)$ such that
$$
\lim_{n\to+\infty} (s_{n+1}-s_n)=+\infty.
$$
Suppose that there exists a geodesic ray $r: [0,+\infty)\to V(0, \simod), r(t)=tu$
with a unit vector $u\in V(0, \simod)$, 
and a constant $C'$ such that 
$$
p([0,+\infty))\subset N_{C'}(r([0,+\infty))). 
$$
Since the sequence of lengths vectors $\overrightarrow{x_n x_{n+1}}$ diverges to infinity and 
the vectors are contained in the $C'$-neighborhood of $r([0,+\infty))$, it follows that the directions of these vectors converge to the direction vector $u$ of the ray $r$. This contradicts the assumption that the sequence of vectors $(v_n)$ does not converge.

\section{Quasiisometric embeddings of spaces and undistorted actions}
\label{sec:appl}

In this chapter 
we prove our main applications of the Morse Lemma: 
Theorems \ref{thm:main1} and \ref{thm:morse subgroups} from the introduction, 
stating the hyperbolicity of quasiisometrically embedded uniformly regular subsets of model spaces and the 
existence of a continuous extension to the Gromov boundary. 

Throughout the chapter, 
we assume that the face type $\taumod\subset\simod$ 
and the subsets $\Theta\subset\ost(\taumod)$ are $\iota$-invariant.
Then $\taumod^{\pm}=\taumod$ and $\Theta_{\pm}=\Theta$,
i.e. directions antipodal to $\taumod$-regular ($\Theta$-regular) directions 
are also $\taumod$-regular ($\Theta$-regular),
and segments satisfying one of these regularity properties keep it when reversing orientation.

\subsection{Regular subsets of euclidean buildings}

Let $X$ be a euclidean building. We now apply our results on regular paths in section~\ref{sec:rectpaths} to regular subsets.

\begin{dfn}[Rectifiably regularly path connected]
A subset $R\subset X$ is called 
{\em rectifiably $\Theta$-regularly path connected}
if any two distinct points in $R$ can be connected by a rectifiable $\Theta$-regular path contained in $R$.
\end{dfn}
Note that such subsets are in particular $\Theta$-regular.
We will now study their geometric and topological properties.

Fix a point $r\in R$.
Then the function
$$r'\mapsto\tau(rr')$$
is well-defined on $R-\{r\}$ and continuous.
(Compare the discussion of $\taumod$-directions in the beginning of section~\ref{sec:rectpaths}.)
\begin{lem}
The function $\tau(r\cdot)$ is locally constant on $R-\{r\}$.
\end{lem}
\proof
The target of the function, the set of type $\taumod$ simplices in $\Si_{r_0}X$, is a discrete space.
\qed

\begin{cor}
If $r_1,r_2\in R-\{r\}$ with $\tau(rr_1)\neq\tau(rr_2)$,
then $r_1$ and $r_2$ lie in different path components of $R-\{r\}$.
\end{cor}

The next observation relies on our main result on regular bilipschitz paths in section~\ref{sec:rectpaths}.
\begin{lem}
\label{lem:midpsependp}
Let $c:[a,b]\to R$ be a rectifiable embedded path.
Then for every $t\in(a,b)$, the points $c(a)$ and $c(b)$ 
lie in different path components of $R-\{c(t)\}$.
\end{lem}
\proof
Theorem~\ref{thm:rrpathsinflats} and Lemma~\ref{lem:opptdir} imply 
that the $\taumod$-directions $\tau_{\pm}(t)$ 
are opposite to each other.
Since $\tau(c(t)c(a))=\tau_-(t)$ and $\tau(c(t)c(b))=\tau_+(t)$,
the previous corollary yields the assertion. 
\qed

\begin{cor}
\label{cor:uniqrectconn}
(i) All embedded paths $c:[a,b]\to R$ are rectifiable.

(ii) Any two embedded paths in $R$ with the same endpoints 
agree up to reparametrization.

(iii) The image of a non-embedded path $c:[a,b]\to R$ 
contains the image of the (up to reparametrization unique) embedded path connecting its endpoints.
\end{cor}
\proof
Let $c_1,c_2:I=[a,b]\to R$ be paths with the same 
endpoints,
and suppose that $c_1$ is embedded and rectifiable.
(By assumption, any two points in $R$ are connected by a rectifiable embedded path in $R$.)
By Lemma~\ref{lem:midpsependp}, $c_2$ must go through every point on $c_1$,
i.e. $c_1(I)\subseteq c_2(I)$.

If $c_2$ is also embedded,
the lemma implies moreover that the order of the points must be preserved, 
i.e. there exists a monotonic injective map $\phi:I\to I$
such that $c_2\circ\phi=c_1$.
The image $\phi(I)\subset I$ is compact, due to the continuity of $c_2$.
If $\phi(I)\neq I$ and $(t,t')$ is a connected component of $I-\phi(I)$,
then necessarily $c_2(t)=c_2(t')$ and we arrive at a contradiction.
Therefore, $\phi$ must be bijective and hence a homeomorphism. 
This shows part (ii), and (i) follows directly.

The initial part of the proof now yields (iii).
\qed

\medskip
Let $d_R$ denote the intrinsic path metric on $R$.
\begin{cor}
\label{cor:intrtr}
$(R,d_R)$ is a metric tree. 
\end{cor}
\proof
Corollary~\ref{cor:uniqrectconn} implies that $(R,d_R)$ is a geodesic metric space,
the distance of two points given by the length of the unique embedded path connecting them. 
This path is also the {\em unique} geodesic segment in $(R,d_R)$ connecting the two points. 

It follows furthermore, that the intersection of any two geodesic segments with the same initial point
is again a geodesic segment (with this initial point),
and that geodesic triangles in $(R,d_R)$ are tripods, i.e. $(R,d_R)$ is 0-hyperbolic.
\qed

\medskip
Note furthermore, that 
the embedding $$(R,d_R)\to(X,d)$$ is $L(\Theta)$-bilipschitz,
cf.\ Corollary~\ref{cor:regpundist}.

We summarize our discussion so far:
\begin{thm}
\label{thm:regsubstr}
Rectifiably $\Theta$-regularly path connected subsets of euclidean buildings are metric trees,
when equipped with their intrinsic path metrics. 
The inclusion is a bilipschitz embedding with bilipschitz constant controlled by $\Theta$.
\end{thm}
We can say more about the extrinsic geometry of $R$ in $X$,
infinitesimally and asymptotically.

Every embedded path $c:[0,\eps)\to R$ 
has a well-defined {\em $\taumod$-initial direction}
$\tau(c)\subset\Si_{c(0)}X$ 
satisfying $$\tau(c(0)c(t))=\tau(c)$$ for all $0<t<\eps$,
cf.\ Lemma~\ref{lem:tmdir}.
\begin{add}[Antipodal infinitesimal branches]
\label{add:infgeomrsubs}
For any two embedded paths $c_1,c_2:[0,\eps)\to R$
with the same initial point $c_1(0)=c_2(0)$,
the $\taumod$-initial directions $\tau(c_1)$ and $\tau(c_2)$
are either equal or antipodal.
In the former case, there exist numbers $\eps_1,\eps_2\in(0,\eps)$ such that 
the subpaths $c_1|_{[0,\eps_1]}$ and $c_2|_{[0,\eps_2]}$ agree up to reparametrization.
\end{add}
\proof 
If the images $c_i((0,\eps))$ are not disjoint,
i.e. if there exist $\eps_1,\eps_2\in(0,\eps)$ such that $c_1(\eps_1)=c_2(\eps_2)$,
then $c_1|_{[0,\eps_1]}$ and $c_2|_{[0,\eps_2]}$ agree up to reparametrization,
cf.\ Corollary~\ref{cor:uniqrectconn}(ii).
In this case, of course, $\tau(c_1)=\tau(c_2)$. 
 
Otherwise,
if $c_1((0,\eps))\cap c_2((0,\eps))=\emptyset$,
then the concatenation 
$c= c_1\star \bar{c}_2: (-\eps, \eps)\to R$
of the path $c_1$ and the reversed path $\bar{c}_2(-t):=c_2(t)$ of $c_2$, 
is a rectifiable embedded path in $R$. 
The statement then follows from Lemma~\ref{lem:opptdir}.  
\qed 

\medskip
Every embedded path $c:[0,+\infty)\to R$ with infinite length
has a well-defined {\em $\taumod$-endpoint at infinity}
$c(+\infty)\in\DtX$
so that $$c([0,+\infty))\subset V(c(0),\st(c(+\infty))),$$
cf.\ Corollary~\ref{cor:rectlocunifregp}.
\begin{add}[Antipodal endpoints at infinity]
\label{add:asygeomrsubs}
For any two embedded paths $c_1,c_2:[0,+\infty)\to R$ with infinite length,
the $\taumod$-endpoints at infinity $c_i(+\infty)\in\DtX$
are either equal or antipodal.
In the former case, there exist numbers $t_1,t_2>0$ such that 
the subpaths $c_1|_{[t_1,+\infty)}$ and $c_2|_{[t_2,+\infty)}$ agree up to reparametrization.
\end{add}
\proof
Since $R$ is intrinsically a metric tree, cf. Corollary~\ref{cor:intrtr},
we may assume after modifying the paths,
that they have the same initial point $c_1(0)=c_2(0)$ and are otherwise disjoint. 
Then the concatenation $c= c_1\star \bar{c}_2:\R\to R$
is an embedded path in $R$,
both of whose ends have infinite length. 
It is in particular uniformly $\Theta$-regular and locally rectifiable, cf. Corollary~\ref{cor:uniqrectconn}.
The assertion then follows from Corollary~\ref{cor:rectlocunifregp}.
\qed

\subsection{Regular maps into euclidean buildings}\label{sec:Regular maps into euclidean buildings}

Let $X$ still be a euclidean building. Our discussion of regular subsets
immediately implies a restriction on the geometry of spaces 
which can be mapped into buildings by regular maps:
\begin{cor}[From tree]
\label{cor:fromtree}
If $Z$ is a path metric space 
and $Z\to X$ is a $\Theta$-regular bilipschitz map,
then $Z$ is a metric tree.  
\end{cor}
\proof
The image of the embedding is a 
rectifiably $\Theta$-regularly path connected subset of $X$ 
and hence, according to Theorem~\ref{thm:regsubstr}, a metric tree.
Thus, $Z$ is bilipschitz homeomorphic to a metric tree.
Lemma~\ref{lem:bilipttr} implies that $Z$ itself is a metric tree. 
\qed

\medskip
Consider now a $\Theta$-regular bilipschitz map
$$b:T\to X$$
from a metric tree $T$.
From our earlier discussion, 
we obtain information on the infinitesimal and asymptotic behavior.

By Addendum~\ref{add:infgeomrsubs}, we have 
in every point $t\in T$
a well-defined induced {\em infinitesimal map}
$$\Si_tb: \Si_tT\to\Si^{\taumod}_{b(t)}X:=\Flagt(\Si_{b(t)}X)$$
such that,
if $c:[0,\eps)\to T$ is a geodesic path starting in $c(0)=t$ in the direction $v\in\Si_tT$,
then the $\Theta$-regular image bilipschitz path $b\circ c$ has the $\taumod$-initial direction 
$$\tau(b\circ c)=(\Si_tb)(v).$$
Furthermore, the infinitesimal maps $\Si_tb$ are {\em antipodal},
i.e. they send distinct directions in $\Si_tT$ to opposite type $\taumod$ simplices in $\Si_{b(t)}X$.

\begin{dfn}[Antipodal map]
A map from a set into the set of simplices of a spherical building 
is called {\em antipodal} if it
sends distinct elements to antipodal simplices.
\end{dfn}

By Addendum~\ref{add:asygeomrsubs},
there is a well-defined {\em boundary map} at infinity
$$ \geo b: \geo T \to \Dt X=\Flagt(\geo X) $$
such that,
if $\rho:[0,+\infty)\to T$ is a unit speed geodesic ray in $T$,
then the $\Theta$-regular image bilipschitz ray $b\circ\rho$ in $X$ 
has the $\taumod$-endpoint at infinity
$$(b\circ\rho)(+\infty)=(\geo b)(\rho(+\infty))$$
Also $\geo b$ is {\em antipodal}.
Let 
$$\bar b:\bar T\to\widetilde{b(T)}\subset\barXt$$
denote the map from the visual compactification $\bar T=T\cup\geo T$
to the subset $\widetilde{b(T)}=b(T)\cup\DtX$ of the 
$\taumod$-bordification $\barXt= X\cup\DtX$,
which combines the map $b$ with the boundary map $\geo b$.
Since the image $b(T)\subset X$ is a $\Theta$-regular subset,
we have a well-defined topology of flag convergence on $\widetilde{b(T)}$
extending the visual topology on $\DtX$ and the subspace topology on $b(T)$,
see our discussion in section~\ref{sec:flagtop}.
It makes therefore sense to speak of the continuity of $\bar b$, 
and we can state:
\begin{thm}[Continuous extension at infinity]
\label{thm:bdregmap}
The extension $\bar b$ of $b$ is continuous with respect to the topology of flag convergence
(on $\widetilde{b(T)}$).
In particular,
the boundary map $\geo b$ is continuous with respect to the visual topology.
\end{thm}
\proof
Trees are rank one euclidean buildings
and we use the description of the visual topology on their visual compactification
as given in Fact~\ref{fact:topviscomp}.
We denote the point shadows and the corresponding basic subsets in $\geo T$ and $\bar T=T\cup\geo T$ by
$\pSh^T_{\cdot,\cdot}$ and $\barpO^T_{\cdot,\cdot}=\pO^T_{\cdot,\cdot}\cup\pSh^T_{\cdot,\cdot}$.

We must show that $\bar b$ is continuous at $\geo T$.
Consider points $t,t'\in T$.
Applying (Theorem~\ref{thm:rrpathsinflats} and) Corollary~\ref{cor:rectlocunifregp}(ii)
to geodesic rays in $T$, which start in $t$ and pass through $t'$, 
we obtain that
\begin{equation*}
\bar b(\barpO^T_{t,t'})\subset\barpOt_{b(t),b(t')} .
\end{equation*}
Let $\rho:[0,+\infty)\to T$ be a geodesic ray.
Then $b\circ\rho$ is a $\Theta$-regular bilipschitz ray in $X$,
which is contained in the Weyl cone $V(b(\rho(0)),\st((\geo b)(\rho(+\infty))))$.
The subsets 
$\barpOt_{b(\rho(0)),b(\rho(u))}\cap\widetilde{b(T)}$ for $u\to+\infty$
therefore form a neighborhood basis of $(\geo b)(\rho(+\infty))$ in $\widetilde{b(T)}$,
cf. Corollary~\ref{cor:basfltop}.
Since their $\bar b$-preimages contain the neighborhoods 
$\barpO^T_{\rho(0),\rho(u)}$ of $\rho(+\infty)$,
it follows that $\bar b$ is continuous at $\rho(+\infty)\in\geo T$.
\qed

\subsection{Regular quasiisometric embeddings into model spaces}

Let $X$ be a model space. From our results on regular maps to euclidean buildings 
we deduce now by an ultralimit argument 
corresponding results for coarsely regular maps to model spaces and isometric actions.
We first show that the large scale geometry of spaces, which admit coarsely regular maps into model spaces,
is restricted:
\begin{thm}[From hyperbolic space]
\label{thm:regQIembeddings}
If $q: Z\to X$ is a coarsely uniformly $\taumod$-regular quasiisometric embedding 
from a quasigeodesic metric space
into a model space,
then $Z$ is Gromov hyperbolic.
\end{thm}
\proof 
Since $Z$ is quasiisometric to its Rips complex $Rips_R(Z)$ for  sufficiently large $R$, 
we can assume without loss of generality
that $Z$ is a geodesic metric space. In order to verify its hyperbolicity, 
it suffices to show that every asymptotic cone of $Z$ is a metric tree, see e.g. \cite{DrutuKapovich}. 

We work with the setup as described in section~\ref{sec:ulim}.
For a sequence of scale factors $\la_n>0$ converging to zero,
a sequence of basepoints $\star_n\in Z$
and the sequence of image points $\star'_n:= q(\star_n)$ in $X$,
we consider the asymptotic cones
$$
(Z_\om, \star_\om)=\ulim (\la_n Z, \star_n), \quad (X_\om, \star'_\om)=\ulim (\la_n X, \star'_n). 
$$ 
Note that they are geodesic spaces, since the original spaces are.
By Lemma \ref{lem:ulimqiebil}, 
compare also the proof of Theorem~\ref{thm:mlemt},
the quasiisometric embedding $q$ gives rise to a uniformly $\taumod$-regular bilipschitz embedding
$$q_\om: Z_\om\to X_\om.$$
Therefore, according to Corollary \ref{cor:fromtree}, $Z_\om$ is a metric tree. 
\qed

\medskip
Now we discuss the asymptotics of coarsely regular maps from hyperbolic spaces.

Let $Z$ be a locally compact geodesic $\delta$-hyperbolic metric space,
and consider its Gromov compactification
$$\bar{Z}=Z\cup \geo Z,$$
where $\geo Z$ is the
space of equivalence classes of geodesic rays in $Z$. 
Here, two rays are called equivalent if they are asymptotic in the sense that 
their images have finite Hausdorff distance. 

The topology on $\bar{Z}$
can be described at infinity as follows,
see \cite{DrutuKapovich}.
Fix a sufficiently large number $r$, say, $r\ge 3\delta$ and define the following basic 
subsets of $\bar{Z}$: 
For points $z,w\in Z$, let the subset $\barbO_{z,w,r}\subset \bar Z$
consist of all points $\bar z\in \bar Z$, 
such that every geodesic (segment or ray)
$z\bar z$ connecting $z$ to $\bar z$ has nonempty intersection with the open ball $B(w,r)$. 
Given an ideal boundary point $\zeta\in\geo Z$ 
and a ray $\rho:[0,+\infty)\to Z$ representing it, 
$\rho(+\infty)=\zeta$, 
then the countable collection of basic subsets 
$\barbO_{\rho(0),\rho(n),r}$ for $n\in {\mathbb N}$
forms a neighborhood basis of $\xi$ in $\bar Z$.
In particular, the topology on $\bar Z$ is first-countable.

Consider now a $(\Theta,B)$-regular $(L,A)$-quasiisometric embedding 
$$q:Z\to X.$$
If $\rho:[0,+\infty)\to Z$ is a geodesic ray,
then $q\circ\rho$ is a $(\Theta,B)$-regular $(L,A)$-quasiray in $X$
and has a well-de\-fi\-ned $\taumod$-endpoint at infinity 
$(q\circ\rho)(+\infty)\in\DtX$,
cf.\ Corollary~\ref{cor:infqgclcp} and Definition~\ref{dfn:endqr}.
More precisely,
the image of $q\circ\rho$ is contained in a tubular neighborhood with uniformly controlled radius 
of the Weyl cone 
$V((q\circ\rho)(0),\st((q\circ\rho)(+\infty)))$.
In particular,
the endpoint $(q\circ\rho)(+\infty)$ depends only on the endpoint $\rho(+\infty)\in\geo Z$.
Hence $q$ induces a well-defined {\em boundary map} at infinity
$$ \geo q: \geo Z \to \Dt X=\Flagt(\geo X) $$
such that,
if $\rho:[0,+\infty)\to Z$ is a ray in $Z$,
then the $(\Theta,B)$-regular image quasiray $q\circ\rho$ in $X$ 
has the $\taumod$-endpoint at infinity
$$(q\circ\rho)(+\infty)=(\geo q)(\rho(+\infty)).$$
Furthermore, 
also as a consequence of Corollary~\ref{cor:infqgclcp},
$\geo q$ is {\em antipodal}.
Let 
$$\bar q:\bar Z\to\widetilde{q(Z)}\subset\barXt$$
denote the map from the visual compactification $\bar Z=Z\cup\geo Z$
to the subset $\widetilde{q(Z)}=q(Z)\cup\DtX$ of the 
$\taumod$-bordification $\barXt=X\cup\DtX$,
which combines the map $q$ with the boundary map $\geo q$.
Since the image $q(Z)\subset X$ is a $(\Theta,B)$-regular subset,
we have a well-defined topology of flag convergence on $\widetilde{q(Z)}$,
and we can state:
\begin{thm}[Continuous extension at infinity]
\label{thm:bdcoaregmap}
The extension $\bar q$ of $q$ is continuous at $\geo Z$ with respect to the topology of flag convergence 
on $\widetilde{q(Z)}$.
In particular, the boundary map $\geo q$ is a topological embedding with respect to the visual topology
on $\DtX$.
\end{thm}
\proof
Suppose that $\bar q$ is not continuous at the ideal point $\zeta\in\geo Z$.
Since the topology on $\bar Z$ is first-countable,
there exists a sequence $\bar z_n\to\zeta$ in $\bar Z$
such that the sequence $(\bar q(\bar z_n))$ in $\widetilde{q(Z)}$ avoids a neighborhood of 
$\geo q(\zeta)\in\DtX$.

Fix a base point $z\in Z$ 
and let $z\bar z_n$ be geodesic segments or rays in $Z$ 
connecting $z$ to $\bar z_n$.
If $\bar z_n\in\geo Z$,
then the $(\Theta,B)$-regular $(L,A)$-quasiray $q(z\bar z_n)$ 
is uniformly close to the Weyl cone $V(q(z),\st(\geo q(\bar z_n)))$,
compare the definition of the boundary map $\geo q$ above. 
If $\bar z_n\in Z$,
then $q(z\bar z_n)$ is a $(\Theta,B)$-regular $(L,A)$-quasigeodesic.
Since $\bar z_n\to\infty$ as $n\to+\infty$,
the $(\Theta,B)$-regular segment $q(z)q(\bar z_n)$ 
is $\Theta'$-regular for all sufficiently large $n$,
cf. Lemma~\ref{lem:lcoaregureg},
and $q(z\bar z_n)$ is then uniformly close to the diamond $\diamot(q(z),q(\bar z_n))$
by Theorem~\ref{thm:mlemt}.

After passing to a subsequence, 
we may assume that the $z\bar z_n$ converge to a ray $z\zeta$.
Again, the quasiray $q(z\zeta)$ 
is uniformly close to the Weyl cone $V(q(z),\st(\geo q(\zeta)))$.

Since $z\bar z_n\to z\zeta$,
there exists a sequence $w_n\to\infty$ of points $w_n\in z\bar z_n$ 
uniformly (arbitrarily) close to $z\zeta$.
Then the asymptotically $\Theta$-regular sequence $(q(w_n))$
is contained in a tubular neighborhood of $V(q(z),\st(\geo q(\zeta)))$, 
i.e. $q(w_n)\to\geo q(\zeta)$ conically.
For a sufficiently large $R>0$ 
independent of $n$, 
the balls $B(q(w_n),R)$ intersect $V(q(z),\st(\geo q(\zeta)))$
and also $V(q(z),\st(\geo q(\bar z_n)))$, respectively, $\diamot(q(z),q(\bar z_n))$.
This means that 
$\geo q(\zeta)\in\bSht_{q(z),q(w_n),R}$ 
and $\bar q(\bar z_n)\in\barbOt_{q(z),q(w_n),R}$.

Let $y_n\in B(q(w_n),R)\cap V(q(z),\st(\geo q(\zeta)))$.
Then also the sequence $(y_n)$ is asymptotically $\Theta$-regular and, 
according to Corollary~\ref{cor:basfltop},
the subsets $\barbOt_{q(z),y_n,2R}\cap\widetilde{q(Z)}$ 
form a neighborhood basis for the point $\geo q(\zeta)$ in $\widetilde{q(Z)}$.
Consequently, also the smaller neighborhoods 
$\barbOt_{q(z),q(w_n),R}\cap\widetilde{q(Z)}$ of $\geo q(\zeta)$ 
form a neighborhood basis.
Thus $\bar q(\bar z_n)\to\geo q(\zeta)$,
a contradiction. 
This shows that $\bar q$ is continuous at $\geo Z$.

That $\geo q$ is a topological embedding follows,
because it is injective (by antipodality), $\geo Z$ is compact and $\DtX$ is Hausdorff.
\qed

\medskip
We now turn to an equivariant setting and specialize the above discussion to group actions. 
We show that the class of groups, which admit asymptotically regular actions on model spaces, is restricted:

\begin{thm}[From hyperbolic group]
\label{thm:frhypgp}
If $\Ga\acts X$ is an asymptotically uniformly $\taumod$-regular undistorted isometric action 
of a finitely generated group on a model space,
then the group $\Ga$ is word hyperbolic.
\end{thm}
\proof 
Asymptotically uniformly $\taumod$-regular actions are coarsely uniformly $\taumod$-regular,
cf. Remark~\ref{rem:asyregact},
i.e. their orbit maps are coarsely uniformly $\taumod$-regular.
By assumption, they are also quasiisometric embeddings. 
The assertion therefore follows from Theorem \ref{thm:regQIembeddings}. 
\qed 

\begin{rem}
The boundary maps at infinity of the orbit maps induce a well-defined 
{\em boundary map} 
$$ \geo\Ga \to \Dt X=\Flagt(\geo X) $$
which is $\Ga$-equivariant.
\end{rem}

\section{NEW STUFF}

\begin{itemize}
\item Define Finsler metrics.
State that diamonds are unions of Finsler geodesics.
\item Finsler formulation of Morse Lemma:
Every (usual) quasigeodesic is uniformly close to a Finsler geodesic.
Note: The Finsler metrics can be pseudo-metrics, and then they are not equivalent to the Riemannian metric.
\item Define Finsler quasi-convex subset and subgroup:
Any two points are connected by {\em some} Finsler geodesic which is uniformly close to the subset.
Of course, ``some'' cannot be replaced by ``all''.
Note that they are not Riemannian/CAT(0) quasiconvex.
\item Applications: Local-to-global principle for coarsely regular quasigeodsics.
List of equivalent definitions, adding to the list from \cite{morse} the items URU and Finsler quasi-convex.
\end{itemize}

Addresses:

\noindent M.K.: Department of Mathematics, \\
University of California, Davis\\
CA 95616, USA\\
email: kapovich@math.ucdavis.edu

\noindent B.L.: Mathematisches Institut\\
Universit\"at M\"unchen \\
Theresienstr. 39\\ 
D-80333, M\"unchen, Germany\\ 
email: b.l@lmu.de

\noindent J.P.: Departament de Matem\`atiques,\\
 Universitat Aut\`onoma de Barcelona,\\ 
 08193 Bellaterra, Spain\\
email: porti@mat.uab.cat

\end{document}